\documentclass[11pt]{article}

\usepackage{graphicx}
\usepackage{natbib}
\usepackage{hyperref}
\usepackage{url}

\usepackage{amstext}
\usepackage{amssymb}
\usepackage{amsmath}
\usepackage{bbm}
\usepackage{booktabs}
\usepackage{caption}
\usepackage{subcaption}
\usepackage{multirow}

\usepackage{comment}

\usepackage{algorithm}
\usepackage{algorithmic}

\DeclareMathOperator*{\argmax}{arg\,max}
\DeclareMathOperator*{\argmin}{arg\,min}

\newcommand{\eg}{\textit{e.g.}}
\newcommand{\ie}{\textit{i.e.}}

\newcommand{\Halmos}{$\hfill \square$}

\newtheorem{lemma}{Lemma}
\newtheorem{theorem}{Theorem}

\newtheorem{proposition}{Proposition}
\newtheorem{corollary}{Corollary}

\newcommand{\boxxx}[1]
 {\begin{center}\fbox{\begin{minipage}{13cm}#1\smallskip\end{minipage}}\end{center}}

\usepackage{tikz}
\usepackage{pgfplots}
\pgfplotsset{compat=1.18}
\usepgfplotslibrary{statistics}
\usetikzlibrary{positioning}

\allowdisplaybreaks

\title{{\bf Stackelberg Dynamic Location Planning under Cumulative Demand}}

\author{}

\author{{\bf Warley Almeida Silva}\\
CIRRELT, Université de Montréal\\
warley.almeida.silva@umontreal.ca
\and
{\bf Margarida Carvalho} \\
    CIRRELT, Université de Montréal \\
carvalho@iro.umontreal.ca\and
{\bf Sanjay Dominik Jena}\\
CIRRELT, Université du Quebec à Montréal\\
jena.sanjay-dominik@uqam.ca}

\date{}

\newcommand{\mathset}[1]{\{{#1}\}}

\newcommand{\ourproblem}[1]{CDFLP-CCD} 

\newcommand{\ceil}[1]{\left\lceil #1  \right\rceil}

\begin{document}
\bibliographystyle{apalike}

\pagenumbering{gobble}

\maketitle

\begin{abstract}
    Dynamic facility location problems predominantly suppose a monopoly over the service or product provided.
    Nonetheless, this premise can be a severe oversimplification in the presence of market competitors, as customers may prefer facilities installed by one of them.
    The monopolistic assumption can particularly worsen planning performance when demand depends on prior location decisions of the market participants, namely, when unmet demand from one period carries over to the next. 
    Such a demand behaviour creates an intrinsic relationship between customer demand and location decisions of  all market participants, and requires the decision-maker to anticipate the competitor's response.
    This work studies a novel competitive facility location problem that combines cumulative demand and market competition to devise high-quality solutions.
    We propose bilevel mixed-integer programming formulations for two variants of our problem, prove that the optimistic variant is $\Sigma^{p}_{2}$-hard, and develop branch-and-cut algorithms with tightened value-function cuts that significantly outperform general-purpose bilevel solvers.
    Our results quantify the severe cost of planning under a monopolistic assumption (profit drops by half on average) and the gains from cooperation over competition (6\% more joint profit), while drawing managerial guidelines on how instance attributes and duopolistic modelling choices shape robust location schedules. 
    
    \,

    \textbf{Keywords:} Facility Location, Stackelberg Competition, Cumulative Demand, Bilevel Programming.

\end{abstract}

\newpage

\section{Introduction}

Dynamic Facility Location Problems (FLPs) are at the heart of strategic production and supply chain planning \citep{nickel2019multi}.
Location decisions typically involve costly site-specific long-term commitments (\eg, to acquire properties and build infrastructure), making relocation or closure during the planning horizon undesirable.
Throughout the last decade, several economic shifts have contributed to the emergence of temporary facilities \citep[see, \eg,][]{rudkowskiHereTodayGone2020, rosenbaumBenefitsPitfallsContemporary2021, CaoQiZhang2024OnlineFacilityLocation}, which remain open at some locations for a short period of time (\eg, weeks or months) before moving to other locations.
For example, in retail, temporary pop-up stores have been integrated into operations for many years and are expected to exceed 95 billion dollars in revenue in 2025 \citep{capital2025}.
By serving different regions at a time, companies may create a sense of urgency among customers currently covered by a facility, while leading to an increase of interest (hype) and demand among those left uncovered. 
%

%
%
%

From a mathematical optimization point of view, the planning for temporary facility location exhibits several unique characteristics. 
Classical mathematical formulations for dynamic FLPs assume that customer demand is independently defined for each time period of the planning horizon \citep[see, \eg,][]{jenaDynamicFacilityLocation2015, marinMultiperiodStochasticCovering2018, alizadehMultiPeriodMaximalCovering2021, vatsaCapacitatedMultiperiodMaximal2021}.
However, for several commodities, this may severely misrepresent the real underlying demand behaviour. 
Indeed, unmet demand for non-perishable items such as personal clothing, beauty products, tech gadgets, and handmade crafts may not simply vanish, but is likely to persist in the future. 
Temporary facilities may therefore cause a more intrinsic relationship between location decisions and customer demand \citep{silvaDynamicFacilityLocation2025}, as demand may carry over to future time periods and build up while facilities are not accessible.
%
%
Such demand behaviour has recently been addressed in the literature, referred to as Cumulative Customer Demand (CCD) \citep[see, \eg,][]{daneshvarTwostageStochasticPostdisaster2023,silvaDynamicFacilityLocation2025}, often exhibiting drastically different structures of optimal solutions that favour the relocation of facilities over the planning horizon.

Indeed, when companies operate in a monopoly, they may intentionally decide to let customer demand build up over time while operating in different regions, maximizing the total profit extracted over the planning horizon \citep{silvaDynamicFacilityLocation2025}.
However, such a strategy is considerably vulnerable when competing firms are present, as the competition may capture the accumulated customer demand before the company.
Explicitly accounting for competition within the planning stage is therefore of crucial importance to ensure profitability in such dynamic and competitive markets.
The literature on competitive FLPs \citep[see, \eg][]{ljubicOuterApproximationSubmodular2018,beresnevApproximationCompetitiveFacility2019,qiSequentialCompetitiveFacility2022}  has modelled market competition mostly as a Stackelberg game \citep{von1934marktform}. 
However, to the best of our knowledge, existing works in this research stream are restricted to a single time period, and therefore cannot readily account for cumulative customer demand.

In this work, we investigate how to handle market competition under cumulative demand.
We consider a company, further referred to as the leader, that must place \textit{her} temporary facilities over the planning horizon to capture customers and maximize her profit.
The leader is aware of a competitor, further referred to as the follower, that will react by locating \textit{his} temporary facilities over the planning horizon to capture customers and maximize his profit.
Customers patronize facilities installed by the leader or the follower at each period based on their preferences, and unmet demand carries over to subsequent periods.
This planning problem, referred to as the \textit{Competitive Dynamic Facility Location Problem under Cumulative Customer Demand} (\ourproblem{}), combines cumulative demand and market competition to produce high-quality location schedules for the leader.
As such, we contribute to the literature on location problems as follows:
\begin{enumerate}
    \item We introduce the \ourproblem{}, a new temporary facility location problem that requires to merge distinctive features which previously have been considered separately. We examine both the optimistic variant of the respective Stackelberh competition, which is more commonly studied in the literature, and the considerably less explored pessimistic variant of the \ourproblem{}.
    \item We explicitly prove the $\Sigma^{p}_{2}$-hardness of the optimistic variant, even if each player is limited to locating a single facility over the planning horizon. To the best of our knowledge, we present the first $\Sigma^{p}_{2}$-hardness result for a non min-max, non zero-sum competitive location problem (\ie, each player maximizes its profit, and the loss of a player does not imply the gain of another). We thus enrich the catalogue of known $\Sigma^{p}_{2}$-complete decision problems, offering a new starting point for complexity reductions in non min-max, non zero-sum bilevel settings.
    \item We propose bilevel mixed-integer programming formulations for two variants of the \ourproblem{}, solvable through branch-and-cut frameworks, and devise a tailored version of the classical value-function cut \citep[see, \eg,][]{lozanoValueFunctionBasedExactApproach2017}. We then present a novel technique to render this tailored value-function cut tight, and prove its theoretical dominance.  On average, the tightened value-function cut solves instances three times faster than the tailored one, and proves optimality gaps half as large, within the same time limit.
    Our tightening technique can be adapted to other competitive problems, including those that are non mix-max, non zero-sum.
    \item We draw managerial insights for our benchmark, inspired by the Canadian province of Quebec. We show that, on average, the leader would obtain half the profit she could have obtained by wrongfully assuming a monopoly at the planning stage, and that players could extract $6\%$ more profit from customers if they cooperated. We also examine the structure of optimal solutions, underlining how each instance attribute shapes the optimal location schedule, and assess the quality of optimal solutions depending on the leader's assumptions about follower's and customers' behaviours. 
\end{enumerate}

The remainder of the paper is organized as follows. 
We review the related literature and position the contributions of our work in Section~\ref{sec:literature}.
Section~\ref{sec:formulation} then defines the optimistic and a pessimistic variants of the \ourproblem{}, and formulates them as bilevel mixed-integer programming formulations.
In Section~\ref{sec:methods}, we characterize the theoretical complexity of the optimistic variant of the \ourproblem{} and present exact methods to tackle each variant.
We finally evaluate the computational performance of our solution methods and devise managerial insights from our computational benchmark in Section~\ref{sec:experiments}, and conclude with final remarks and future work in Section~\ref{sec:conclusion}.

\section{Literature Review}
\label{sec:literature}

%
To the best of our knowledge, dynamic FLPs in competitive settings have not yet been tackled in the literature.
We therefore focus here on (single-period) competitive FLPs, which have received a lot of attention in the last decade \citep{mishraLocationCompetitiveFacilities2022}, and are computationally significantly more complex than classical FLPs without competition.
Most related works \citep[see, \eg][]{hemmatiMixedintegerBilevelProgramming2016, beresnevApproximationCompetitiveFacility2019, qiSequentialCompetitiveFacility2022} assume a sequential decision-making process, where the company (\ie, the leader) decides where to place facilities in anticipation of an optimal reaction of the competitor (\ie, the follower).
This setting leads to Stackelberg games \citep{von1934marktform} that are naturally modelled as bilevel mixed-integer programming formulations \citep{kleinertSurveyMixedIntegerProgramming2021}.
Existing works differ from each other on key assumptions about how players and customers behave in the duopoly, which we discuss in Section~\ref{sec:single-period}.
We then review solution methods for bilevel mixed-integer programming formulations in Section~\ref{sec:bilevel-programming} and summarize how our contributions bridge existing gaps in the literature on location problems in Section~\ref{sec:gap-contributions}.

\subsection{Core Assumptions in Competitive FLPs}
\label{sec:single-period}

We first focus on assumptions about players and their behaviour.
With respect to the players' objectives, in zero-sum games \citep[see, \eg, ][]{noltemeierMultipleVotingLocation2007,kucukaydinLeaderFollowerGame2012, alekseevaExactMethodDiscrete2015, qiSequentialCompetitiveFacility2022},  maximizing the profit of the follower requires minimizing the profit of the leader.
We may therefore derive min-max formulations, where players have strictly opposing objective functions.
Non zero-sum games \citep[see, \eg, ][]{hemmatiMixedintegerBilevelProgramming2016,gentileIntegerProgrammingFormulations2018,beresnevApproximationCompetitiveFacility2019, linLocatingFacilitiesCompetition2022} are more complex to model given that different optimal solutions for the follower may yield considerably different outcomes for the leader.
In this context, we may consider either an optimistic follower behaviour, 
in which he breaks ties in favour of the leader \citep[see, \eg,][]{gentileIntegerProgrammingFormulations2018}, or a pessimistic follower behaviour, 
in which he breaks ties against the leader \citep[see, \eg,][]{beresnevApproximationCompetitiveFacility2019}.
While the Stackelberg structure of competitive FLPs is generally indicative of the problem's computational complexity, we are only aware of one $\Sigma^{p}_{2}$-hardness result for the \textit{zero-sum, min-max} $(r|p)$-centroid problem \citep{noltemeierMultipleVotingLocation2007}.
It is well-known that $\Sigma^{p}_{2}$-hardness results are scarce \citep{grune2025completeness} and, to the best of our knowledge, there are none for \textit{non min-max, non zero-sum} competitive location problems like the \ourproblem{}.
With respect to the players' strategy sets, the vast majority of competitive FLPs assumes that they have access to disjoint sets of locations \citep[see, \eg,][]{plastriaDiscreteModelsCompetitive2008, kucukaydinLeaderFollowerGame2012,dreznerLeaderFollowerModel2015,qiCompetitiveFacilityLocation2017,gentileIntegerProgrammingFormulations2018,linLocatingFacilitiesCompetition2022} or that co-location is explicitly forbidden, akin to interdiction games \citep[see, \eg,][]{roboredoBranchandcutAlgorithmDiscrete2013,beresnevBranchandboundAlgorithmCompetitive2013, alekseevaExactMethodDiscrete2015, beresnevApproximationCompetitiveFacility2019, qiSequentialCompetitiveFacility2022}.
To the best of our knowledge, \cite{hemmatiMixedintegerBilevelProgramming2016} are the only ones to address co-location by supposing that customers patronizing a location occupied by both players split their demand proportionally to their perception of each player (\eg, brand recognition).

We now turn to assumptions on preference behaviour of customers, which may have to decide at which facility to satisfy demand when more than one available facility is deemed appropriate.
Implemented as Discrete Choice Models, the literature models customer behaviour either as deterministic, choosing exactly one facility among those made available \citep[see, \eg,][]{
roboredoBranchandcutAlgorithmDiscrete2013,hemmatiMixedintegerBilevelProgramming2016, gentileIntegerProgrammingFormulations2018, beresnevApproximationCompetitiveFacility2019}, or as probabilistic, patronizing each available facility with a non-zero probability \citep[see, \eg,][]{dreznerLeaderFollowerModel2015,qiCompetitiveFacilityLocation2017,  linLocatingFacilitiesCompetition2022, qiSequentialCompetitiveFacility2022}.
%
%
While classical probabilistic choice models, such as the Multinomial Logit Model and its extensions, have been integrated into FLPs without and with competition, such models induce nonlinear terms in the objective function and have been found to have limited predictive performance \citep{berbegliaComparativeEmpiricalStudy2022} due to their fixed model size. 
Non-parametric choice models, i.e., those that can adjust its complexity to the amount and variety of available data, such as the rank-based choice model  \citep[see, \eg,][]{fariasNonparametricApproachModeling2013, vanryzinMarketDiscoveryAlgorithm2015, jenaPartiallyRankedChoice2020}, have shown superior predictive performance \citep{berbegliaComparativeEmpiricalStudy2022} and can be integrated into optimization models using linear terms. 
While rank-based choice models consist of a probability distribution of strictly ranked preference lists, 
it has recently been shown that they can be represented through a single deterministic ranked list for each customer by duplicating customers for each probabilistic preference list \citep[see ][Appendix A]{silvaDynamicFacilityLocation2025}.
%
They thus allow for tractable formulations, while preserving predictive accuracy and a realistic choice behaviour. 
Finally, mathematical formulations embedding customer rankings can be easily adjusted to handle co-location \citep{hemmatiMixedintegerBilevelProgramming2016}.
Such a feature is essential in the \ourproblem{}, since companies may effectively deploy their temporary facilities at the same location (\eg, the same park or shopping mall).
%

%
%
%
%

\subsection{Bilevel Programming Methodologies}
\label{sec:bilevel-programming}

Problem-specific algorithms \citep[see, \eg,][]{roboredoBranchandcutAlgorithmDiscrete2013, hemmatiMixedintegerBilevelProgramming2016,beresnevApproximationCompetitiveFacility2019,qiSequentialCompetitiveFacility2022} and general-purpose bilevel solvers \citep[see, \eg,][]{mooreMixedIntegerLinear1990, denegre2009branch,xuExactAlgorithmBilevel2014,wangWatermelonAlgorithmBilevel2017, fischettiNewGeneralPurposeAlgorithm2017, lozanoValueFunctionBasedExactApproach2017, tahernejadBranchandcutAlgorithmMixed2020} 
employ branch-and-cut frameworks to solve bilevel mixed-integer programming formulations.
Specifically, the works cited above relax the optimality of the lower level (\ie, the follower) and progressively add cuts to remove solutions that are not bilevel feasible (\ie, where the reaction of the follower is suboptimal).
Each work then exploits core premises about players and customers to further refine their solution methods.
This is the case for the well-known interdiction games \citep[see, \eg,][]{fischettiInterdictionGamesMonotonicity2019, taninmics2022branch}, which are zero-sum, min-max and where (the equivalent of) co-location is forbidden.

Overall, the literature on competitive FLPs \citep[see, \eg,][]{roboredoBranchandcutAlgorithmDiscrete2013, hemmatiMixedintegerBilevelProgramming2016,beresnevApproximationCompetitiveFacility2019,qiSequentialCompetitiveFacility2022} favours algorithms tailored to specific problems as they scale better than general-purpose bilevel solvers. 
Problem-specific algorithms predominantly devise tailored variants of general value-function cuts \citep[see, \eg,][]{lozanoValueFunctionBasedExactApproach2017}, which impose globally valid lower (respectively, upper) bounds to the minimization (respectively, maximization) problem solved by the follower based on problem-specific assumptions.

\subsection{Relation to Our Work}
\label{sec:gap-contributions}

We conclude by positioning our contributions within the literature. On the modelling side, we bridge the disconnect between dynamic and competitive FLPs by integrating cumulative customer demand ~\citep{silvaDynamicFacilityLocation2025} into a Stackelberg competition. Unlike the majority of the literature on bilevel programs, we address a general non min-max setting. This creates a distinct challenge regarding the follower decision selection which is absent in min-max games. Furthermore, we generalize the spatial restrictions found in typical competitive FLPs. Rather than strictly forbidding co-location (a common assumption in interdiction), we adopt a flexible approach akin to \cite{hemmatiMixedintegerBilevelProgramming2016}, where co-location is permissible and governed by demand-splitting parameters.

From a theoretical perspective, we address the scarcity of complexity results for non-zero-sum games. While the hardness of min-max interdiction games problems is well-documented \citep[see, \eg,][]{furini2019maximum,frohlich2021hardness,boggio2025completeness}, we provide the first $\Sigma^{p}_{2}$-completeness result for a competitive location problem where objectives are not strictly opposing. Finally, regarding solution methodology, standard value-function cuts cannot be directly applied due to the cumulative nature of the demand. Consequently, we depart from general-purpose bilevel solvers, notably the MIX++ Solver \citep{fischettiNewGeneralPurposeAlgorithm2017} which we show in Appendix~\ref{apx:preliminary} to be computationally prohibitive for our benchmark, and instead develop a specialized branch-and-cut algorithm. This approach extends the intuition of previous value-function cuts but requires a novel tightening strategy to remain effective under the proposed dynamic structure.

\section{Mathematical Formulation} 
\label{sec:formulation}

We define the \ourproblem{} in Section~\ref{sec:definition} and present stylized mathematical formulations for both optimistic and pessimistic variants in Section~\ref{sec:stylized}.
We then discuss each component of these stylized mathematical formulations in Section~\ref{sec:components} and devise detailed mathematical formulations in Section~\ref{sec:detailed}.
In the remainder of this paper, bold letters are used to denote vectors.

\subsection{Problem Definition}
\label{sec:definition}

Let $\mathcal{I}$ be the set of candidate locations, $\mathcal{J}$ be the set of targeted customers, and $\mathcal{T} = \{1, \ldots, T\}$ be the set of time periods.
First, the leader decides where to locate her facilities over the planning horizon.
Let $y^{t}_{i} \in \mathset{0,1}$ take value $1$ if the leader opens a facility at location $i$ during time period $t$, and $0$ otherwise.
The leader then announces location schedule $\boldsymbol{y}$ to the public. 
In reaction, the follower decides where to place his facilities over the planning horizon.
To this end, let $z^{t}_{i} \in \mathset{0,1}$ equal $1$ if the follower opens a facility at location $i$ during time period $t$, and $0$ otherwise.
Once both location schedules $\boldsymbol{y}$ and $\boldsymbol{z}$ are made public, customers decide which facility they will attend at each period of the planning horizon.
Each customer satisfies their entire demand at the facility location that they rank highest in their preference model. If none of their preferred locations holds a facility, their demand remains unserved and carries over to the next time, accumulating until satisfied by one of the players.
The player that satisfies the accumulated customer demand collects the corresponding profit. 
The leader wants to maximize her profit $\pi^{L}(\boldsymbol{y}, \boldsymbol{z})$, which explicitly depends on the location schedule $\boldsymbol{z}$ of the follower, as customers may prefer to patronize facilities installed by the follower. 
Similarly, the follower desires to maximize his profit $\pi^{F}(\boldsymbol{y}, \boldsymbol{z})$.
%






\subsection{Stylized Formulations}
\label{sec:stylized}

Let $\mathcal{Y}$ (respectively, $\mathcal{Z}$) be the set of feasible location schedules for the leader (respectively, the follower), which includes  previously defined binary constraints for variables $y^{t}_{i}$ and $z^{t}_{i}$.
We write the stylized \textit{optimistic} formulation of the \ourproblem{} as follows:
\begin{subequations}
    \label{mdl:stylized-optimistic-formulation}
    \begin{align}
        \mathcal{L}^{optimistic}:  \quad 
        \max_{\boldsymbol{y}, \boldsymbol{z}} \quad \left\{
          \pi^{L} \left(\boldsymbol{y}, \boldsymbol{z}\right)
            \, \, \middle| \, \,
         \boldsymbol{y} \in \mathcal{Y}, \boldsymbol{z} \in \mathcal{F}(\boldsymbol{y}) \right\} \nonumber 
        \text{ where } \mathcal{F}(\boldsymbol{y}) =  \argmax_{\boldsymbol{\hat{z}} \in \mathcal{Z}} \quad
         \pi^{F} \left(\boldsymbol{y}, \boldsymbol{\hat{z}}\right). 
          \nonumber 
    \end{align}
\end{subequations}

In the presence of multiple optimal solutions for the follower, the optimistic formulation assumes that the follower will break ties in favour of the leader.
One may argue that a pessimistic follower behaviour is more realistic, where the follower breaks ties in detriment of the leader.
In this context, we can write the stylized \textit{pessimistic} formulation of the \ourproblem{} as follows:
\begin{subequations}
    \label{mdl:stylized-pessimistic-formulation}
    \begin{align}
        \mathcal{L}^{pessimistic}:  
        \max_{\boldsymbol{y}} \min_{\boldsymbol{z}}  \quad
        \left\{ \pi^{L} \left(\boldsymbol{y}, \boldsymbol{z}\right)
       \, \, \middle| \, \, 
        \boldsymbol{y} \in \mathcal{Y}, \boldsymbol{z} \in \mathcal{F}(\boldsymbol{y}) \right\} \nonumber 
        \text{ where } \mathcal{F}(\boldsymbol{y}) = \argmax_{\boldsymbol{\hat{z}} \in \mathcal{Z}} \quad
         \pi^{F} \left(\boldsymbol{y}, \boldsymbol{\hat{z}} \right). 
        \nonumber
        \end{align}
\end{subequations}

In the remainder of the paper, we mainly focus on the optimistic variant of the \ourproblem{}.
Nevertheless, we explain how to adapt our solutions methods for the pessimistic variant in Section~\ref{sec:pessimistic}, briefly discuss their computational performance in Section~\ref{sec:performance}, and consider pessimistic bilevel solutions to devise managerial insights about the \ourproblem{} in Section~\ref{sec:insights}.

\subsection{Problem Components}
\label{sec:components}

We now explain each component of the stylized formulations in more detail.

\subsubsection{Location Decisions}
\label{sub:locations}
Feasible sets $\mathcal{Y}$ and $\mathcal{Z}$ may contain a wide range of constraints on locations decisions.
For example, they may include phase-in (respectively, phase-out) constraints that require players to only open (respectively, only close) facilities over the planning horizon, or budget constraints that restrict the number of facilities at each period.
In this paper, we allow each player to locate a single facility at each period regardless of where the facility was open at the previous period.
Specifically, we consider
\begin{align*}
    \mathcal{Y} = \left\{y^{t}_{i} \in \mathset{0,1}^{I \times T} \, \, \middle| \, \,\sum_{i \in \mathcal{I}} y^{t}_{i} \leq 1, \forall t \in \mathcal{T}\right\} \textrm{ and }
    \mathcal{Z} = \left\{z^{t}_{i} \in \mathset{0,1}^{I \times T}  \, \, \middle| \, \, \sum_{i \in \mathcal{I}} z^{t}_{i} \leq 1, \forall t \in \mathcal{T}\right\}.
\end{align*}
Such feasible sets allow us to focus on how competition influences location decisions, particularly under cumulative customer demand, rather than analyzing the trade-off with problem-specific constraints on location decisions.
We highlight, however, that our solution methods can be applied to different feasible sets $\mathcal{Y}$ and $\mathcal{Z}$, as long as they can be represented through linear constraints, thus including applications where players have more than a single facility available at each period.

\subsubsection{Profit Functions}

Profit functions $\pi^{L} (\boldsymbol{y}, \boldsymbol{z})$ and $\pi^{F} (\boldsymbol{y}, \boldsymbol{z})$ depend on (\textit{i}) how customers patronize facilities throughout the planning horizon, (\textit{ii}) how much demand has accumulated since their last capture, \textit{(iii)} whether their profit will be split among players, and \textit{(iv)} other logistical costs.

\paragraph{Customer Preferences.}
We employ rankings to represent customer behaviour \citep[see, \eg,][]{
roboredoBranchandcutAlgorithmDiscrete2013,hemmatiMixedintegerBilevelProgramming2016, gentileIntegerProgrammingFormulations2018, beresnevApproximationCompetitiveFacility2019}, which allow the representation of rank-based choice models \citep[see, \eg,][]{fariasNonparametricApproachModeling2013, vanryzinMarketDiscoveryAlgorithm2015, jenaPartiallyRankedChoice2020}.
Each customer $j$ patronizes only a subset of locations called the consideration set.
Facilities in this set are ranked by preference, and customers patronize the most preferred location in their consideration set.
More specifically, we denote the ranking of customer $j$ over the set of candidate locations $\mathcal{I}$ as $\succ_{j}$, where $i \succ_{j} k$ indicates that customer $j$ prefers location $i$ over location $k$.
We employ an artificial location $0$ to allow customers to rank the choice of no service at all (\ie, $0 \succ_{j} i$ stipulates that customer $j$ prefers no service over location $i$).
We assume that customers always look for service unless there is no facility available in their consideration set.

\paragraph{Demand Accumulation.}
We adapt the cumulative demand behaviour of \cite{silvaDynamicFacilityLocation2025} to the competitive setting. 
Let $d^{t}_{j} \in \mathbb{R}^{+}$ be the spawning demand of customer $j$ at period $t \geq 1$.
We formally define the accumulated demand of customer $j$ at the beginning of period $t$ as $c^{t}_{j} (\boldsymbol{y}, \boldsymbol{z}) = u^{t-1}_{j} (\boldsymbol{y}, \boldsymbol{z}) + d^{t}_{j}$, and the unmet demand $u^{t}_{j} (\boldsymbol{y}, \boldsymbol{z})$ of customer $j$ at period $t$ as
    $$u^{t}_{j} (\boldsymbol{y}, \boldsymbol{z}) = \begin{cases}
        \left(1 - \max_{\substack{i \in \mathcal{I}: \\ i \succ_{j} 0}} \mathset{y^{t}_{i}, z^{t}_{i}} \right) c^{t}_{j} (\boldsymbol{y}, \boldsymbol{z}), & \text{ if } t \in \mathcal{T}, \\
        0, & \text{ if } t = 0.
    \end{cases}$$
Note that the term $\max_{\substack{i \in \mathcal{I}: \\ i \succ_{j} 0}} \mathset{y^{t}_{i}, z^{t}_{i}}$ equals $1$ if and only if customer $j$ has been captured by the leader or the follower through some location $i$ at period $t$.
Cumulative customer demand introduces an additional layer of complexity to the leader's planning problem, as the customer demand available at period $t$ is not only affected by the leader's location decisions from periods $1,2, \ldots, t-1$, but also the respective location decisions of the follower, over which the leader has no control.

\paragraph{Profit Splitting.}
We employ the profit splitting introduced by \cite{hemmatiMixedintegerBilevelProgramming2016}.
At period $t$, customer $j$ may be (\textit{i}) captured only by the leader, (\textit{ii}) captured only by the follower, (\textit{iii}) captured by both players, or (\textit{iv}) not captured at all.
Let $r_{i} \in \mathbb{R}^{+}$ be the reward per demand unit captured through location $i$.
If customer $j$ is captured only by the leader (respectively, the follower) through location $i$, she (respectively, he)  receives the \textit{integral} contribution $r_{i} c^{t}_{j} (\boldsymbol{y}, \boldsymbol{z})$.
Both players can simultaneously capture customer $j$ at period $t$ only if they do it through the same most preferred location $i$.
In this case, the leader and the follower receive \textit{partial} contributions $\rho r_{i}  c^{t}_{j} (\boldsymbol{y}, \boldsymbol{z})$ and $(1 - \rho) r_{i} c^{t}_{j} (\boldsymbol{y}, \boldsymbol{z})$, respectively, where $\rho \in [0,1]$ is the splitting factor.
This splitting factor may be set, for example, based on brand recognition of the leader and the follower among customers.
If customer $j$ is not captured at all, there is no contribution to account for in the profit of each player.

\paragraph{Logistics Costs.}
Location schedules $\boldsymbol{y}$ and $\boldsymbol{z}$ may also incur logistics costs $\gamma^{L} (\boldsymbol{y})$ and $\gamma^{F} (\boldsymbol{z})$, respectively, that directly influence the profit of each player.
These logistics costs are related to assembling, disassembling, and transportation costs, among others.
In this paper, we assume that logistics costs are marginal in comparison to revenues, and therefore set cost functions $\gamma^{L} (\boldsymbol{y}) = 0, \forall \boldsymbol{y} \in \mathcal{Y}$ and $\gamma^{F} (\boldsymbol{z}) = 0, \forall \boldsymbol{z} \in \mathcal{Z}$.
We highlight, however, that our solution methods can be applied to different cost functions $\gamma^{L} (\boldsymbol{y})$ and $\gamma^{F} (\boldsymbol{z})$, as long as they can be expressed linearly.

\subsection{Detailed Formulations}
\label{sec:detailed}

We now present the detailed \textit{optimistic} formulation of the \ourproblem{}, where we write components $\pi^{L} (\boldsymbol{y}, \boldsymbol{z})$ and $\pi^{F} (\boldsymbol{y}, \boldsymbol{z})$ linearly.
Let $u^{\ell{}t}_{ij} \in [0,1]$ (respectively, $v^{\ell{}t}_{ij} \in [0,1]$) be the percentage of demand from customer $j$ accumulated from period $\ell{} + 1$ to period $t$ and captured through location $i$ installed by the leader (respectively, the follower).
These decision variables represent the accumulated demand function $c^{t}_{j} (\boldsymbol{y}, \boldsymbol{z})$. 
If customer $j$ was last captured at period $\ell{}$, we can compute its accumulated demand at period $t$ as $D^{\ell{}t}_{j} = \sum_{\substack{s \in \mathcal{T}: \\ s > \ell{} \\ s \leq t}} d^{s}_{j}$.
Let $\mathcal{T}^{S} = \{0\} \cup \mathcal{T}$ and $\mathcal{T}^{F} = \mathcal{T} \cup \{T+1\}$ be the set of periods, additionally including either the start period $0$ and the final period $T + 1$, respectively.
We write the detailed optimistic formulation of the \ourproblem{} as follows:
\begin{subequations}
    \label{mdl:detailed-optimistic-formulation}
    \begin{align}
        \mathcal{L}^{optimistic}:  
        \max_{\boldsymbol{y}, \boldsymbol{z}, \boldsymbol{u}, \boldsymbol{v}} \quad
        &  \sum_{t \in \mathcal{T}} \sum_{\substack{\ell{} \in \mathcal{T}^{S}: \\ \ell{} < t}}\sum_{j \in \mathcal{J}} \sum_{\substack{i \in \mathcal{I}: \\ i \succ_{j} 0}} r_{i} D^{\ell{}t}_{j}  u^{\ell{}t}_{ij}
        && \label{eq:upper-detailed-optimistic-obj}\\
        \text{s.t.} \quad
        & \boldsymbol{y} \in \mathcal{Y}
        &&  \label{eq:upper-detailed-optimistic-ct1} \\
        %
        %
        & (\boldsymbol{z}, \boldsymbol{u}, \boldsymbol{v}) \in \mathcal{F}(\boldsymbol{y})
        &&  \label{eq:upper-detailed-optimistic-dm2} \\
        \text{ where } \mathcal{F}(\boldsymbol{y}) = 
        %
        %
        \argmax_{\boldsymbol{\hat{z}}, \boldsymbol{\hat{u}}, \boldsymbol{\hat{v}}} \quad
        &  \sum_{t \in \mathcal{T}} \sum_{\substack{\ell{} \in \mathcal{T}^{S}: \\ \ell{} < t}} \sum_{j \in \mathcal{J}} \sum_{\substack{i \in \mathcal{I}: \\ i \succ_{j} 0}} r_{i} D^{\ell{}t}_{j}  \hat{v}^{\ell{}t}_{ij}
        && \label{eq:lower-detailed-optimistic-obj}\\
        \text{s.t.} \quad
        %
        & \boldsymbol{\hat{z}} \in \mathcal{Z}
        &&  \label{eq:lower-detailed-optimistic-ct1} \\
        & \sum_{\substack{\ell{} \in \mathcal{T}: \\ \ell{} < t}} \hat{u}^{\ell{} t}_{ij} \leq y^{t}_{i} - \left(1 - \rho\right) y^{t}_{i} \hat{z}^{t}_{i}
        && \substack{\forall i \in \mathcal{I}, \\ \forall j \in \mathcal{J}, \\ \forall t \in \mathcal{T} : \\ i \succ_{j} 0} \label{eq:third-detailed-optimistic-ct1} \\
        & \sum_{\substack{\ell{} \in \mathcal{T}: \\ \ell{} < t}} \hat{v}^{\ell{} t}_{ij} \leq \hat{z}^{t}_{i} - \rho y^{t}_{i} \hat{z}^{t}_{i}
        && \substack{\forall i \in \mathcal{I}, \\ \forall j \in \mathcal{J}, \\ \forall t \in \mathcal{T} : \\ i \succ_{j} 0} \label{eq:third-detailed-optimistic-ct2} \\
        & y^{t}_{i} \leq \sum_{\substack{\ell{} \in \mathcal{T}^{S}:\\\ell{} < t}} \sum_{\substack{k \in \mathcal{I}: \\k \succeq_{j} i}} \left(\hat{u}^{\ell{} t}_{kj} + \hat{v}^{\ell{} t}_{kj}\right)
        && \substack{\forall i \in \mathcal{I}, \\ \forall j \in \mathcal{J}, \\ \forall t \in \mathcal{T} : \\ i \succ_{j} 0} \label{eq:third-detailed-optimistic-ct3} \\  
        & \hat{z}^{t}_{i} \leq \sum_{\substack{\ell{} \in \mathcal{T}^{S}:\\\ell{} < t}} \sum_{\substack{k \in \mathcal{I}: \\k \succeq_{j} i}} \left(\hat{u}^{\ell{} t}_{kj} + \hat{v}^{\ell{} t}_{kj}\right)
        && \substack{\forall i \in \mathcal{I}, \\ \forall j \in \mathcal{J}, \\ \forall t \in \mathcal{T} : \\ i \succ_{j} 0} \label{eq:third-detailed-optimistic-ct4} \\  
        & y^{t}_{i} + \sum_{\substack{k \in \mathcal{I}: \\ i \succ_{j} k}} \sum_{\substack{\ell{} \in \mathcal{T}: \\ \ell{} < t}} \hat{u}^{\ell{} t}_{kj} \leq 1
        && \substack{\forall i \in \mathcal{I}, \\ \forall j \in \mathcal{J}, \\ \forall t \in \mathcal{T} : \\ i \succ_{j} 0} \label{eq:third-detailed-optimistic-ct5} \\ 
        & \hat{z}^{t}_{i} + \sum_{\substack{k \in \mathcal{I}: \\ i \succ_{j} k}} \sum_{\substack{\ell{} \in \mathcal{T}: \\ \ell{} < t}} \hat{u}^{\ell{} t}_{kj} \leq 1
        && \substack{\forall i \in \mathcal{I}, \\ \forall j \in \mathcal{J}, \\ \forall t \in \mathcal{T} : \\ i \succ_{j} 0} \label{eq:third-detailed-optimistic-ct6} \\ 
        & y^{t}_{i}  + \sum_{\substack{k \in \mathcal{I}: \\ i \succ_{j} k}} \sum_{\substack{\ell{} \in \mathcal{T}: \\ \ell{} < t}} \hat{v}^{\ell{} t}_{kj} \leq 1
        && \substack{\forall i \in \mathcal{I}, \\ \forall j \in \mathcal{J}, \\ \forall t \in \mathcal{T} : \\ i \succ_{j} 0} \label{eq:third-detailed-optimistic-ct7} \\ 
        & \hat{z}^{t}_{i} + \sum_{\substack{k \in \mathcal{I}: \\ i \succ_{j} k}} \sum_{\substack{\ell{} \in \mathcal{T}: \\ \ell{} < t}} \hat{v}^{\ell{} t}_{kj} \leq 1
        && \substack{\forall i \in \mathcal{I}, \\ \forall j \in \mathcal{J}, \\ \forall t \in \mathcal{T} : \\ i \succ_{j} 0} \label{eq:third-detailed-optimistic-ct8} \\               
        & \sum_{s \in \mathcal{T}^{F}} \sum_{i \in \mathcal{I}} \left(\hat{u}^{0s}_{ij} + \hat{v}^{0s}_{ij}\right) = 1
        && \substack{\forall j \in \mathcal{J}} \label{eq:third-detailed-optimistic-ct9} \\  
        & \sum_{\substack{s \in \mathcal{T}^{S}: \\ s < t}} \sum_{i \in \mathcal{I}} \left(\hat{u}^{st}_{ij} + \hat{v}^{st}_{ij}\right) =  
        \sum_{\substack{s \in \mathcal{T}^{F}: \\ s > t}}  \sum_{i \in \mathcal{I}} \left(\hat{u}^{ts}_{ij} + \hat{v}^{ts}_{ij}\right) 
        && \substack{\forall j \in \mathcal{J}, \\ \forall t \in \mathcal{T}} \label{eq:third-detailed-optimistic-ct10} \\
        & \hat{u}^{\ell{}t}_{ij}, \hat{v}^{\ell{}t}_{ij} \in [0,1]
        && \substack{\forall i \in \mathcal{I}, \\ \forall j \in \mathcal{J}, \\ \forall \ell{} \in \mathcal{T}^{S}, \\ \forall t \in \mathcal{T} :\\ i \succ_{j} 0, \ell{} < t}. \label{eq:third-detailed-optimistic-dm1}
        %
        %
        \end{align}
\end{subequations}

Objective Function \eqref{eq:upper-detailed-optimistic-obj} maximizes the profit of the leader obtained by serving customer demand.
Constraints \eqref{eq:upper-detailed-optimistic-ct1} guarantee that location decisions of the leader respect application-specific feasibility.
Constraints \eqref{eq:upper-detailed-optimistic-dm2} are variable domain constraints, and enforce an optimal reaction of the follower.
Objective Function \eqref{eq:lower-detailed-optimistic-obj} maximizes the profit of the follower against the location schedule $\boldsymbol{y}$ of the leader.
Constraints \eqref{eq:lower-detailed-optimistic-ct1} guarantee that location decisions of the follower are feasible.
Constraints \eqref{eq:third-detailed-optimistic-ct1}--\eqref{eq:third-detailed-optimistic-ct2} prevent customers from patronizing facilities that are not open by the leader or the follower, and force these players to split demand if they capture customers through the same location.
Constraints \eqref{eq:third-detailed-optimistic-ct3}--\eqref{eq:third-detailed-optimistic-ct4} force customers to patronize one of the available facilities if at least one of them is preferred over no service.
Bilinear terms $y^{t}_{i} z^{t}_{i}$ are linearized through additional variables and McCormick envelope constraints.
Constraints \eqref{eq:third-detailed-optimistic-ct5}--\eqref{eq:third-detailed-optimistic-ct8} enforce customer preferences (\ie, customer $j$ cannot patronize a location $k$ less preferred than location $i$ if location $i$ holds a facility during period $t$).
Constraints \eqref{eq:third-detailed-optimistic-ct9}--\eqref{eq:third-detailed-optimistic-ct10} preserve the captured flow for customers throughout the planning horizon. 
Finally, Constraints \eqref{eq:third-detailed-optimistic-dm1} define feasible variable domains.
We remark that, once $\boldsymbol{y}$ and $\boldsymbol{z}$ are fixed to binary values,  we can unambiguously compute optimal values of  variables $u^{\ell{}t}_{ij}$ and $v^{\ell{}t}_{ij}$. 
We therefore refer to solutions of the \ourproblem{} solely as the pair $(\boldsymbol{y}, \boldsymbol{z})$.

We can derive the detailed \textit{pessimistic} formulation of the \ourproblem{} by replacing the maximum in Objective Function~\eqref{eq:upper-detailed-optimistic-obj} by a max-min as explained in Section~\ref{sec:stylized}.
Similarly, one might easily account for time- or customer-dependent parameters in the \ourproblem{} (\eg, when they vary with the season or when customers have different perceptions about each player).
More specifically, we may consider time-dependent customer rankings $\succ^{t}_{j}$, time- and customer-dependent rewards $r^{t}_{ij}$,  as well as time- and customer-dependent splitting factors $\rho^{t}_{j}$.
Our mathematical formulations and solution methods can readily account for such extensions.

\section{Solution Methods}
\label{sec:methods}

In this section, we propose exact  solution methods for the optimistic and pessimistic variants of the \ourproblem{}.
We first study the theoretical complexity of the optimistic variant in Section~\ref{sec:tractability}, and propose an exact branch-and-cut algorithm for it in Section~\ref{sec:branch-n-cut}.
We then discuss how to adapt the exact branch-and-cut algorithm for the pessimistic variant in Section~\ref{sec:pessimistic}. 
For the sake of conciseness, we present the proofs of theoretical results in Appendix~\ref{apx:mathematical}.

\subsection{Theoretical Complexity}
\label{sec:tractability}

We first investigate the theoretical complexity of the optimistic variant of the \ourproblem{} to clarify the algorithmic limitations that we may face when devising solution methods.
To this end, we define the decision version of the problem under consideration as follows:
\vspace{0.3cm}
\boxxx{
{\bf Decision version of the optimistic variant of the \ourproblem{}}: \\
{\sc instance}: Finite sets $\mathcal{T} = \mathset{1, \ldots,T}$, $\mathcal{I}$,  $\mathcal{J}$, positive rational numbers $\{r_{i}\}_{i \in \mathcal{I}}$ and $\{d_{j}^t\}_{j \in \mathcal{J}, t\in \mathcal{T}}$,  rankings $\{\succ_{j}\}_{j \in \mathcal{J}}$, and a positive rational number $\Pi$.\\
{\sc question}: $\exists (\boldsymbol{y}, \boldsymbol{z}) \in \mathcal{Y} \times \mathcal{Z}, \forall \boldsymbol{z}^{\prime} \in \mathcal{Z}: \pi^{L}(\boldsymbol{y}, \boldsymbol{z}) \geq \Pi\land \pi^{F}(\boldsymbol{y}, \boldsymbol{z}) \geq \pi^{F}(\boldsymbol{y}, \boldsymbol{z}^{\prime})$, where feasible sets $\mathcal{Y}$ and $\mathcal{Z}$ are defined as in Section~\ref{sub:locations}?
}
\vspace{0.3cm}

Since the follower problem resembles an NP-hard problem \citep{silvaDynamicFacilityLocation2025}, one may naturally expect the optimistic variant of the \ourproblem{} to be $\Sigma^{p}_{2}$-hard.
Theorem~\ref{thm:sigmap2} states that this is indeed the case, proven through a reduction from the $\exists\forall$3SAT, which is $\Sigma^{p}_{2}$-complete \citep{wrathall1976complete}.

\begin{theorem}
    \label{thm:sigmap2}
    The decision version of the optimistic variant of the \ourproblem{} is $\Sigma^{p}_{2}$-complete.
\end{theorem}

In contrast to the majority of previously studied $\Sigma^{p}_{2}$-hard competitive problems \citep{grune2025completeness}, the novel $\Sigma^{p}_{2}$-complete competitive problem introduced in Theorem~\ref{thm:sigmap2} is neither a min-max nor a zero-sum game.
From a theoretical point of view, Theorem~\ref{thm:sigmap2} shows that the \ourproblem{} is harder than problems in the well-known NP-hard class, unless the second level of the polynomial hierarchy collapses \citep{stockmeyer1976polynomial}.
From a methodological point of view, Theorem~\ref{thm:sigmap2} reveals that there is no single-level mixed-integer linear reformulation of the bilevel formulation whose size is polynomial in the problem input, unless P = NP \citep{jeroslowPolynomialHierarchySimple1985}.
We therefore rely on bilevel programming methodologies to tackle the \ourproblem{}.
Lastly, from a practical perspective, the proof of Theorem~\ref{thm:sigmap2} suggests that the larger the number of periods, the harder it is.
We should thus expect exact solution methods to face scalability challenges for instances with many periods.

\subsection{Exact Branch-and-Cut Algorithm}
\label{sec:branch-n-cut}

Constraints~\eqref{eq:upper-detailed-optimistic-dm2} are the main source of intractability in the detailed optimistic formulation of the \ourproblem{}.
The standard approach in the literature \citep[see, \eg,][]{hemmatiMixedintegerBilevelProgramming2016, gentileIntegerProgrammingFormulations2018, beresnevApproximationCompetitiveFacility2019} is to relax these constraints and progressively add value-function cuts to approximate optimal follower decisions. 
Such an approach yields a cutting-plane framework based on the high-point relaxation $\mathcal{H}$, which we write as follows:
\begin{subequations}
    \label{mdl:relaxed-formulation}
    \begin{align}
        \mathcal{H}:  
        \max_{\boldsymbol{y}, \boldsymbol{z}, \boldsymbol{u}, \boldsymbol{v}} \quad
        \left\{   \sum_{t \in \mathcal{T}} \sum_{\substack{\ell{} \in \mathcal{T}^{S}: \\ \ell{} < t}}\sum_{j \in \mathcal{J}} \sum_{\substack{i \in \mathcal{I}: \\ i \succ_{j} 0}} r_{i} D^{\ell{}t}_{j}  u^{\ell{}t}_{ij} \, \, 
       \middle| \, \, 
       \boldsymbol{y} \in \mathcal{Y},
         \boldsymbol{z} \in \mathcal{Z},
        \eqref{eq:third-detailed-optimistic-ct1}-\eqref{eq:third-detailed-optimistic-dm1} \right\}. \nonumber 
        \end{align}
\end{subequations}

We implement an exact branch-and-cut algorithm through callbacks.
These callbacks are executed whenever the branch-and-bound solver finds an integer-feasible solution $(\boldsymbol{y}^{\prime}, \boldsymbol{z}^{\prime})$.
Such solution may be bilevel feasible ($\boldsymbol{z}^{\prime} \in \mathcal{F} (\boldsymbol{y}^{\prime})$) or bilevel infeasible ($\boldsymbol{z}^{\prime} \notin \mathcal{F} (\boldsymbol{y}^{\prime})$).
If it is the former, the solver can accept this integer-feasible solution as bilevel feasible and update the best incumbent; if it is the latter, the solver deems this integer-feasible solution as bilevel infeasible and removes it by adding a value-function cut. We highlight that these value-function cuts are globally valid despite being added at a specific node of the branch-and-bound tree.
Algorithm~\ref{alg:callback-optimistic} presents the described callback procedure for an integer-feasible solution $(\boldsymbol{y}^{\prime}, \boldsymbol{z}^{\prime})$.
\begin{algorithm}
    \caption{Optimistic Branch-and-Cut Callback}
    \label{alg:callback-optimistic}
    \begin{algorithmic}
        \REQUIRE Integer-feasible solution $(\boldsymbol{y}^{\prime}, \boldsymbol{z}^{\prime})$ at node $n$ and its constraint set $\mathcal{H}_n$ 
        \STATE Solve $\mathcal{F} (\boldsymbol{y}^{\prime})$ to obtain an optimal reaction $\boldsymbol{z}^{\star}$
        \IF{$\pi^{F} (\boldsymbol{y}^{\prime}, \boldsymbol{z}^{\prime}) = \pi^{F}(\boldsymbol{y}^{\prime}, \boldsymbol{z}^{\star})$}
            \STATE Update the best incumbent and move to the next node of the branch-and-bound tree \label{alg:adaptation}
        \ELSE
            \STATE Add value-function cut built with solution $(\boldsymbol{y}^{\prime}, \boldsymbol{z}^{\star})$  to $\mathcal{H}_s$ that cuts off solution $(\boldsymbol{y}^{\prime}, \boldsymbol{z}^{\prime})$
        \ENDIF
    \end{algorithmic}
\end{algorithm}

\begin{proposition}
    \label{prp:termination}
   The branch-and-cut framework with Algorithm~\ref{alg:callback-optimistic} as a callback terminates in a finite number of steps and returns an optimal solution to the optimistic variant of the \ourproblem{}.
\end{proposition}

In what follows, we discuss how to generate a value-function cut through a bilevel feasible solution $(\boldsymbol{y}^{\prime}, \boldsymbol{z}^{\star})$.
This cut, defined by some vector $(\alpha_0,\boldsymbol{\alpha})$, establishes $\pi^{F}(\boldsymbol{y},\boldsymbol{z}) \geq \boldsymbol{\alpha}^\top \boldsymbol{y}+\alpha_0$ for all bilevel feasible solutions $(\boldsymbol{y},\boldsymbol{z})$ with the inequality tight for solution $(\boldsymbol{y}^\prime,\boldsymbol{z}^\star)$; in simple terms, the cut enforces valid lower bounds on the profit of the follower (thus, the designation value-function cut).
Recall that we can compute optimal values for $u^{\ell{}t}_{ij}$ and $v^{\ell{}t}_{ij}$ once $\boldsymbol{y}$ and $\boldsymbol{z}$  are fixed to binary values.
We further refer to these optimal values as $u^{\ell{}t}_{ij}(\boldsymbol{y}, \boldsymbol{z})$ and  $v^{\ell{}t}_{ij}(\boldsymbol{y}, \boldsymbol{z})$ for a given pair $(\boldsymbol{y}, \boldsymbol{z})$.

\subsubsection{Decomposition of Profit in Intervals}
\label{sec:decomposition}

Before moving forward, we highlight an important property of the \ourproblem{} that is helpful to write value-function cuts.
Let $(\boldsymbol{y}^{\prime}, \boldsymbol{z}^{\star})$ be a bilevel feasible solution.
Once location schedules are fixed, the profit of the follower becomes separable by customer $j$.
Assume that the follower captures customer $j$ a total of $n$ times over the planning horizon -- through location $i^j_1$ at period $t^j_1$,  through location $i^j_2$ at period $t^j_2$, and so on until location $i^j_n$ at period $t^j_n$.
Since we consider a single customer $j$ at a time, in the following, we omit index $j$ from these location and period indexes.
We rewrite the profit of the follower from to customer $j$ as:
\begin{align*}
    \sum_{t \in \mathcal{T}} \sum_{\substack{\ell{} \in \mathcal{T}^{S}: \\ \ell{} < t}} \sum_{\substack{i \in \mathcal{I}: \\ i \succ_{j} 0}} r_{i} D^{\ell{}t}_{j}  v^{\ell{}t}_{ij} (\boldsymbol{y}^{\prime}, \boldsymbol{z}^{\star}) = 
        r_{i_1} D^{0t_1}_{j}  v^{0t_1}_{i_1j} (\boldsymbol{y}^{\prime}, \boldsymbol{z}^{\star}) +
        [\ldots] +
        r_{i_n} D^{t_{n-1}t_n}_{j}  v^{t_{n-1}t_n}_{i_nj} (\boldsymbol{y}^{\prime}, \boldsymbol{z}^{\star}).
\end{align*}

If the leader deviates to a location schedule $\boldsymbol{y} \neq \boldsymbol{y}^{\prime}$, her location decisions from period $t_{k-1}$ to period $t_{k}$ only affect 
the term  $r_{i_k} D^{t_{k-1}t_k}_{j}  v^{t_{k-1}t_k}_{i_kj} (\boldsymbol{y}^{\prime}, \boldsymbol{z}^{\star})$.
We can therefore analyze the profit of each customer $j$ at each period $t_k$ separately, and properly adjust the lower bound provided by the bilevel feasible solution $(\boldsymbol{y}^{\prime}, \boldsymbol{z}^{\star})$ to the profit of the follower.

\subsubsection{Tailored Value-Function Cut}

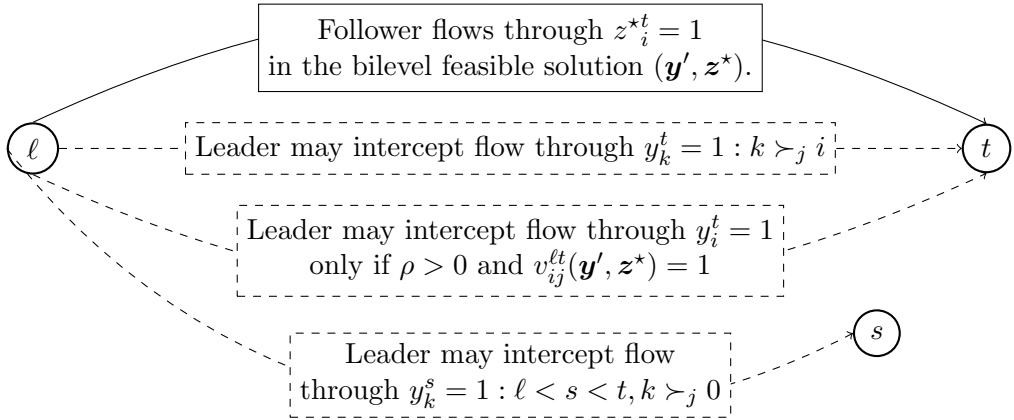
\begin{figure}[!ht]
    \centering
    \caption{Potential interferences of the leader on the profit of the follower under bilevel feasible solution $(\boldsymbol{y}^{\prime}, \boldsymbol{z}^{\star})$. Nodes represent periods, while arcs represent feasible transitions between periods, according to Constraints~\eqref{eq:third-detailed-optimistic-ct1}-\eqref{eq:third-detailed-optimistic-ct10}.}
    \begin{tikzpicture}[
    period/.style={circle, draw=black!100, thick},
    every text node part/.style={align=center}
    ]
    \begin{scope}
        \node[period] (before2) {$\ell{}$}; 
        
        \node[period] (after2) [right = 12cm of before2] {$t$}; 
    
        \node[period] (between2) [below left = 2cm and 1cm of after2] {$s$}; 
        
        \draw[->, dashed] (before2.east) -- (after2.west) node [midway, draw, color = black, fill = white, text = black] (label2) {Leader may intercept flow through $y^{t}_{k} = 1 : k \succ_{j} i$};

        \draw[->, solid, bend left = 25] (before2.north) to (after2.north) node [draw, color = black, fill = white, text = black, above = 0.4cm of label2]  {Follower flows through ${z^{\star}}^{t}_{i} = 1$ \\ in the bilevel feasible solution $(\boldsymbol{y}^{\prime}, \boldsymbol{z}^{\star})$. }; 

        \draw[->, dashed, bend right = 25] (before2.south) to (after2.south) node [draw, color = black, fill = white, text = black, below = 0.4cm of label2] {Leader may intercept flow through $y^{t}_{i} = 1$ \\ only if $\rho > 0$ and $v^{\ell{}t}_{ij} (\boldsymbol{y}^{\prime}, \boldsymbol{z}^{\star}) = 1$};
        
        \draw[->, dashed, bend right = 40] (before2.west) to (between2.west) node [draw, color = black, fill = white, text = black, below = 2.1cm of label2] {Leader may intercept flow \\ through $y^{s}_{k} = 1 : \ell < s < t, k \succ_{j} 0$};
    \end{scope}
    \end{tikzpicture}
    \label{fig:tailord-interdiction}
\end{figure}

We first propose a tailored value-function cut inspired by general value-function cuts \citep{lozanoValueFunctionBasedExactApproach2017}.
Let $(\boldsymbol{y}^{\prime}, \boldsymbol{z}^{\star})$ be a bilevel feasible solution.
Consider a customer $j$ captured by the follower through location $i$ at period $t$ after last being captured at period $\ell{}$ (\ie,  $v^{\ell{}t}_{ij} (\boldsymbol{y}^{\prime}, \boldsymbol{z}^{\star}) > 0$).
Figure~\ref{fig:tailord-interdiction} shows how the leader can interfere with contribution $r_{i} D^{\ell{}t}_{j} v^{\ell{}t}_{ij} (\boldsymbol{y}^{\prime}, \boldsymbol{z}^{\star})$, which can be partial (\ie, $v^{\ell{}t}_{ij} (\boldsymbol{y}^{\prime}, \boldsymbol{z}^{\star}) = 1 - \rho$) or integral (\ie, $v^{\ell{}t}_{ij} (\boldsymbol{y}^{\prime}, \boldsymbol{z}^{\star}) = 1$).

We therefore write the tailored value-function cut as follows, where $\mathbb{I}$ is an indicator function:
\begin{align}
      \sum_{t \in \mathcal{T}} \sum_{\substack{\ell{} \in \mathcal{T}^{S}: \\ \ell{} < t}} \sum_{j \in \mathcal{J}} \sum_{\substack{i \in \mathcal{I}: \\ i \succ_{j} 0}} r_{i} D^{\ell{}t}_{j}  v^{\ell{}t}_{ij} \geq \nonumber\\
        \sum_{t \in \mathcal{T}} \sum_{\substack{\ell{} \in \mathcal{T}^{S}: \\ \ell{} < t}} \sum_{j \in \mathcal{J}} \sum_{\substack{i \in \mathcal{I}: \\ i \succ_{j} 0}}  r_{i} D^{\ell{}t}_{j} v^{\ell{}t}_{ij} (\boldsymbol{y}^{\prime}, \boldsymbol{z}^{\star}) \left( 1 - \sum_{\substack{k \in \mathcal{I}: \\ k \succ_{j} 0}} \sum_{\substack{s \in \mathcal{T}: \\ \ell{} < s < t}} y^{s}_{k} - \sum_{\substack{k \in \mathcal{I}: \\ k \succ_{j} i}} y^{t}_{k} - \mathbb{I} \left[ v^{\ell{}t}_{ij} (\boldsymbol{y}^{\prime}, \boldsymbol{z}^{\star}) = 1 \land \rho > 0 \right] y^{t}_{i} \right) \label{eq:tailored-cut}
\end{align}

Inequality~\eqref{eq:tailored-cut} differs from general value-function cuts \citep{lozanoValueFunctionBasedExactApproach2017} by avoiding the use of large constants to subtract intercepted contributions.
Lemma~\ref{lem:tailored} guarantees that Inequality~\eqref{eq:tailored-cut} cuts the bilevel infeasible solution $(\boldsymbol{y}^{\prime}, \boldsymbol{z}^{\prime})$ within Algorithm~\ref{alg:callback-optimistic} and does not remove bilevel feasible solutions from the search space, ensuring finite convergence towards an optimal solution.

\begin{lemma}
    \label{lem:tailored}
    Inequality~\eqref{eq:tailored-cut}, built from a bilevel feasible solution $(\boldsymbol{y}^{\prime}, \boldsymbol{z}^{\star})$, holds for all bilevel feasible solutions of $\mathcal{L}^{optimstic}$ and is violated by any bilevel-infeasible (but integer-feasible) solution  $(\boldsymbol{y}^{\prime}, \boldsymbol{z}^{\prime})$.
\end{lemma}

We now discuss drawbacks of Inequality~\eqref{eq:tailored-cut}.
First, note that the term with binary values on the right-hand side can be negative.
This is the case, for example, when the leader captures customer $j$ at periods $s_1: \ell{} < s_1 < t$ and $s_2: \ell{} < s_2 < t$ in an alternative location schedule $\boldsymbol{y} \neq \boldsymbol{y}^{\prime}$.
Such a structure underestimates the profit of the follower by subtracting $r_{i} D^{\ell{}t}_{j} v^{\ell{}t}_{ij} (\boldsymbol{y}^{\prime}, \boldsymbol{z}^{\star})$ twice from the right-hand side.
Second, the follower may keep a portion of contribution $r_{i} D^{\ell{}t}_{j} v^{\ell{}t}_{ij} (\boldsymbol{y}^{\prime}, \boldsymbol{z}^{\star})$ after the interference of the leader.
For instance, if the leader does not capture customer $j$ between periods $s + 1$ and $t$, the follower still maintains contribution $r_{i} D^{st}_{j} v^{st}_{ij} (\boldsymbol{y}^{\prime\prime}, \boldsymbol{z}^{\star})$.
The right-hand side of Inequality~\eqref{eq:tailored-cut} cannot handle such subtlety, deducting contribution $r_{i} D^{\ell{}t}_{j} v^{\ell{}t}_{ij} (\boldsymbol{y}^{\prime}, \boldsymbol{z}^{\star})$ completely. 
These loose lower bounds may cause Algorithm~\ref{alg:callback-optimistic} to produce multiple value-function cuts for the same location schedule $\boldsymbol{z}^{\star}$ of the follower, each against a different location schedule $\boldsymbol{y}$ of the leader, thus slowing convergence of the exact branch-and-cut algorithm to an optimal bilevel solution.

\subsubsection{Tightened Value-Function Cut}
\label{sec:improved}

We now introduce a general tightening technique that addresses shortcomings such as the two sources of underestimation in the tailored value-function cut.
Let $(\boldsymbol{y}^{\prime}, \boldsymbol{z}^{\star})$ be a bilevel feasible solution.
The tightened value-function cut considers the profit of the follower with location schedule $\boldsymbol{z}^{\star}$ in the ideal scenario (\ie, when the leader deploys an empty location schedule $\boldsymbol{0}$), and discounts exact amounts based on potential interferences of the leader.
More precisely, we write the tightened value-function cut as follows:
\begin{align} 
    \sum_{t \in \mathcal{T}} \sum_{\substack{\ell{} \in \mathcal{T}^{S}: \\ \ell{} < t}} \sum_{j \in \mathcal{J}} \sum_{\substack{i \in \mathcal{I}: \\ i \succ_{j} 0}} r_{i} D^{\ell{}t}_{j}  v^{\ell{}t}_{ij}  \geq  \nonumber \\
        \sum_{t \in \mathcal{T}} \sum_{\substack{\ell{} \in \mathcal{T}^{S} :\\ \ell{} < t}} \sum_{j \in \mathcal{J}} \sum_{\substack{i \in \mathcal{I}: \\ i \succ_{j} 0}}  r_{i} v^{\ell{}t}_{ij} (\boldsymbol{0}, \boldsymbol{z}^{\star}) \left( D^{\ell{}t}_{j} - \sum_{\substack{s \in \mathcal{T}: \\ \ell{} < s \leq t}} d^{s}_{j} o^{\ell{}st}_{ij} - \rho y^{t}_{i}\sum_{\substack{s \in \mathcal{T}: \\ \ell{} < s \leq t}} d^{s}_{j} (1 - o^{\ell{}st}_{ij}) \right),
        \label{eq:improved}
\end{align}
where $o^{\ell{}st}_{ij} = \max \left\{ \max_{\substack{s^{\prime} \in \mathcal{T} \\ k \in \mathcal{I}: \\ s \leq s^{\prime} < t \\k \succ_{j} 0}} \left\{ y^{s^{\prime}}_{k}\right\}, \max_{\substack{k \in \mathcal{I}: \\ k \succ_{j} i}} \left\{ y^{t}_{k}\right\} \right\}$ computes discounts triggered by the leader.

The right-hand side of Inequality~\eqref{eq:improved} can be interpreted as follows. 
Consider a customer $j$ captured by the follower through location $i$ at period $t$ after last being captured at period $\ell{}$ (\ie,  $v^{\ell{}t}_{ij} (\boldsymbol{0}, \boldsymbol{z}^{\star}) = 1$).
Recall that contribution $r_{i} D^{\ell{}t}_{j} v^{\ell{}t}_{ij} (\boldsymbol{0}, \boldsymbol{z}^{\star}) = r_{i} D^{\ell{}t}_{j}$ can be decomposed as $\sum_{\substack{s \in \mathcal{T} \\ \ell{} < s \leq t}} r_{i} d^{s}_{j}$.
In this sense, maximization terms $o^{\ell{}st}_{ij} \in \mathset{0,1}$ equal $1$ if the leader integrally intercepts portion $r_{i} d^{s}_{j}$ of contribution $r_{i} D^{\ell{}t}_{j} $, and $0$ otherwise.
We can then see these maximization terms $o^{\ell{}st}_{ij}$ as auxiliary variables and 
enforce their values through the following constraints:
\begin{align}
    o^{\ell{}st}_{ij} \leq \sum_{\substack{s^{\prime} \in \mathcal{T}: \\s \leq s^{\prime} < t}} \sum_{\substack{k \in \mathcal{I}: \\ k \succ_{j} 0}} y^{s^{\prime}}_{k} + \sum_{\substack{k \in \mathcal{I}: \\ k \succ_{j} i}} y^{t}_{k} && \substack{\forall  i \in \mathcal{I}, \forall j \in \mathcal{J}, \forall \ell{} \in \mathcal{T}^{S}, \\ \forall s \in \mathcal{S}, \forall t \in \mathcal{T}: \ell{} < t, \ell{} < s \leq t} \label{eq:upper-auxiliary1} \\
    o^{\ell{}st}_{ij} \in [0,1] && \substack{\forall  i \in \mathcal{I}, \forall j \in \mathcal{J}, \forall \ell{} \in \mathcal{T}^{S}, \\ \forall s \in \mathcal{S}, \forall t \in \mathcal{T}: \ell{} < t, \ell{} < s \leq t}. \label{eq:upper-auxiliary2}
\end{align}
Note that we do not need to impose the equality to the maximization term nor binary constraints on auxiliary variables $o^{\ell{}st}_{ij}$, as the leader will always set these variables to $1$ whenever feasible to avoid loosing customers to the follower.
In turn, bilinear terms $y^{t}(1 - o^{\ell{}st}_{ij}) \in \mathset{0,1}$ equal $1$ if the leader did not integrally intercept portion $r_{i} d^{s}_{j}$ of contribution $r_{i} D^{\ell{}t}_{j}$, but rather splits it with the follower through location $i$ at period $t$, $0$ otherwise.
These bilinear terms can be easily linearized through additional variables and McCormick envelope constraints, so that Inequality~\eqref{eq:improved} can be replaced by a set of linear constraints.
Similarly to variables $u^{\ell{}t}_{ij}$ and $v^{\ell{}t}_{ij}$, once $\boldsymbol{y}$ is fixed to binary values, variables $o^{\ell{}st}_{ij}$ can be computed unambiguously.
We further refer to these optimal values as $o^{\ell{}st}_{ij} (\boldsymbol{y})$.

Lemma~\ref{lem:improved} guarantees that Inequality~\eqref{eq:improved} cuts the bilevel infeasible solution $(\boldsymbol{y}^{\prime}, \boldsymbol{z}^{\prime})$ within Algorithm~\ref{alg:callback-optimistic} and does not remove bilevel feasible solutions from the search space.

\begin{lemma}
    \label{lem:improved}
    Inequality~\eqref{eq:improved}, built from a bilevel feasible solution $(\boldsymbol{y}^{\prime}, \boldsymbol{z}^{\star})$, holds for all bilevel feasible solutions of $\mathcal{L}^{optimstic}$ and is violated by any bilevel-infeasible (but integer-feasible) solution $(\boldsymbol{y}^{\prime}, \boldsymbol{z}^{\prime})$.
\end{lemma}

The major advantage of the tightened value-function cut is the construction of a tight lower bound to the profit of the follower with location schedule $\boldsymbol{z}^{\star}$ against any location schedule $\boldsymbol{y}$ of the leader.
In this sense, Algorithm~\ref{alg:callback-optimistic} equipped with Inequality~\eqref{eq:improved} will never enumerate the same location schedule $\boldsymbol{z}^{\star}$ more than once, and may converge to an optimal bilevel solution faster than with Inequality~\eqref{eq:tailored-cut}.
Theorem~\ref{thm:dominance} formally states that the former theoretically dominates the latter.

\begin{theorem}
    \label{thm:dominance}
    Inequality~\eqref{eq:improved} dominates Inequality~\eqref{eq:tailored-cut} for the follower problem $\mathcal{F}(\boldsymbol{y})$.
\end{theorem}

From a practical perspective, Inequality~\eqref{eq:improved} requires the addition of a large number of variables and constraints to the high-point relaxation $\mathcal{H}$, leading to large linear programs at each node of the branch-and-cut tree.
Within a short time limit, this may prevent the tightened value-function cut from finding bilevel feasible solutions of quality as high as the tailored value-function cut because the former explores less nodes than the latter -- however, as the time limit increases, the former should outperform the latter in terms of solution quality.
We conclude by mentioning that the detailed optimistic formulation could theoretically be rewritten as a single-level formulation with an exponential number of constraints by adding Inequality~\eqref{eq:improved}, along with its additional variables and constraints, for each location schedule $\boldsymbol{z} \in \mathcal{Z}$ to the high-point relaxation $\mathcal{H}$.
However, the exponential size of the feasible set $\mathcal{Z}$ would render the model intractable, which makes a branch-and-cut approach preferable in practice, particularly as the number of periods increases.

Although Inequality~\eqref{eq:improved} is problem-specific and cannot be directly applied to other competitive problems, we can use the same intuition to design value-function cuts tighter than general value-function cuts \citep{lozanoValueFunctionBasedExactApproach2017} or tailored value-function cuts \citep[see, \eg,][]{hemmatiMixedintegerBilevelProgramming2016, gentileIntegerProgrammingFormulations2018, beresnevApproximationCompetitiveFacility2019} for other planning problems under competition.
This involves three main steps.
The first step is to verify whether the profit of the follower is fixed to some value once decisions of the leader and the follower have been fixed, as discussed at the end of Section~\ref{sec:detailed}.
This is often the case for competitive FLPs, where the profit of each customer to the follower is determined once location decisions of both players are fixed. 
The second step is to identify which decisions of the leader may trigger discounts on the profit of the follower, and group them in a maximization term, as explained in Section~\ref{sec:improved}.
In competitive FLPs, this is usually the opening of a location that is strictly more preferred by the customer.
The third step is to linearize the maximization terms, similarly to what is done in Section~\ref{sec:improved}.

\subsection{Pessimistic Variant of the \ourproblem{}}
\label{sec:pessimistic}

We now explain how to adapt Algorithm~\ref{alg:callback-optimistic} for the pessimistic variant of the \ourproblem{}.
Let $(\boldsymbol{y}^{\prime}, \boldsymbol{z}^{\prime})$ be the integer-feasible solution obtained by the callback.
In addition to checking whether $\boldsymbol{z}^{\prime} \in \mathcal{F} (\boldsymbol{y}^{\prime})$, we need to ensure that $\boldsymbol{z}^{\prime} \in \argmin_{\boldsymbol{z} \in\mathcal{F} (\boldsymbol{y}^{\prime})} \mathset{\pi^{L} (\boldsymbol{y}^{\prime}, \boldsymbol{z})}$ in Line \ref{alg:adaptation}.
In other words, we need to verify whether $\boldsymbol{z}^{\prime}$ is a follower-optimal location schedule minimizing  the profit of the leader with location schedule $\boldsymbol{y}^{\prime}$.
If this is the case, then solution $(\boldsymbol{y}^{\prime}, \boldsymbol{z}^{\prime})$ is pessimistically bilevel feasible; otherwise, solution $(\boldsymbol{y}^{\prime}, \boldsymbol{z}^{\prime})$ is optimistically bilevel feasible, but pessimistically bilevel infeasible.
In the latter case, we add a no-good cut \citep[see, \eg,][]{beresnevApproximationCompetitiveFacility2019} to remove $(\boldsymbol{y}^{\prime}, \boldsymbol{z}^{\prime})$:
\begin{align*}
    \sum_{t \in \mathcal{T}} \sum_{\substack{i \in \mathcal{I} \\ {y^{\prime}}^{t}_{i} = 1}} (1 - y^{t}_{i}) + \sum_{t \in \mathcal{T}} \sum_{\substack{i \in \mathcal{I} \\ {y^{\prime}}^{t}_{i} = 0}} y^{t}_{i} + \sum_{t \in \mathcal{T}} \sum_{\substack{i \in \mathcal{I} \\ {z^{\prime}}^{t}_{i} = 1}} (1 - z^{t}_{i}) + \sum_{t \in \mathcal{T}} \sum_{\substack{i \in \mathcal{I} \\ {z^{\prime}}^{t}_{i} = 0}} z^{t}_{i} \geq 1
\end{align*}

Note that pessimistically bilevel feasible solutions must also be optimistically bilevel feasible, as 
both detailed formulations have the same set of feasible solutions and only differ in their objective functions.
In other words, pessimistically bilevel feasible solutions must respect both the tailored and the tightened value-function cuts for the follower problem $\mathcal{F}(\boldsymbol{y})$, and the dominance relationship established by  Theorem~\ref{thm:dominance} also applies for the pessimistic variant of the \ourproblem{}.
\begin{corollary}
    Inequalities \eqref{eq:tailored-cut} and \eqref{eq:improved}, built from a bilevel feasible solution $(\boldsymbol{y}^{\prime}, \boldsymbol{z}^{\star})$, hold for all pessimistically bilevel feasible solutions of $\mathcal{L}^{pessimistic}$ and Inequalitiy \eqref{eq:improved} dominates Inequality \eqref{eq:tailored-cut} for the follower problem $\mathcal{F}(\boldsymbol{y}^{\prime})$.
\end{corollary}

Contrary to Lemmas~\ref{lem:tailored} and ~\ref{lem:improved}, value-function cuts cannot solely remove pessimistically bilevel infeasible solutions .
However, their combination with no-good cuts guarantees that the branch-and-cut algorithm with the adapted version of Algorithm~\ref{alg:callback-optimistic} never visits the same integer-feasible solution twice, thus converging towards an optimal solution of the pessimistic variant of the \ourproblem{}.

\begin{corollary}
    The branch-and-cut framework with (the adapted version of) Algorithm~\ref{alg:callback-optimistic} as a callback terminates in a finite number of steps and returns an optimal solution to the pessimistic variant of the \ourproblem{}.
\end{corollary}

\section{Computational Experiments}
\label{sec:experiments}

In this section, we study the performance of our solution methods and draw managerial insights about the \ourproblem{}.
We first present the computational benchmark based on real-world data in Section~\ref{sec:benchmark}. 
Section~\ref{sec:performance} then evaluates the performance of the exact solution methods for our computational benchmark in terms of solution quality and computing times, paying particular attention to the benefits of the tightened value-function cut.
Lastly, Section~\ref{sec:insights} focuses on managerial insights on the impact of cumulative demand, duopolistic and optimistic behaviour assumptions, as well as the revenue splitting factor $\rho$.


We implement our solution methods in Python (version \texttt{3.10}), and solve the mixed-integer programming formulations with Gurobi (version \texttt{12.0}).
All jobs were processed on the Nibi cluster of the Digital Research Alliance of Canada with a maximum RAM of 30GB and a single thread. 
Throughout our computational experiments, we refer to the branch-and-cut algorithm with Inequality~\eqref{eq:tailored-cut} as Tailored B\&C, and to the same algorithm with Inequality~\eqref{eq:improved} as Tightened B\&C.

\subsection{Computational Benchmark}
\label{sec:benchmark}

Since the \ourproblem{} is an extension of competitive location problems considered in the literature, no benchmark instances are readily available. 
We therefore generate synthetic instances based on real-world data, 
whose sizes are comparable to those studied in other (bilevel) competitive FLPs \citep[see, \eg,][]{roboredoBranchandcutAlgorithmDiscrete2013, hemmatiMixedintegerBilevelProgramming2016, gentileIntegerProgrammingFormulations2018, qiSequentialCompetitiveFacility2022}. 
We restrict each player to a single facility per period, as suggested in Section~\ref{sub:locations}, which renders the optimistic variant of the \ourproblem{} $\Sigma^{p}_{2}$-hard and allows for focusing on the relationship between cumulative demand and market competition rather than analyzing the trade-off with problem-specific constraints on location decisions. 
%
%
%

We create our computational benchmark instances using federal electoral districts of the Canadian province of Quebec, which provide a realistic and diverse spatial distribution of customers and potential facility locations while relying solely on publicly available data.
We consider two geographic scopes: the entirety of Quebec with $|\mathcal{J}| = 78$ customers, and the Montreal region  with $|\mathcal{J}| = 40$ customers.
We then sample a subset of $\frac{|\mathcal{J}|}{2}$ customers to build the set of locations. 
We create instances with $|\mathcal{T}| \in \mathset{3, 5, 7}$ periods to examine the impact of cumulative customer demand across different planning horizons.
As more than 8 out of 10 Canadians use a vehicle to commute \citep{statcan2023_commuting}, we consider that customers are willing to access facilities that are within a maximum driving time of $M \in \mathset{15, 30,45}$ minutes.
%
Customers prefer nearby locations rather than distant ones.
We examine identical $\left(r_{i} = |\mathcal{I}|, \forall i \in \mathcal{I}\right)$ and inverse $\left(r_{i} = \ceil{\frac{|\mathcal{I}|}{\sum_{\substack{j \in \mathcal{J} \\ i \succ_{j} 0}} 1}}, \forall i \in \mathcal{I}\right)$ rewards per unit of captured demand.
Intuitively, the former describes applications where the reward is independent of location, whereas the latter describes applications where popular locations tend to have larger costs (\eg, higher rent) and, consequently, smaller rewards.
We also examine constant ($d^{t}_{j} = P_{j}, \forall j \in \mathcal{J}, \forall t \in \mathcal{T}$) and sparse ($d^{t}_{j} \sim \mathcal{U} \mathset{0,P_{j}}, \forall j \in \mathcal{J}, \forall t \in \mathcal{T}$) spawning demands, where $P_{j} \in \mathbb{R}^{+}$ is the population size of customer $j$ and $\mathcal{U}$ is a discrete uniform distribution.
The former exemplifies scenarios where customers have demand appearing throughout the entire planning horizon, whereas the latter exemplifies scenarios where customers may not have demand appearing at some time periods.
Lastly, we consider splitting factors $\rho \in \mathset{0, 0.25, 0.5, 0.75, 1}$ to understand how brand recognition influences the structure of optimal solutions.
Since some parameters are generated randomly, we consider seed values $S \in \mathset{1,2,3,4,5}$.
The combination of the aforementioned parameters yields a benchmark with a total of $2 \cdot 3^2 \cdot 2^2 \cdot 5^2 = 1800$ instances.
Details on how to manipulate raw data from Statistics Canada and OpenStreetMap are provided in Appendix~\ref{apx:benchmark}.

\subsection{Computational Performance}
\label{sec:performance}

In preliminary results, we compared our solution methods with the state-of-the-art general-purpose bilevel solver MIX++ \citep{fischettiNewGeneralPurposeAlgorithm2017}. 
The MIX++ Solver employs locally valid value-function cuts, which may yield more conservative lower bounds on the follower's objective than the globally valid cuts of~\cite{lozanoValueFunctionBasedExactApproach2017}, and combines these cuts with intersection cuts and other enhancements, resulting in effective general performance. 
%
%
We found, however, that our solution methods significantly outperform this general-purpose bilevel solver within a time limit of $1$ hour.
%
%
More specifically, the MIX++ Solver manages to solve to optimality only $2\%$ of an applicable subset of 720 instances, whereas the Tightened B\&C succeeds to solve $53\%$ of them, within the same time limit. 
Moreover, the MIX++ Solver can only tackle an optimistic follower behaviour, preventing us from drawing insights on the pessimistic variant of the \ourproblem{}, and requires a particular structure on the lower level, keeping us from considering instances with splitting factors $\rho \in \{0.25, 0.5, 0.75\}$.
We therefore only consider the Tailored B\&C and the Tightened B\&C in our computational experiments  (for more details on the preliminary results, see Appendix~\ref{apx:preliminary}).

We now first evaluate the performance of the branch-and-cut algorithm for the optimistic variant of the \ourproblem{} within a time limit of $24$ hours.
Such a time limit is reasonable for real-world applications with a planning horizon of several weeks or months, and allows us to find optimal solutions for approximately $74\%$ of the computational benchmark. 
Table~\ref{tab:summary-optimistic} presents the percentage of the benchmark solved to optimality by each exact method, grouped by dimensional instance attributes (\ie, number of locations $|\mathcal{I}|$, customers $|\mathcal{J}|$, and periods $|\mathcal{T}|$).
The Tightened B\&C solves more instances to optimality than the Tailored B\&C, and that for different instance dimensions. 
Recall that the proof of Theorem~\ref{thm:sigmap2} suggests that the larger the number of periods, the harder (the optimistic variant of) the \ourproblem{} is to solve.
Table~\ref{tab:summary-optimistic} confirms the practical implications of such theoretical hardness -- we solve all instances with $|\mathcal{T}| = 3$ periods to optimality, but gradually less instances with an increasing length of the planning horizon. 
%
%
\begin{table}[!ht]
    \centering
    \small
    \caption{Percentage of the benchmark solved to optimality by each exact method, grouped by dimensional instance attributes (\ie, number of locations $|\mathcal{I}|$, customers $|\mathcal{J}|$, and periods $|\mathcal{T}|$). [Optimistic variant of the \ourproblem{}]}
    \begin{tabular}{ccccc}
    \toprule
        \small
        \multirow{2}{*}{$|\mathcal{T}|$}& \multicolumn{2}{c}{Montreal ($|\mathcal{I}| = 20, |\mathcal{J}| = 40$)} & \multicolumn{2}{c}{Quebec ($|\mathcal{I}| = 39, |\mathcal{J}| = 78$)}  \\ \cmidrule(lr){2-5}
        & Tailored B\&C& Tightened B\&C&  Tailored B\&C & Tightened B\&C\\ \midrule
        3&100.00\%&100.00\%&100.00\%&100.00\%\\ \midrule
5&66.67\%&82.67\%&70.67\%&86.33\%\\ \midrule
7&12.67\%&40.33\%&15.00\%&31.33\%\\ \bottomrule

    \end{tabular}
    \label{tab:summary-optimistic}
\end{table}

A comparison on the $1092$ instances solved to optimality by both approaches has shown that the Tightened B\&C solves instances approximately three times faster than the Tailored B\&C, on average ($36.90 \pm 130.25$ versus $97.55 \pm 238.11$ minutes). For the remaining $708$ instances that were not solved to optimality by both exact methods, the Tightened B\&C proves optimality gaps that are almost half as large, on average, than those proven by the Tailored B\&C ($ 16.97\% \pm 20.91\%$ versus $28.94\% \pm 19.10\%$). Note that we are bound to obtain large values for standard deviations, as solution times vary strongly with instance attributes. Detailed results, separated by instance characteristics, can be found in Appendix~\ref{apx:supplementary-optimistic}.


At first glance, 
these results provide empirical evidence that the tightened value-function cut should be strictly preferred over the tailored value-function cut, in line with the theoretical results.
%
%
We now analyze performance graphs for computing time and solution quality to compare these methods on an instance-to-instance basis. 
More specifically, we compute the \textit{runtime ratio} of each exact method as $\frac{\Delta^{\prime}}{\Delta^{b}}$, where $\Delta^{b}$ is the lowest computing time among the exact methods, and $\Delta^{\prime}$ is the computing time taken by the exact method at hand \footnote{\label{not:ratios} Small (respectively, large) runtime ratios indicate that the exact method at hand has a computing time closer to (respectively, farther from) the fastest one. Similarly, small (respectively, large) objective ratios indicate that the exact method at hand finds a solution with an objective value closer to (respectively, farther from) the best one.}.
Similarly, we compute the \textit{objective ratio} of each exact method as $\frac{\Pi^{b}}{\Pi^{\prime}}$, where $\Pi^{b}$ is the highest objective value found among the exact methods, and $\Pi^{\prime}$ is the objective value obtained by the exact method at hand \footref{not:ratios}.
Figures~\ref{fig:runtimes-optimistic} and \ref{fig:objectives-optimistic} present both ratios, respectively, comparing both exact methods.
\begin{figure}[!ht]
     \centering
     \caption{Runtime ratios for exact methods, where the $y$ axis presents the percentage of instances with a ratio smaller than or equal to the reference value on the $x$ axis. [Optimistic variant of the \ourproblem{}]}
    \begin{tikzpicture}
\begin{axis}[name = 3periods, xmin = 1, xmax = 2, xscale = 0.5, ymin = -5, ymax = 105, yscale = 0.75, xtick = {1,1.25,1.5,1.75,2}, ytick = {0, 25, 50, 75, 100}, xlabel = {Runtime ratio}, ylabel = {Instances (\%)}, legend pos = south west]
\addplot[line width = 2pt, color = gray, style = dotted] coordinates {(1,17) (1.1,38)};
\addplot[line width = 2pt, color = gray, style = dotted] coordinates {(1.1,38) (1.2,42)};
\addplot[line width = 2pt, color = gray, style = dotted] coordinates {(1.2,42) (1.3,46)};
\addplot[line width = 2pt, color = gray, style = dotted] coordinates {(1.3,46) (1.4,50)};
\addplot[line width = 2pt, color = gray, style = dotted] coordinates {(1.4,50) (1.5,53)};
\addplot[line width = 2pt, color = gray, style = dotted] coordinates {(1.5,53) (1.6,55)};
\addplot[line width = 2pt, color = gray, style = dotted] coordinates {(1.6,55) (1.7,58)};
\addplot[line width = 2pt, color = gray, style = dotted] coordinates {(1.7,58) (1.8,61)};
\addplot[line width = 2pt, color = gray, style = dotted] coordinates {(1.8,61) (1.9,63)};
\addplot[line width = 2pt, color = gray, style = dotted] coordinates {(1.9,63) (2.0,65)};
\addplot[line width = 2pt, color = orange, style = solid] coordinates {(1,36) (1.1,74)};
\addplot[line width = 2pt, color = orange, style = solid] coordinates {(1.1,74) (1.2,79)};
\addplot[line width = 2pt, color = orange, style = solid] coordinates {(1.2,79) (1.3,82)};
\addplot[line width = 2pt, color = orange, style = solid] coordinates {(1.3,82) (1.4,86)};
\addplot[line width = 2pt, color = orange, style = solid] coordinates {(1.4,86) (1.5,90)};
\addplot[line width = 2pt, color = orange, style = solid] coordinates {(1.5,90) (1.6,92)};
\addplot[line width = 2pt, color = orange, style = solid] coordinates {(1.6,92) (1.7,94)};
\addplot[line width = 2pt, color = orange, style = solid] coordinates {(1.7,94) (1.8,95)};
\addplot[line width = 2pt, color = orange, style = solid] coordinates {(1.8,95) (1.9,96)};
\addplot[line width = 2pt, color = orange, style = solid] coordinates {(1.9,96) (2.0,97)};

\end{axis}
\begin{axis}[name = 5periods, at={(3periods.east)}, anchor=west, xshift = 1.5cm, xmin = 1, xmax = 2, xscale = 0.5, ymin = -5, ymax = 105, yscale = 0.75, xtick = {1,1.25,1.5,1.75,2}, ytick = {0, 25, 50, 75, 100}, xlabel = {Runtime ratio}, ylabel = {Instances (\%)}, legend pos = south west]
\addplot[line width = 2pt, color = gray, style = dotted] coordinates {(1,15) (1.1,32)};
\addplot[line width = 2pt, color = gray, style = dotted] coordinates {(1.1,32) (1.2,35)};
\addplot[line width = 2pt, color = gray, style = dotted] coordinates {(1.2,35) (1.3,38)};
\addplot[line width = 2pt, color = gray, style = dotted] coordinates {(1.3,38) (1.4,40)};
\addplot[line width = 2pt, color = gray, style = dotted] coordinates {(1.4,40) (1.5,43)};
\addplot[line width = 2pt, color = gray, style = dotted] coordinates {(1.5,43) (1.6,44)};
\addplot[line width = 2pt, color = gray, style = dotted] coordinates {(1.6,44) (1.7,45)};
\addplot[line width = 2pt, color = gray, style = dotted] coordinates {(1.7,45) (1.8,47)};
\addplot[line width = 2pt, color = gray, style = dotted] coordinates {(1.8,47) (1.9,48)};
\addplot[line width = 2pt, color = gray, style = dotted] coordinates {(1.9,48) (2.0,49)};
\addplot[line width = 2pt, color = orange, style = solid] coordinates {(1,42) (1.1,86)};
\addplot[line width = 2pt, color = orange, style = solid] coordinates {(1.1,86) (1.2,88)};
\addplot[line width = 2pt, color = orange, style = solid] coordinates {(1.2,88) (1.3,89)};
\addplot[line width = 2pt, color = orange, style = solid] coordinates {(1.3,89) (1.4,91)};
\addplot[line width = 2pt, color = orange, style = solid] coordinates {(1.4,91) (1.5,91)};
\addplot[line width = 2pt, color = orange, style = solid] coordinates {(1.5,91) (1.6,92)};
\addplot[line width = 2pt, color = orange, style = solid] coordinates {(1.6,92) (1.7,93)};
\addplot[line width = 2pt, color = orange, style = solid] coordinates {(1.7,93) (1.8,93)};
\addplot[line width = 2pt, color = orange, style = solid] coordinates {(1.8,93) (1.9,93)};
\addplot[line width = 2pt, color = orange, style = solid] coordinates {(1.9,93) (2.0,94)};

\end{axis}
\begin{axis}[name = 7periods, at={(5periods.east)}, anchor=west, xshift = 1.5cm, xmin = 1, xmax = 2, xscale = 0.5, ymin = -5, ymax = 105, yscale = 0.75, xtick = {1,1.25,1.5,1.75,2}, ytick = {0, 25, 50, 75, 100}, xlabel = {Runtime ratio}, ylabel = {Instances (\%)}, legend pos = south west]
\addlegendimage{line width = 2, color = gray, style = dotted}
\addlegendentry{Tailored B\&C}
\addlegendimage{line width = 2, color = orange, style = solid}
\addlegendentry{Tightened B\&C}
\addplot[line width = 2pt, color = gray, style = dotted] coordinates {(1,29) (1.1,65)};
\addplot[line width = 2pt, color = gray, style = dotted] coordinates {(1.1,65) (1.2,66)};
\addplot[line width = 2pt, color = gray, style = dotted] coordinates {(1.2,66) (1.3,67)};
\addplot[line width = 2pt, color = gray, style = dotted] coordinates {(1.3,67) (1.4,68)};
\addplot[line width = 2pt, color = gray, style = dotted] coordinates {(1.4,68) (1.5,69)};
\addplot[line width = 2pt, color = gray, style = dotted] coordinates {(1.5,69) (1.6,69)};
\addplot[line width = 2pt, color = gray, style = dotted] coordinates {(1.6,69) (1.7,70)};
\addplot[line width = 2pt, color = gray, style = dotted] coordinates {(1.7,70) (1.8,70)};
\addplot[line width = 2pt, color = gray, style = dotted] coordinates {(1.8,70) (1.9,71)};
\addplot[line width = 2pt, color = gray, style = dotted] coordinates {(1.9,71) (2.0,72)};
\addplot[line width = 2pt, color = orange, style = solid] coordinates {(1,28) (1.1,98)};
\addplot[line width = 2pt, color = orange, style = solid] coordinates {(1.1,98) (1.2,98)};
\addplot[line width = 2pt, color = orange, style = solid] coordinates {(1.2,98) (1.3,98)};
\addplot[line width = 2pt, color = orange, style = solid] coordinates {(1.3,98) (1.4,98)};
\addplot[line width = 2pt, color = orange, style = solid] coordinates {(1.4,98) (1.5,98)};
\addplot[line width = 2pt, color = orange, style = solid] coordinates {(1.5,98) (1.6,99)};
\addplot[line width = 2pt, color = orange, style = solid] coordinates {(1.6,99) (1.7,99)};
\addplot[line width = 2pt, color = orange, style = solid] coordinates {(1.7,99) (1.8,99)};
\addplot[line width = 2pt, color = orange, style = solid] coordinates {(1.8,99) (1.9,99)};
\addplot[line width = 2pt, color = orange, style = solid] coordinates {(1.9,99) (2.0,99)};

\end{axis}
\node[align=center, yshift=1cm] at (3periods.north) {Instances with $|\mathcal{T}| = 3$.};
\node[align=center, yshift=0.5cm] at (3periods.north) {(total of $600$ instances)};
\node[align=center, yshift=1cm] at (5periods.north) {Instances with $|\mathcal{T}| = 5$.};
\node[align=center, yshift=0.5cm] at (5periods.north) {(total of $600$ instances)};
\node[align=center, yshift=1cm] at (7periods.north) {Instances with $|\mathcal{T}| = 7$.};
\node[align=center, yshift=0.5cm] at (7periods.north) {(total of $600$ instances)};
\end{tikzpicture}
    \label{fig:runtimes-optimistic}
\end{figure}
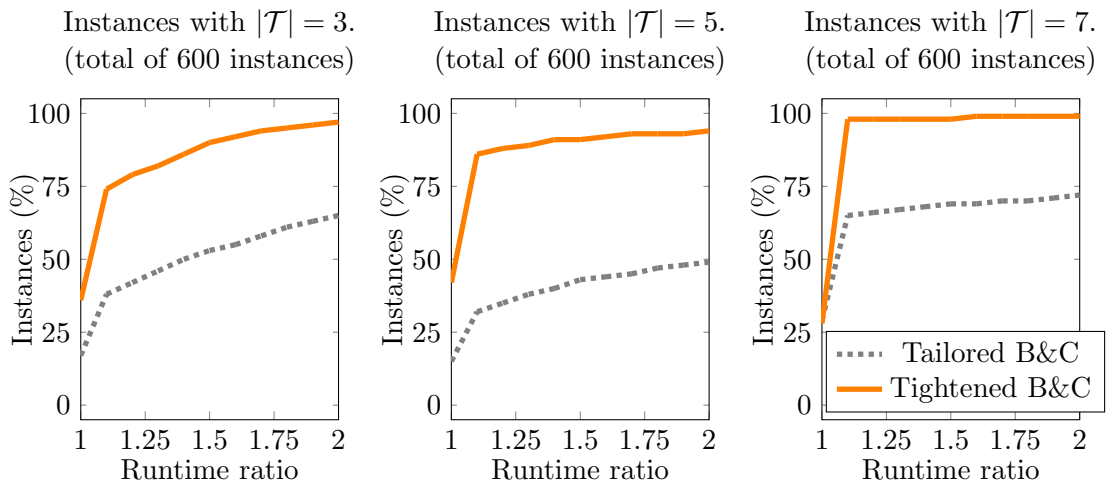
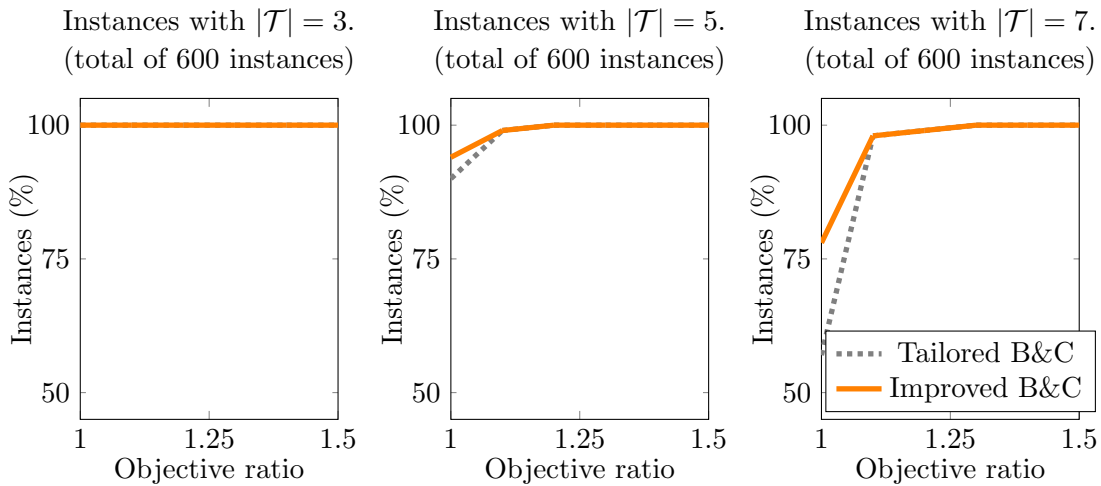
\begin{figure}[!ht]
     \centering
     \caption{Objective ratios for exact methods, where the $y$ axis presents the percentage of instances with a ratio smaller than or equal to the reference value on the $x$ axis. [Optimistic variant of the \ourproblem{}]}    
    \begin{tikzpicture}
\begin{axis}[name = 3periods, xmin = 1, xmax = 1.5, xscale = 0.5, ymin = 45, ymax = 105, yscale = 0.75, xtick = {1,1.25,1.5,1.75,2}, ytick = {50, 75, 100}, xlabel = {Objective ratio}, ylabel = {Instances (\%)}, legend pos = south west]
\addplot[line width = 2pt, color = gray, style = dotted] coordinates {(1,100) (1.1,100)};
\addplot[line width = 2pt, color = gray, style = dotted] coordinates {(1.1,100) (1.2,100)};
\addplot[line width = 2pt, color = gray, style = dotted] coordinates {(1.2,100) (1.3,100)};
\addplot[line width = 2pt, color = gray, style = dotted] coordinates {(1.3,100) (1.4,100)};
\addplot[line width = 2pt, color = gray, style = dotted] coordinates {(1.4,100) (1.5,100)};
\addplot[line width = 2pt, color = orange, style = solid] coordinates {(1,100) (1.1,100)};
\addplot[line width = 2pt, color = orange, style = solid] coordinates {(1.1,100) (1.2,100)};
\addplot[line width = 2pt, color = orange, style = solid] coordinates {(1.2,100) (1.3,100)};
\addplot[line width = 2pt, color = orange, style = solid] coordinates {(1.3,100) (1.4,100)};
\addplot[line width = 2pt, color = orange, style = solid] coordinates {(1.4,100) (1.5,100)};

\end{axis}
\begin{axis}[name = 5periods, at={(3periods.east)}, anchor=west, xshift = 1.5cm, xmin = 1, xmax = 1.5, xscale = 0.5, ymin = 45, ymax = 105, yscale = 0.75, xtick = {1,1.25,1.5,1.75,2}, ytick = {50, 75, 100}, xlabel = {Objective ratio}, ylabel = {Instances (\%)}, legend pos = south west]
\addplot[line width = 2pt, color = gray, style = dotted] coordinates {(1,90) (1.1,99)};
\addplot[line width = 2pt, color = gray, style = dotted] coordinates {(1.1,99) (1.2,100)};
\addplot[line width = 2pt, color = gray, style = dotted] coordinates {(1.2,100) (1.3,100)};
\addplot[line width = 2pt, color = gray, style = dotted] coordinates {(1.3,100) (1.4,100)};
\addplot[line width = 2pt, color = gray, style = dotted] coordinates {(1.4,100) (1.5,100)};
\addplot[line width = 2pt, color = orange, style = solid] coordinates {(1,94) (1.1,99)};
\addplot[line width = 2pt, color = orange, style = solid] coordinates {(1.1,99) (1.2,100)};
\addplot[line width = 2pt, color = orange, style = solid] coordinates {(1.2,100) (1.3,100)};
\addplot[line width = 2pt, color = orange, style = solid] coordinates {(1.3,100) (1.4,100)};
\addplot[line width = 2pt, color = orange, style = solid] coordinates {(1.4,100) (1.5,100)};

\end{axis}
\begin{axis}[name = 7periods, at={(5periods.east)}, anchor=west, xshift = 1.5cm, xmin = 1, xmax = 1.5, xscale = 0.5, ymin = 45, ymax = 105, yscale = 0.75, xtick = {1,1.25,1.5,1.75,2}, ytick = {50, 75, 100}, xlabel = {Objective ratio}, ylabel = {Instances (\%)}, legend pos = south west]
\addlegendimage{line width = 2, color = gray, style = dotted}
\addlegendentry{Tailored B\&C}
\addlegendimage{line width = 2, color = orange, style = solid}
\addlegendentry{Improved B\&C}
\addplot[line width = 2pt, color = gray, style = dotted] coordinates {(1,57) (1.1,98)};
\addplot[line width = 2pt, color = gray, style = dotted] coordinates {(1.1,98) (1.2,99)};
\addplot[line width = 2pt, color = gray, style = dotted] coordinates {(1.2,99) (1.3,100)};
\addplot[line width = 2pt, color = gray, style = dotted] coordinates {(1.3,100) (1.4,100)};
\addplot[line width = 2pt, color = gray, style = dotted] coordinates {(1.4,100) (1.5,100)};
\addplot[line width = 2pt, color = orange, style = solid] coordinates {(1,78) (1.1,98)};
\addplot[line width = 2pt, color = orange, style = solid] coordinates {(1.1,98) (1.2,99)};
\addplot[line width = 2pt, color = orange, style = solid] coordinates {(1.2,99) (1.3,100)};
\addplot[line width = 2pt, color = orange, style = solid] coordinates {(1.3,100) (1.4,100)};
\addplot[line width = 2pt, color = orange, style = solid] coordinates {(1.4,100) (1.5,100)};

\end{axis}
\node[align=center, yshift=1cm] at (3periods.north) {Instances with $|\mathcal{T}| = 3$.};
\node[align=center, yshift=0.5cm] at (3periods.north) {(total of $600$ instances)};
\node[align=center, yshift=1cm] at (5periods.north) {Instances with $|\mathcal{T}| = 5$.};
\node[align=center, yshift=0.5cm] at (5periods.north) {(total of $600$ instances)};
\node[align=center, yshift=1cm] at (7periods.north) {Instances with $|\mathcal{T}| = 7$.};
\node[align=center, yshift=0.5cm] at (7periods.north) {(total of $600$ instances)};
\end{tikzpicture}
    \label{fig:objectives-optimistic}
\end{figure}
On the one hand, Figure~\ref{fig:runtimes-optimistic} shows that the Tightened B\&C is considerably faster than the Tailored B\&C across our computational benchmark.
For $|\mathcal{T}| = 7$ periods, both methods hit the time limit for about $25\%$ of the instances, resulting in a less pronounced gap between the runtime ratios.
%
%
On the other hand, Figure~\ref{fig:objectives-optimistic} reveals that both the Tightened B\&C and the Tailored B\&C find solutions of similar quality within the time limit of $24$ hours.
In particular, the former finds solutions of better quality than the latter within the same time limit as the size of the planning horizon increases, 
%
empirically confirming its superior performance. 

Finally, we note that detailed performance results on the pessimistic variant of the \ourproblem{} can be found in Appendix~\ref{apx:supplementary-pessimistic}.
Both branch-and-cut algorithms struggle more to solve the pessimistic variant, since it needs to additionally exclude solutions that are optimistically bilevel feasible but pessimistically bilevel infeasible.
Nevertheless, we observe the same trends between Tightened B\&C  and Tailored B\&C as in the optimistic variant, favouring the former over the latter.

\subsection{Managerial Insights}
\label{sec:insights}

As the \ourproblem{} aims to provide decision-support for rather recently introduced application domains, such as temporary pop-up stores, we now derive managerial insights in regards to the key features of this young class of planning problems. 
We focus on the optimistic variant of the \ourproblem{} and consider the subset of $1325$ instances (approximately $74\%$ of the computational benchmark) that have been solved to optimality by at least one of the exact methods.
%
%

\subsubsection{Impact of a Monopolistic Assumption}

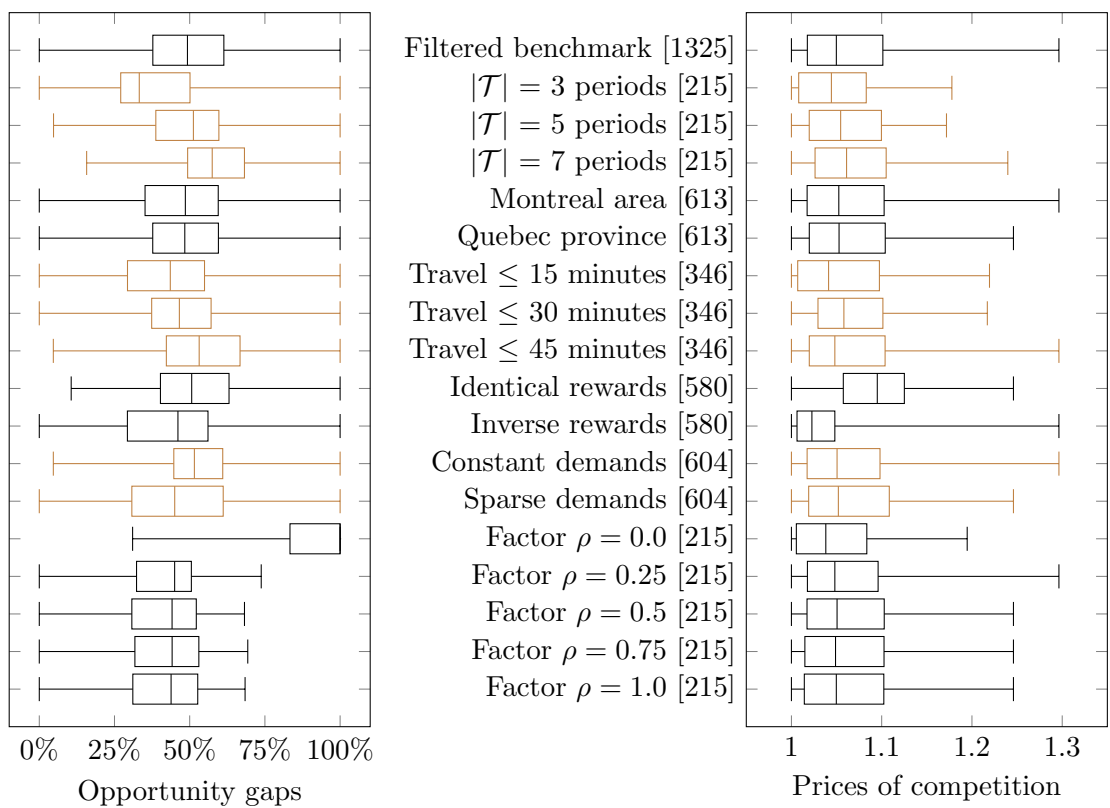
\begin{figure}
    \centering
    \caption{Opportunity gaps under the monopolistic assumption (left) and prices of competition (right) grouped by different instance attributes, with the number of instances between square brackets. [Optimistic variant of the \ourproblem{}]}
    \begin{tikzpicture}
\begin{axis}[name = oppgaps, y=0.5cm, ymin=0, ymax=19, x = 0.04cm, xmin = -10, xmax = 110, xlabel = {Opportunity gaps}, xtick = {0,25,50,75,100},  xticklabels = {0\%, 25\%, 50\%, 75\%, 100\%}, ytick = {1,2,3,4,5,6,7,8,9,10,11,12,13,14,15,16,17,18}, yticklabels = {}]
\addplot+[
	boxplot prepared={
	lower whisker=0.0,
	lower quartile=31.0948,
	median=43.8129,
	upper quartile=52.67505,
	upper whisker=68.41
}, color = black, style = solid
] coordinates {};
\addplot+[
	boxplot prepared={
	lower whisker=0.0,
	lower quartile=31.783549999999998,
	median=44.1595,
	upper quartile=53.04145,
	upper whisker=69.3041
}, color = black, style = solid
] coordinates {};
\addplot+[
	boxplot prepared={
	lower whisker=0.0,
	lower quartile=30.7692,
	median=44.14,
	upper quartile=52.19265,
	upper whisker=68.2105
}, color = black, style = solid
] coordinates {};
\addplot+[
	boxplot prepared={
	lower whisker=0.0,
	lower quartile=32.33385,
	median=44.9974,
	upper quartile=50.55135,
	upper whisker=73.7762
}, color = black, style = solid
] coordinates {};
\addplot+[
	boxplot prepared={
	lower whisker=31.0246,
	lower quartile=83.3333,
	median=100.0,
	upper quartile=100.0,
	upper whisker=100.0
}, color = black, style = solid
] coordinates {};
\addplot+[
	boxplot prepared={
	lower whisker=0.0,
	lower quartile=30.7692,
	median=45.0202,
	upper quartile=61.1413,
	upper whisker=100.0
}, color = brown, style = solid
] coordinates {};
\addplot+[
	boxplot prepared={
	lower whisker=4.668,
	lower quartile=44.718,
	median=51.557649999999995,
	upper quartile=60.985125,
	upper whisker=100.0
}, color = brown, style = solid
] coordinates {};
\addplot+[
	boxplot prepared={
	lower whisker=0.0,
	lower quartile=29.28975,
	median=46.09135,
	upper quartile=56.0837,
	upper whisker=100.0
}, color = black, style = solid
] coordinates {};
\addplot+[
	boxplot prepared={
	lower whisker=10.6061,
	lower quartile=40.272400000000005,
	median=50.66605,
	upper quartile=63.074225,
	upper whisker=100.0
}, color = black, style = solid
] coordinates {};
\addplot+[
	boxplot prepared={
	lower whisker=4.668,
	lower quartile=42.2399,
	median=53.15195,
	upper quartile=66.73712499999999,
	upper whisker=100.0
}, color = brown, style = solid
] coordinates {};
\addplot+[
	boxplot prepared={
	lower whisker=0.0,
	lower quartile=37.363075,
	median=46.5646,
	upper quartile=57.114050000000006,
	upper whisker=100.0
}, color = brown, style = solid
] coordinates {};
\addplot+[
	boxplot prepared={
	lower whisker=0.0,
	lower quartile=29.3226,
	median=43.549549999999996,
	upper quartile=54.927,
	upper whisker=100.0
}, color = brown, style = solid
] coordinates {};
\addplot+[
	boxplot prepared={
	lower whisker=0.0,
	lower quartile=37.7149,
	median=48.3699,
	upper quartile=59.5429,
	upper whisker=100.0
}, color = black, style = solid
] coordinates {};
\addplot+[
	boxplot prepared={
	lower whisker=0.0,
	lower quartile=35.1605,
	median=48.5537,
	upper quartile=59.4855,
	upper whisker=100.0
}, color = black, style = solid
] coordinates {};
\addplot+[
	boxplot prepared={
	lower whisker=15.7449,
	lower quartile=49.3158,
	median=57.4863,
	upper quartile=68.20085,
	upper whisker=100.0
}, color = brown, style = solid
] coordinates {};
\addplot+[
	boxplot prepared={
	lower whisker=4.7355,
	lower quartile=38.7396,
	median=51.229,
	upper quartile=59.70075,
	upper whisker=100.0
}, color = brown, style = solid
] coordinates {};
\addplot+[
	boxplot prepared={
	lower whisker=0.0,
	lower quartile=27.0703,
	median=33.237,
	upper quartile=50.104299999999995,
	upper whisker=100.0
}, color = brown, style = solid
] coordinates {};
\addplot+[
	boxplot prepared={
	lower whisker=0.0,
	lower quartile=37.7446,
	median=49.24895,
	upper quartile=61.3558,
	upper whisker=100.0
}, color = black, style = solid
] coordinates {};
\end{axis}

\begin{axis}[at={(oppgaps.east)}, anchor=west, xshift = 5cm,y=0.5cm, ymin=0, ymax=19, x=12cm, xmin = 0.95, xmax = 1.35, xlabel = {Prices of competition}, xtick = {1,1.1,1.2,1.3},  ytick = {1,2,3,4,5,6,7,8,9,10,11,12,13,14,15,16,17,18}, yticklabels = {Factor $\rho = 1.0$ [215],Factor $\rho = 0.75$ [215],Factor $\rho = 0.5$ [215],Factor $\rho = 0.25$ [215],Factor $\rho = 0.0$ [215],Sparse demands [604],Constant demands [604],Inverse rewards [580],Identical rewards [580],Travel $\leq$ 45 minutes [346],Travel $\leq$ 30 minutes [346],Travel $\leq$ 15 minutes [346],Quebec province [613],Montreal area [613],$|\mathcal{T}|$ = 7 periods [215],$|\mathcal{T}|$ = 5 periods [215],$|\mathcal{T}|$ = 3 periods [215],Filtered benchmark [1325]}]
\addplot+[
	boxplot prepared={
	lower whisker=1.0,
	lower quartile=1.0142720483508616,
	median=1.049727463312369,
	upper quartile=1.1023648648648647,
	upper whisker=1.2459677419354838
}, color = black, style = solid
] coordinates {};
\addplot+[
	boxplot prepared={
	lower whisker=1.0,
	lower quartile=1.0149038542148945,
	median=1.0488676996424315,
	upper quartile=1.1025900900900902,
	upper whisker=1.2459677419354838
}, color = black, style = solid
] coordinates {};
\addplot+[
	boxplot prepared={
	lower whisker=1.0,
	lower quartile=1.0172794162900314,
	median=1.050561797752809,
	upper quartile=1.1028146699708303,
	upper whisker=1.2459677419354838
}, color = black, style = solid
] coordinates {};
\addplot+[
	boxplot prepared={
	lower whisker=1.0,
	lower quartile=1.0177077209484615,
	median=1.0481048577798797,
	upper quartile=1.0960391850550022,
	upper whisker=1.2963619021665076
}, color = black, style = solid
] coordinates {};
\addplot+[
	boxplot prepared={
	lower whisker=1.0,
	lower quartile=1.005444968007685,
	median=1.0380952380952382,
	upper quartile=1.0835254529567324,
	upper whisker=1.1947261663286004
}, color = black, style = solid
] coordinates {};
\addplot+[
	boxplot prepared={
	lower whisker=1.0,
	lower quartile=1.0191137395283612,
	median=1.0519932541533603,
	upper quartile=1.1084272315105954,
	upper whisker=1.2459677419354838
}, color = brown, style = solid
] coordinates {};
\addplot+[
	boxplot prepared={
	lower whisker=1.0,
	lower quartile=1.017391304347826,
	median=1.050561797752809,
	upper quartile=1.0981449095520843,
	upper whisker=1.2963619021665076
}, color = brown, style = solid
] coordinates {};
\addplot+[
	boxplot prepared={
	lower whisker=1.0,
	lower quartile=1.0061996280223187,
	median=1.0226914591868894,
	upper quartile=1.0481048577798797,
	upper whisker=1.2963619021665076
}, color = black, style = solid
] coordinates {};
\addplot+[
	boxplot prepared={
	lower whisker=1.0,
	lower quartile=1.0574712643678161,
	median=1.0951194147341967,
	upper quartile=1.125101543460601,
	upper whisker=1.2459677419354838
}, color = black, style = solid
] coordinates {};
\addplot+[
	boxplot prepared={
	lower whisker=1.0,
	lower quartile=1.0197870102550617,
	median=1.0481048577798797,
	upper quartile=1.1037735849056605,
	upper whisker=1.2963619021665076
}, color = brown, style = solid
] coordinates {};
\addplot+[
	boxplot prepared={
	lower whisker=1.0,
	lower quartile=1.0294218268901814,
	median=1.0581280704117684,
	upper quartile=1.1013443640124094,
	upper whisker=1.217032967032967
}, color = brown, style = solid
] coordinates {};
\addplot+[
	boxplot prepared={
	lower whisker=1.0,
	lower quartile=1.0069629612199378,
	median=1.0411717056887515,
	upper quartile=1.0974683544303798,
	upper whisker=1.2194323144104804
}, color = brown, style = solid
] coordinates {};
\addplot+[
	boxplot prepared={
	lower whisker=1.0,
	lower quartile=1.0197870102550617,
	median=1.05277152791982,
	upper quartile=1.103942652329749,
	upper whisker=1.2459677419354838
}, color = black, style = solid
] coordinates {};
\addplot+[
	boxplot prepared={
	lower whisker=1.0,
	lower quartile=1.017391304347826,
	median=1.0524956970740103,
	upper quartile=1.1027027027027028,
	upper whisker=1.2963619021665076
}, color = black, style = solid
] coordinates {};
\addplot+[
	boxplot prepared={
	lower whisker=1.0,
	lower quartile=1.0261845386533666,
	median=1.0611033519553073,
	upper quartile=1.1049737771452666,
	upper whisker=1.239733629300777
}, color = brown, style = solid
] coordinates {};
\addplot+[
	boxplot prepared={
	lower whisker=1.0,
	lower quartile=1.019736842105263,
	median=1.0545454545454545,
	upper quartile=1.0997054272708637,
	upper whisker=1.171785028790787
}, color = brown, style = solid
] coordinates {};
\addplot+[
	boxplot prepared={
	lower whisker=1.0,
	lower quartile=1.0080624355005159,
	median=1.044280442804428,
	upper quartile=1.0829015544041452,
	upper whisker=1.1777777777777778
}, color = brown, style = solid
] coordinates {};
\addplot+[
	boxplot prepared={
	lower whisker=1.0,
	lower quartile=1.017633009994121,
	median=1.0497581202487907,
	upper quartile=1.1013001018030815,
	upper whisker=1.2963619021665076
}, color = black, style = solid
] coordinates {};
\end{axis}
\end{tikzpicture}
    \label{fig:oppgaps-prices-optimistic}
\end{figure}

We first examine whether the leader can simply ignore market competition at the planning stage, even though competition is present. 
The leader would solve the \ourproblem{} assuming that the follower does not exist (\ie, $\boldsymbol{z} = \boldsymbol{0}$)  to obtain a location schedule $\boldsymbol{y}^{\prime}$. 
Once the leader announces her decisions to the public, she would then face an optimistic reaction $\boldsymbol{z}^{\star} \neq \boldsymbol{0}$ of the follower. 
We refer to this solution method as the monopolistic heuristic.
Let $\frac{\Pi^{\star} - \Pi^{\prime}}{\Pi^{\star}}$ be the opportunity gap for each instance, where $\Pi^{\star}$ is the optimal objective value obtained through an exact method and $\Pi^{\prime}$ is the objective value obtained by the monopolistic heuristic solution\footnote{Small (respectively, large) opportunity gaps indicate that the leader obtains a profit closer to (respectively, farther from) the optimal profit in a duopoly even if she ignores competition at the planning stage.}.
The box plots on the left of Figure~\ref{fig:oppgaps-prices-optimistic} present these gaps grouped by different instance attributes.
We find an average opportunity gap of $52.19\% \pm 23.20\%$, meaning that the leader could obtain twice the profit by properly accounting for market competition at the planning stage rather than ignoring it.
Falsely assuming a monopoly is particularly damaging to the leader with a splitting factor $\rho = 0$, as the follower may steal the entire profit of the leader by replicating her location schedule.
Furthermore, ignoring the duopoly is significantly suboptimal for instances with large planning horizons, long maximum travel times, constant spawning demands, and identical rewards, as they have larger opportunity gaps than their counterparts.
%

\subsubsection{Cooperation as an Alternative}

Since the \ourproblem{} is not a zero-sum game, players may jointly collect more profit if they cooperate rather than compete with each other.
We therefore assess the advantage of cooperation from the perspective of the players through the price of competition for each instance, defined as their (joint) profit under cooperation divided by the sum of their (separate) profits under competition\footnote{Small (respectively, large) prices of competition indicate that players could extract close to the same (respectively, considerably more) profit from customers if they were to cooperate rather than compete with each other.}.
The box plots on the right of Figure~\ref{fig:oppgaps-prices-optimistic} present the price of competition  grouped by different instance attributes.
We obtain an average price of competition of $1.06 \pm 0.05$, indicating that players could jointly extract $6\%$ more revenue from customers if they cooperate rather than compete with each other. Players would split this surplus according to a preset rule.
%
%
Note that competition seems particularly damaging for instances with identical rewards.
In a competitive setting, identical rewards seem to encourage players to locate their temporary facilities close to each other throughout the planning horizon, thus capturing the same subset of customers multiple times over the planning horizon.
In a cooperative setting, players would rather choose locations that are farther from each other as to capture a more comprehensive subset of customers throughout the planning horizon, which explains the high price of competition for those instances.
%

Cooperating rather than competing also affects customer service quality, represented by two metrics: the average number of captures and the average percentage of captured demand.
%
%
%
%
We report these two metrics in Appendix~\ref{apx:managerial-optimistic}.
On average, cooperation seems to decrease the average number of captures throughout the planning horizon, while increasing the percentage of demand captured by the end of the planning horizon.
%
%
%
These insights are valuable not only for players deciding how to operate their temporary facilities, but also for policy makers that may consider imposing regulations.
For example, if players know that customers value continuous customer engagement (\ie, being visited often), they would not actively consider cooperating with each other.
In turn, if the government prioritizes long-term customer satisfaction (i.e., ensuring most demand is eventually served), then it should not only facilitate but also encourage cooperation (e.g., providing a framework and benefits for companies to come to a mutually benefic agreement).

\subsubsection{Structure of Optimal Location Schedules}

We now analyze how instance attributes affect the structure of the location schedules.
For the sake of conciseness, we chose an \textit{illustrative instance} of moderate difficulty: Montreal region, $|\mathcal{T}| = 5$ periods, maximum travel time $M = 30$ minutes, identical rewards, constant spawning demands, and splitting factor $\rho = 0.5$.
Figure~\ref{fig:illustrative} presents optimal solutions of the illustrative instance and some of its variations, where gray regions are federal election districts, each representing a customer, and coloured squares denote locations chosen by the leader and the follower at each of the five periods.
Before discussing it, we address a potentially counterintuitive characteristic of Stackelberg competitions. 
Although one may naturally expect that the leader has a first-mover advantage \citep[page 78]{10.1561/2400000040} over the follower, \ie, that by moving first she always earns a higher profit him, such an advantage is not ensured in a Stackelberg competition \cite[see, \eg,][]{hu2024whoplaysfirst}.
We may thus see optimal solutions where, even though the leader solves the \ourproblem{} to establish a location schedule under competition, she obtains a profit smaller than the one of the follower. 
\begin{figure}[ht!]
\caption{Location schedules of the illustrative instance and its variations, where $\Pi^{L}$ and $\Pi^{F}$ indicate the objective function values of the leader and the follower, respectively.
    [Optimistic variant of the \ourproblem{}]} 
\vspace{-0.3cm}

\begin{subfigure}{0.48\textwidth}
\caption{Leader wrongly assumes monopoly.}
\include{maps_solutions/map-qc_1-5-1-30-montreal-identical-constant-driving-0.5-monopoly}
\vspace{-0.8cm}
\label{fig:illustrative-monopoly}
\end{subfigure}
\hfill
\begin{subfigure}{0.48\textwidth}
\caption{Leader assumes an optimistic duopoly.}

\include{maps_solutions/map-qc_1-5-1-30-montreal-identical-constant-driving-0.5-duopoly}
\label{fig:illustrative-duopoly}
\vspace{-0.8cm}
\end{subfigure}
\begin{subfigure}{0.48\textwidth}
\caption{
Inverse rewards rather than identical ones.
}
\include{maps_solutions/map-qc_1-5-1-30-montreal-inversely-constant-driving-0.5}
\label{fig:illustrative-different}
\vspace{-0.8cm}
\end{subfigure}
\hfill
\begin{subfigure}{0.48\textwidth}
\caption{
Sparse spawning demands rather than constant ones.
}
\include{maps_solutions/map-qc_1-5-1-30-montreal-identical-sparsed-driving-0.5}
\label{fig:illustrative-sparse}
\vspace{-0.8cm}
\end{subfigure}
\begin{subfigure}{0.48\textwidth}
\caption{
Splitting factor $\rho = 0.25$ rather than $\rho = 0.5$.
}
\include{maps_solutions/map-qc_1-5-1-30-montreal-identical-constant-driving-0.25}
\label{fig:illustrative-0.25}
\vspace{-1cm}
\end{subfigure}
\hfill
\begin{subfigure}{0.48\textwidth}
\caption{
Splitting factor $\rho = 0.75$ rather than $\rho = 0.5$.
}
\include{maps_solutions/map-qc_1-5-1-30-montreal-identical-constant-driving-0.75}
\label{fig:illustrative-0.75}
\vspace{-1cm}
\end{subfigure}
\label{fig:illustrative}
\end{figure}

We first compare the solution structure when the leader assumes to either operate within a monopoly or within a duopoly.
Figure~\ref{fig:illustrative-monopoly} reports the optimal location schedule for the leader (blue) within a monopoly, and the respective reaction of the follower (red).
The leader places the single facility in suburban areas early in the planning horizon and leaves central districts build up demands for late periods.
Such a strategy is bound to underperform under a duopoly, as the follower captures central customers early in the planning horizon and considerably decreases the profit of the leader.
Figure~\ref{fig:illustrative-duopoly} then presents optimal locations schedules for the leader (blue) and the follower (red) within a duopoly.
Both players move their respective temporary facilities around central locations, capturing similar customers in early periods and distancing from each other in later periods. 
In this context, some customers tend to be captured multiple times, whereas others are neglected for most of the planning horizon.
%
%
This comparison underlines once again the importance of accounting for market competition for the leader -- not only the respective location schedules are structurally quite different, the leader obtains approximately half the profit that she could have obtained by properly accounting for market competition.

We now investigate how instance attributes change optimal location schedules within a duopoly.
Figure~\ref{fig:illustrative-different} presents the optimal solution for the illustrative instance, but with inverse rewards rather than identical ones (see Figure~\ref{fig:illustrative-duopoly}).
This is the case, for example, when opening a facility requires renting out a space, as higher rents of central locations implies lower profits.
In this context, both players seek to locate their temporary facilities in suburban locations 
rather than in central ones. 
As a result, more customers tend to be captured within the planning horizon, although less often than with identical rewards.
As a rule of thumb, the leader should prefer suburban locations over central ones when rewards are inverse, guided by the stark difference in terms of their profitability.
While Figure~\ref{fig:illustrative-duopoly} assumed a constant demand over the planning horizon, 
Figure~\ref{fig:illustrative-sparse} reports the optimal solution with sparse spawning demands, describing, for instance, a company selling a new product for which demand among the population is not yet established.
Players are similarly motivated to be less conservative (\ie, move the single facility far from the previous locations), since capturing many customers no longer necessarily guarantees the capture of a large amount of demand, but not as intensely.
Therefore, when spawning demands are not constant, the leader should expect to move the temporary facility often as to follow demand shifts throughout the planning horizon.
The splitting factor $\rho$ also naturally impacts the structure of optimal solutions.
Figures~\ref{fig:illustrative-0.25} and~\ref{fig:illustrative-0.75} present the optimal solution for splitting factors $\rho= 0.25$ and $\rho= 0.75$ rather than $\rho= 0.5$, respectively (see Figure~\ref{fig:illustrative-duopoly}). 
When $\rho= 0.75$, the leader retains a larger percentage of customer demand split with the follower, so there is less incentive to leave central locations.
As a consequence, she barely moves her temporary facility throughout the planning horizon.
In turn, when $\rho= 0.25$, the leader loses a larger percentage of customer demand split with the follower, so there is more incentive to explore alternative locations.
Brand recognition plays thus an important role in the optimal solution, prompting a less preferred leader to relocate more often than a established one.

%
%

\subsubsection{Impact of Duopolistic Assumptions}

\begin{figure}[!ht]
     \centering
     \caption{Resilience ratios averaged over our computational benchmark. [Optimistic variant of the \ourproblem{}]}
     \label{fig:heatmaps-o}
     \begin{subfigure}{0.44\textwidth}
        \caption{Under an optimistic ground truth.}
        \begin{tikzpicture}
\begin{axis}[
    colormap/blackwhite,
    xlabel = {Ground truth $(optimistic, \rho_2)$},
    ylabel = {Assumption $(optimistic, \rho_1)$},
    point meta min=0,
    point meta max=100,
    xmin=-0.125, xmax=1.125, xscale=0.8,
    xtick style={color=black},
    xtick={0,0.25,0.5,0.75,1},
    ymin=-0.125, ymax=1.125, yscale=0.8,
    ytick style={color=black},
    ytick={0,0.25,0.5,0.75,1},
    every node near coord/.append style={font=\footnotesize, color=black, anchor=center,
        /pgf/number format/.cd, fixed, fixed zerofill, precision=2
    }
    ]
    \addplot[%
        matrix plot*,
        point meta=explicit,
        mesh/cols=5,
        visualization depends on=\thisrow{Deviation}\as\std,,
        nodes near coords={\pgfmathprintnumber{\pgfplotspointmeta}},
    ] table[meta=Mean] {
    x y Mean Deviation
0.0 0.0 100.00 0.00
0.25 0.0 88.26 10.98
0.5 0.0 81.41 17.18
0.75 0.0 80.53 18.02
1.0 0.0 80.50 18.03
0.0 0.25 44.00 31.97
0.25 0.25 100.00 0.00
0.5 0.25 92.71 9.82
0.75 0.25 91.56 11.10
1.0 0.25 91.53 11.13
0.0 0.5 40.22 31.83
0.25 0.5 83.03 17.07
0.5 0.5 100.00 0.00
0.75 0.5 98.96 2.86
1.0 0.5 98.93 2.99
0.0 0.75 42.22 32.07
0.25 0.75 83.54 16.93
0.5 0.75 97.91 5.55
0.75 0.75 100.00 0.00
1.0 0.75 100.00 0.00
0.0 1.0 42.18 31.91
0.25 1.0 83.28 17.09
0.5 1.0 98.04 5.14
0.75 1.0 100.00 0.00
1.0 1.0 100.00 0.00
    };
    \end{axis}
    \end{tikzpicture}
        \label{fig:heatmap-oo}
     \end{subfigure}
     \hfill
     \begin{subfigure}{0.55\textwidth}
        \caption{Under a pessimistic ground truth.}
        \begin{tikzpicture}
\begin{axis}[
    colormap/blackwhite,
    colorbar,
    xlabel = {Ground truth $(pessimistic, \rho_2)$},
    ylabel = {Assumption $(optimistic, \rho_1)$},
    point meta min=0,
    point meta max=100,
    xmin=-0.125, xmax=1.125, xscale=0.8,
    xtick style={color=black},
    xtick={0,0.25,0.5,0.75,1},
    ymin=-0.125, ymax=1.125, yscale=0.8,
    ytick style={color=black},
    ytick={0,0.25,0.5,0.75,1},
    every node near coord/.append style={font=\footnotesize, color=black, anchor=center,
        /pgf/number format/.cd, fixed, fixed zerofill, precision=2
    }
    ]
    \addplot[%
        matrix plot*,
        point meta=explicit,
        mesh/cols=5,
        visualization depends on=\thisrow{Deviation}\as\std,,
        nodes near coords={\pgfmathprintnumber{\pgfplotspointmeta}},
    ] table[meta=Mean] {
    x y Mean Deviation
0.0 0.0 81.93 20.63
0.25 0.0 81.80 14.16
0.5 0.0 75.52 17.78
0.75 0.0 74.75 18.31
1.0 0.0 74.73 18.30
0.0 0.25 33.41 28.84
0.25 0.25 85.50 16.43
0.5 0.25 79.77 15.38
0.75 0.25 78.88 15.48
1.0 0.25 78.87 15.48
0.0 0.5 30.80 27.44
0.25 0.5 72.77 16.03
0.5 0.5 86.62 16.15
0.75 0.5 85.92 15.85
1.0 0.5 85.91 15.86
0.0 0.75 32.23 28.34
0.25 0.75 73.82 16.04
0.5 0.75 85.58 15.81
0.75 0.75 86.73 16.13
1.0 0.75 86.72 16.13
0.0 1.0 32.49 28.34
0.25 1.0 73.63 16.13
0.5 1.0 85.48 15.65
0.75 1.0 86.48 16.00
1.0 1.0 86.48 16.01
    };
    \end{axis}
    \end{tikzpicture}
        \label{fig:heatmap-op}
    \end{subfigure}
\end{figure}
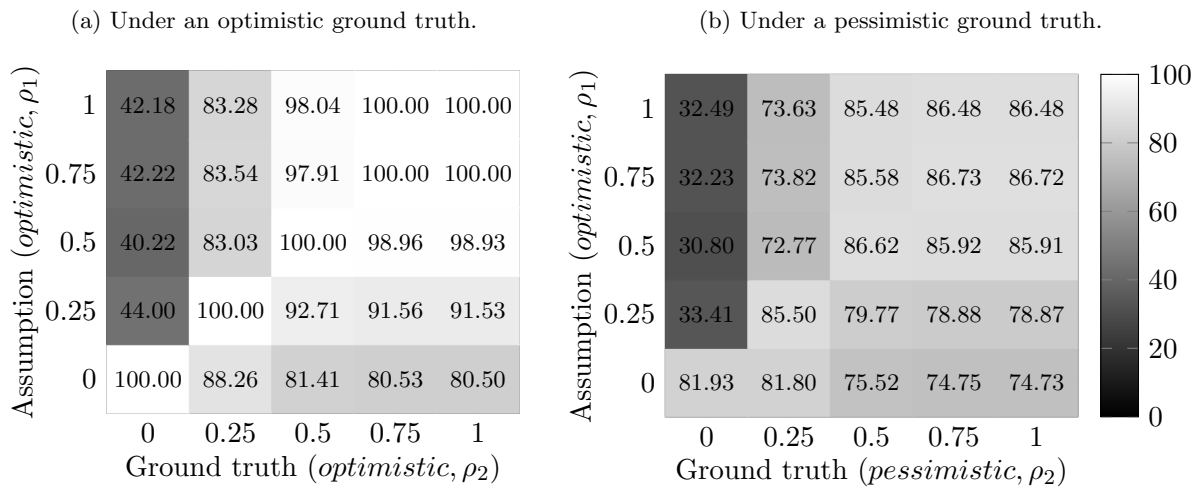

We conclude by conducting a sensitivity analysis on the assumptions made about the follower within the \ourproblem{}.
Naturally, a location schedule devised under an optimistic follower behaviour may underperform against a pessimistic one, and vice-versa.
Similarly, a wrong assumption on the revenue splitting factor $\rho$ may render the location schedules highly ineffective.
We therefore study how location schedules devised under a certain assumption perform under a different ground truth.
%
%
More specifically, let $(\boldsymbol{y}^{\star}_{b, \rho}, \boldsymbol{z}^{\star}_{b_1, \rho_1})$ and $(\boldsymbol{y}^{\star}_{b, \rho}, \boldsymbol{z}^{\star}_{b_2, \rho_2})$ be optimal solutions obtained  under some assumption $b_1$ and $\rho_1$ and under some ground truth $b_2$ and $\rho_2$, where $b_1, b_2 \in \mathset{optimistic, pessimistic}$ and $\rho_1, \rho_2 \in \mathset{0, 0.25, 0.5, 0.75, 1}$.
We compute the \textit{resilience ratio}, i.e., the relative revenue loss due to wrong parameter assumptions, as $100 \cdot \frac{\pi^{L}(\boldsymbol{y}^{\star}_{b_1, \rho_1}, \boldsymbol{z}^{\star}_{b_1, \rho_1})}{\pi^{L}(\boldsymbol{y}^{\star}_{b_2, \rho_2}, \boldsymbol{z}^{\star}_{b_2, \rho_2})}$ for each instance, where $\pi^{L}$ is the objective function of the instance under the ground truth $b_2$ and $\rho_2$ \footnote{Small (respectively, large) resilience ratios indicate that the optimal solution devised under the assumption has a weak (respectively, strong) performance when compared to what could be obtained with the ground truth.}.
Figure~\ref{fig:heatmaps-o} 
summarizes the resilience ratios for optimistic and pessimistic ground truths, averaged over the computational benchmark.
%
%
%


Figure~\ref{fig:heatmap-oo} indicates that assuming a splitting factor $\rho_1 > 0$ is bound to provide low-quality solutions to the leader whenever the ground-truth splitting factor $\rho_2 = 0$, while solutions derived with $\rho_1 = 0$ remain of moderate quality against $\rho_2 > 0$.
As a result, if the leader is even slightly suspicious that $\rho_2 = 0$, then assuming $\rho_1 = 0$ (i.e., revenue is split in favour of the follower) is a reasonably priced hedge against being too vulnerable to the follower.
However, if the leader is sufficiently confident that $\rho_2$ is at least $0.25$, then guessing the wrong value of $\rho_1 > 0$ seems less problematic, since resilience ratios remain above $80\%$, on average, for $\rho_2 \in \mathset{0.25, 0.5, 0.75, 1}$.
Figure~\ref{fig:heatmap-oo} reveals a similar pattern, the main difference being that resilience ratios are overall smaller because the follower behaves pessimistically rather than optimistically.
We can see on the diagonal that, if we estimate the splitting factor correctly (\ie, $\rho_1 = \rho_2$), then the pessimistic behaviour results in approximately $15\% -20\%$ less profit than what the leader could have obtained.
These results provide a guideline for companies facing the \ourproblem{}, particularly on which assumptions are generally more robust than others when modelling the follower.

We have found similar conclusions for the pessimistic variant of the \ourproblem{}.
We here highlight that the opportunity gap for instances with splitting factor $\rho = 0$ is always $100\%$, as the follower can simply replicate the location schedule of the leader and steal her profit completely.
We also found that solutions devised under a pessimistic assumption perform better under an optimistic ground truth, on average, than solutions devised under an optimistic assumption perform under a pessimistic ground truth.
As a result, if the leader is unsure about the expected behaviour of the follower, assuming a pessimistic one is likely to devise more robust location schedules.
Details for the pessimistic variant of the \ourproblem{} can be found in Appendix~\ref{apx:managerial-pessimistic}.

\section{Conclusion}
\label{sec:conclusion}

This paper introduces a novel competitive facility location problem, capturing salient features of emerging applications such as temporary retail.
While existing literature has been restricted to competitive location problems with a single time period, here we assume that the leader makes location decisions over a planning horizon with multiple time periods.
Each player maximizing their own profit, the follower reacts to the leader's decisions to compete over customer demand.
We here assume that unmet demand carries over from one time period to the next, a feature that is highly relevant in the context of temporary retail, for example, but has  only been sparsely addressed in the literature. 
Bridging the gap between dynamic and competitive facility location, the resulting non min-max, non zero-sum competitive planning problem requires an adapted planning approach to a delicate trade-off: letting customer demand accumulate sufficiently long, but capturing the demand before the competitor does.

%
%

On the theoretical side, we establish that the optimistic variant of our problem is $\Sigma^{p}_{2}$-hard, thereby extending the frontier of known hardness results for competitive location problems beyond zero-sum interdiction models and clarifying the inherent algorithmic difficulty induced by Stackelberg competition     and cumulative demand. 
On the methodological front, we propose bilevel mixed-integer programming formulations for the optimistic and pesssimistic variants of our problem, and present branch-and-cut algorithms to solve them.
In particular, we devise two value-function cuts: general value-function cuts tailored to our formulation, as well as their tight counterparts. 

Computational experiments on a large benchmark inspired by the Canadian province of Quebec show that the proposed branch-and-cut algorithms  strongly outperform the state-of-the-art general-purpose bilevel solver, that the tightened value-function cut improves computing time and solution quality relative to the tailored cuts. 
The managerial study revealed that ignoring competition may halve the leader’s profit, that cooperation between firms can substantially increase joint profit (by ~$6\%$), and that competition enhances service availability while cooperation broadens service coverage.
%
In addition, modelling assumptions regarding optimistic versus pessimistic follower behaviour and brand-based revenue splitting have a pronounced impact on the robustness of location schedules and their profitability. When in doubt, our results suggest that it is most beneficial to assume a pessimistic follower behaviour, and that revenue will be split in favour of the follower.

%
The combination of theoretical results, solution methods and managerial insights offers a structured approach for tackling competitive planning problems in practice.
Future research directions include analyzing simultaneous decision-making between the company and its competitor (Nash equilibria) and designing fair profit-sharing mechanisms under cooperation. Finally, the tightening technique proposed here for the tailored general value-function cuts can be applied in other problem settings. Given their effectiveness, this constitutes a promising research direction for competitive planning problems in general.


\section*{Acknowledgements} \label{sec:acknowledgments}

This work was funded by the FRQNT Doctoral Scholarship No. B2X-328911,  and the NSERC Grants No. 2017-05224 and 2024-04051. This research was also enabled in part by support provided by \href{https://www.calculquebec.ca}{Calcul Québec} and the \href{https://alliancecan.ca/}{Digital Research Alliance of Canada}.

\bibliography{references}

\appendix

\section{Mathematical Proofs}\label{apx:mathematical}

In this appendix, we prove the theoretical results presented throughout the paper.

\subsection{Hardness Proof} \label{apx:hardness}

We define the $\Sigma^{p}_{2}$ class, and then present the proof of Theorem~\ref{thm:sigmap2}.
For the sake of conciseness, we refer to the optimistic variant of the \ourproblem{} simply as \ourproblem{} throughout the proof.

A decision problem is in $\Sigma^{p}_{2}$ if we can write it in the form $\exists \boldsymbol{x}^{1} \in \mathcal{X}^{1}, \forall \boldsymbol{x}^{2} \in \mathcal{X}^{2}: P(\boldsymbol{x}^{1}, \boldsymbol{x}^{2})$, where $P(\boldsymbol{x}^{1}, \boldsymbol{x}^{2})$ is a logical predicate verifiable in polynomial time.
The decision version of the \ourproblem{}, presented in Section~\ref{sec:tractability}, is clearly in $\Sigma^{p}_{2}$, as the predicate $\pi^{L}(\boldsymbol{y}, \boldsymbol{z}) \geq \Pi\land \pi^{F}(\boldsymbol{y}, \boldsymbol{z}) \geq \pi^{F}(\boldsymbol{y}, \boldsymbol{z}^{\prime})$ can be verified in polynomial time.
We now show that the  \ourproblem{} is $\Sigma^{p}_{2}$-hard through a reduction from the $\exists \forall$3SAT, which is $\Sigma^{p}_{2}$-complete \citep{wrathall1976complete}.

\vspace{0.3cm}
\boxxx{
{\bf $\exists \forall$3SAT}: \\
{\sc instance}: Two disjoint non-empty sets of Boolean variables $\mathcal{X}^{1} = \mathset{x_{1}, \ldots, x_{k}}$ and $\mathcal{X}^{2} = \mathset{x_{k+1}, \ldots, x_{n}}$, a Boolean formula $E(x_{1}, \ldots, x_{n})$ in 3 Conjunctive Normal Form (3CNF) with exactly $m$ clauses. \\
{\sc question}: $\exists x_{1}, \ldots, x_{k}, \forall x_{k+1}, \ldots, x_{n}: E(x_{1}, \ldots, x_{n})$ is not satisfied?
}
\vspace{0.3cm}

The proof consists of two directions. If a YES instance of the \ourproblem{} implies a YES instance of the $\exists\forall$3SAT (positive direction) and a NO instance of the \ourproblem{} implies a NO instance of the $\exists\forall$3SAT (negative direction), then the \ourproblem{} is $\Sigma^{p}_{2}$-hard -- otherwise, we could answer the decision question of the $\exists\forall$3SAT through the \ourproblem{}.

\paragraph{Reduction.} 
Consider the following \ourproblem{} instance built from a $\exists\forall$3SAT instance:
\begin{itemize}
    \item We create $n$ periods.
    Each period represents the choice of a literal to a variable.
    \item We create $4n$ locations.
    Each location represents the assignment of a literal to a variable in the upper or lower level of the \ourproblem{} (\eg, $[x^{upper}_{1} = true]$, $[x^{lower}_{1} = true]$, $[x^{upper}_{1} = false]$, $[x^{lower}_{1} = false]$ are four different locations).
    We set identical rewards $r_{i} = 1, \forall i \in \mathcal{I}$.
    \item We create $4n + m$ customers.
    The first $2n$ represent variable assignments in the upper or lower level of the \ourproblem{}, regardless of their literal; the next $2n$ represent variable assignments to true or false; regardless of their level in the \ourproblem{}; and the last $m$ represent clauses.
    \begin{itemize}
        \item \textbf{Literal-customers.} For the literal-customer $[x_{s} = true]$ (respectively, $[x_{s} = false]$) such that $x_{s} \in \mathcal{X}^{1}$, we build customer rankings as $[x^{upper}_{s} = true] \succ [x^{lower}_{s} = false] \succ 0$ (respectively, $[x^{upper}_{s} = false] \succ [x^{lower}_{s} = true] \succ 0$). For the literal-customer $[x_{s} = true]$ (respectively, $[x_{s} = false]$) such that $x_{s} \in \mathcal{X}^{2}$, we build customer rankings as $[x^{lower}_{s} = true] \succ [x^{upper}_{s} = false] \succ 0$ (respectively, $[x^{lower}_{s} = false] \succ [x^{upper}_{s} = true] \succ 0$). In both cases, we set spawning demands as $d^{t}_{j} =  \begin{cases}
            M^1, \text{ if } t = s, \\
            0, \text{ otherwise }
        \end{cases} \hspace{-0.4cm} \forall t \in \mathcal{T}$, where constant $M^1$ is defined below.
        \item \textbf{Level-customers.} For the level-customer $x^{upper}_{s}$, we build customer rankings as $[x^{upper}_{s} = true] \succ [x^{upper}_{s} = false] \succ 0$ and set spawning demands as $d^{t}_{j} = \begin{cases}
        M^2, \text{ if } t = s, \\
        0, \text{ otherwise }
        \end{cases} \hspace{-0.4cm} \forall t \in \mathcal{T}$.
        For the level-customer $x^{lower}_{s}$, we build customer rankings as $[x^{lower}_{s} = true] \succ [x^{lower}_{s} = false] \succ 0$ and set spawning demands as $d^{t}_{j} = \begin{cases}
        M^3, \text{ if } t = s, \\
        0, \text{ otherwise }
        \end{cases} \hspace{-0.4cm} \forall t \in \mathcal{T}$, where constants $M^2$ and $M^3$ are defined below
        \item \textbf{Clause-customers.} For a clause of the form $(x_{s_{1}} \lor x_{s_{2}} \lor x_{s_{3}})$, we build customer rankings as $[x^{lower}_{s_{1}} = true] \succ [x^{lower}_{s_{2}} = true] \succ [x^{lower}_{s_{3}} = true] \succ [x^{upper}_{s_{3}} = false] \succ 0$. If a clause contains $\neg{x_{s}}$ rather than $x_{s}$, we can change locations in the customer ranking from $[x^{lower}_{s} = true]$ to $[x^{lower}_{s} = false]$ accordingly. We then set spawning demands as $d^{t}_{j} = \begin{cases} 1, \text{ if } t = 1 \\ 0, \text{ otherwise }\end{cases} \hspace{-0.4cm} \forall t \in \mathcal{T}$.
    \end{itemize}
    \item We set $\rho = 0.5$, $M^1 = 3m + 3$, $M^2 = 2m + 2$, $M^3 = m + 1$, where $m$ is the number of clauses.
    \item We define $\mathcal{Y}$ and $\mathcal{Z}$ as in Section~\ref{sec:components} (\ie, each player places a single facility per period).
\end{itemize}

\begin{figure}[!ht]
    \centering
    \caption{Representation of clause $x_{s_{1}} \lor \neg x_{s_{2}} \lor x_{s_{3}}$ in the \ourproblem{} instance.}
    \begin{tikzpicture}[
    customer/.style={rectangle, draw=blue!60, fill=blue!5, rounded corners, very thick},
    location/.style={rectangle, draw=gray!60, fill=gray!5, very thick},
    ]
    \node[location] (s1) {$x^{lower}_{s_{1}} = true$};
    \node[location] (s2) [below = of s1]{$x^{lower}_{s_{2}} = false$};
    \node[location] (s3) [right =7cm of s1] {$x^{lower}_{s_{3}} = true$};
    \node[location] (s4) [below = of s3] {$x^{upper}_{s_{3}} = false$};    
    \node[customer] (clause) [right =2cm of s1] {$x_{s_{1}} \lor \neg x_{s_{2}} \lor x_{s_{3}}$};
    \draw[blue] (clause.west) -- (s1.east) node [midway, fill = white] {\#1};
    \draw[blue] (clause.west) -- (s2.east) node [midway, fill = white] {\#2};
    \draw[blue] (clause.east) -- (s3.west) node [midway, fill = white] {\#3};
    \draw[blue] (clause.east) -- (s4.west) node [midway, fill = white] {\#4};
    \end{tikzpicture}
    \label{fig:gadget-clause}
\end{figure}
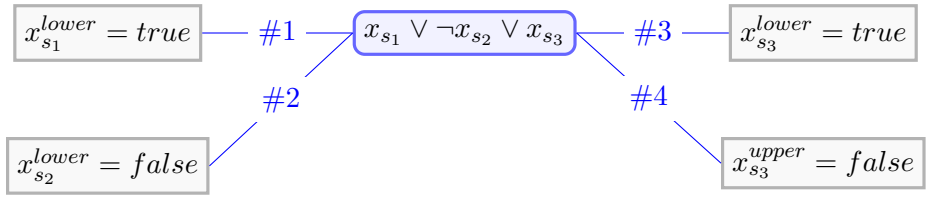

Figures~\ref{fig:gadget-exists}, \ref{fig:gadget-forall}, and \ref{fig:gadget-clause} present customer preferences associated with variables $x_{s} \in \mathcal{X}^{1}$, variables $x_{s} \in \mathcal{X}^{2}$, and clauses from the $\exists\forall$3SAT instance in  the \ourproblem{} instance.
Rectangles with sharp corners indicate locations, whereas those with round corners indicate customers.
Edges between a location and customer contain a number indicating the ranking of the former to the latter.

\begin{figure}[!ht]
    \centering
    \caption{Representation of variable $x_{s} \in \mathcal{X}^{1}$ in the \ourproblem{} instance.}
    \begin{tikzpicture}[
    cupper/.style={rectangle, draw=purple!60, fill=purple!5, rounded corners, very thick},
    clower/.style={rectangle, draw=olive!60, fill=olive!5, rounded corners, very thick},
    ctrue/.style={rectangle, draw=teal!60, fill=teal!5, rounded corners, very thick},
    cfalse/.style={rectangle, draw=orange!60, fill=orange!5, rounded corners, very thick},
    location/.style={rectangle, draw=gray!60, fill=gray!5, very thick},
    ]
    \node[location] (uppertrue) {$x^{upper}_{s} = true$};
    \node[location] (upperfalse) [below=1cm of uppertrue] {$x^{upper}_{s} = false$};
    \node[location] (lowertrue) [below=1cm of upperfalse] {$x^{lower}_{s} = true$};
    \node[location] (lowerfalse) [below=1cm of lowertrue] {$x^{lower}_{s} = false$};
    \node[ctrue] (true) [left=3cm of upperfalse] {$x_{s} = true$};
    \node[cfalse] (false) [below=1cm of true] {$x_{s} = false$};
    \node[cupper] (upper) [right=3cm of upperfalse]{$x^{upper}_{s}$};
    \node[clower] (lower) [below=1cm of upper]{$x^{lower}_{s}$};
    \draw[teal] (uppertrue.west) -- (true.east) node [midway, fill = white] {\#1};
    \draw[teal] (lowerfalse.west) .. controls +(left:60mm) .. (true.west) node [midway, fill = white] {\#2};
    \draw[orange] (upperfalse.west) -- (false.east) node [midway, fill = white] {\#1};
    \draw[orange] (lowertrue.west) -- (false.east) node [midway, fill = white] {\#2};
    \draw[purple] (uppertrue.east) -- (upper.west) node [midway, fill = white] {\#1};
    \draw[purple] (upperfalse.east) -- (upper.west) node [midway, fill = white] {\#2};
    \draw[olive] (lowertrue.east) -- (lower.west) node [midway, fill = white] {\#1};
    \draw[olive] (lowerfalse.east) -- (lower.west) node [midway, fill = white] {\#2};
    \end{tikzpicture}
    \label{fig:gadget-exists}
\end{figure}
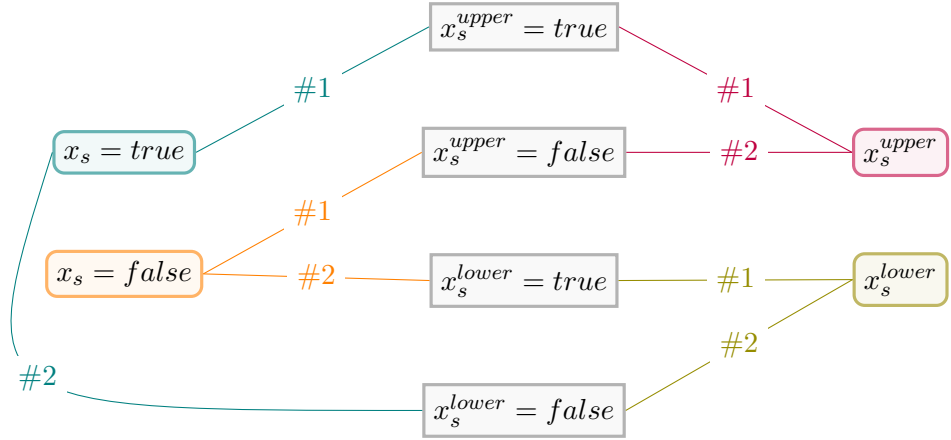

\begin{figure}[!ht]
    \centering
    \caption{Representation of variable $x_{s} \in \mathcal{X}^{2}$ in the \ourproblem{} instance.}
    \begin{tikzpicture}[
    cupper/.style={rectangle, draw=purple!60, fill=purple!5, rounded corners, very thick},
    clower/.style={rectangle, draw=olive!60, fill=olive!5, rounded corners, very thick},
    ctrue/.style={rectangle, draw=teal!60, fill=teal!5, rounded corners, very thick},
    cfalse/.style={rectangle, draw=orange!60, fill=orange!5, rounded corners, very thick},
    location/.style={rectangle, draw=gray!60, fill=gray!5, very thick},
    ]
    \node[location] (uppertrue) {$x^{upper}_{s} = true$};
    \node[location] (upperfalse) [below=1cm of uppertrue] {$x^{upper}_{s} = false$};
    \node[location] (lowertrue) [below=1cm of upperfalse] {$x^{lower}_{s} = true$};
    \node[location] (lowerfalse) [below=1cm of lowertrue] {$x^{lower}_{s} = false$};
    \node[ctrue] (true) [left=3cm of upperfalse] {$x_{s} = true$};
    \node[cfalse] (false) [below=1cm of true] {$x_{s} = false$};
    \node[cupper] (upper) [right=3cm of upperfalse]{$x^{upper}_{s}$};
    \node[clower] (lower) [below=1cm of upper]{$x^{lower}_{s}$};
    \draw[teal] (lowertrue.west) -- (true.east) node [midway, fill = white] {\#1};
    \draw[teal] (upperfalse.west) -- (true.east) node [midway, fill = white] {\#2};
    \draw[orange] (lowerfalse.west) -- (false.east) node [midway, fill = white] {\#1};
    \draw[orange] (uppertrue.west) .. controls +(left:60mm) .. (false.west) node [midway, fill = white] {\#2};
    \draw[purple] (uppertrue.east) -- (upper.west) node [midway, fill = white] {\#1};
    \draw[purple] (upperfalse.east) -- (upper.west) node [midway, fill = white] {\#2};
    \draw[olive] (lowertrue.east) -- (lower.west) node [midway, fill = white] {\#1};
    \draw[olive] (lowerfalse.east) -- (lower.west) node [midway, fill = white] {\#2};
    \end{tikzpicture}
    \label{fig:gadget-forall}
\end{figure}
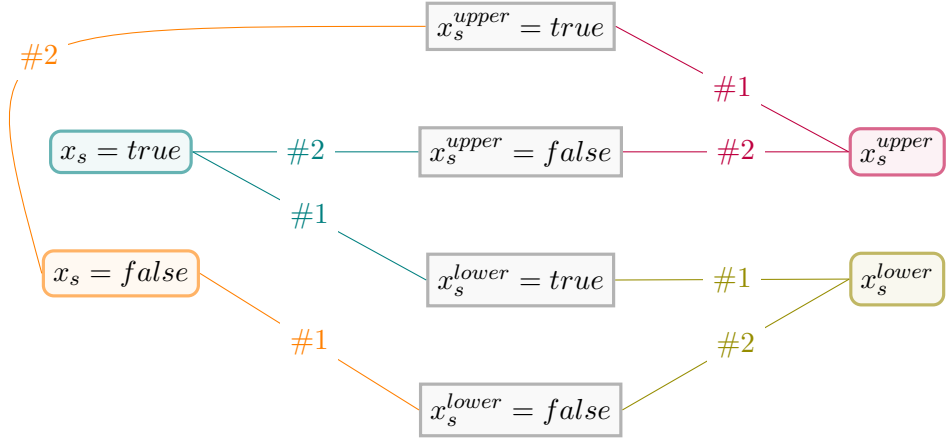

We now highlight some properties of this \ourproblem{} instance.
First, note that the leader chooses between locations $[x^{upper}_{s} = true]$ and $[x^{upper}_{s} = false]$ at period $s$.
In fact, if the leader were to choose between locations $[x^{lower}_{s} = true]$ and $[x^{lower}_{s} = false]$ at period $s$, she would swap a guaranteed marginal contribution between $M^1 + M^2 = 5m + 5$ and $M^1 + M^2 + m = 6m + 5$ for a marginal contribution between $M^1 + M^3 = 4m + 4$ and $M^1 + M^3 + m = 5m + 4$, which is suboptimal.
Similarly, the follower chooses between locations $[x^{lower}_{s} = true]$ and $[x^{lower}_{s} = false]$ at period $s$.
Although $\rho = 0.5$, if the follower were to choose between locations $[x^{upper}_{s} = true]$ and $[x^{upper}_{s} = false]$ at period $s$, he would swap a guaranteed marginal contribution between $M^1 + M^3 = 4m + 4$ and $M^1 + M^3 + m = 5m + 4$ for a marginal contribution between $\frac{M^1 + M^2}{2} = \frac{5m + 5}{2}$ and $\frac{M^1 + M^2 + m}{2} = \frac{6m + 5}{2}$, which is suboptimal.

Second, note that the leader has the priority of assigning a literal to variables $x_{s} \in \mathcal{X}^{1}$.
If the leader chooses $[x^{upper}_{s} = true]$ (respectively, $[x^{upper}_{s} = false]$), the follower chooses $[x^{lower}_{s} = true]$ (respectively, $[x^{lower}_{s} = false]$), otherwise he obtains a suboptimal marginal contribution from splitting customer $x_s = true$ (respectively, $x_s = false$).
In turn, the follower has the priority of assigning a literal to variables $x_{s} \in \mathcal{X}^{2}$.
If the leader anticipates the follower to choose $[x^{lower}_{s} = true]$ (respectively, $[x^{lower}_{s} = false]$), she chooses $[x^{upper}_{s} = true]$ (respectively, $[x^{uper}_{s} = false]$), otherwise she obtains a suboptimal marginal contribution from splitting customer $x_s = true$ (respectively, $x_s = false$).
On the one hand, the leader assigns literals to variables $x_{s} \in \mathcal{X}^{1}$ to prevent one or more clauses from being satisfiable with variables $x_{s} \in \mathcal{X}^{2}$ (\ie, prevent one or more clause-customers from being captured by the follower).
On the other hand, the follower assigns literals to variables $x_{s} \in \mathcal{X}^{2}$ to satisfy as many clauses as possible given literals of variables $x_{s} \in \mathcal{X}^{1}$ (\ie, capture as many clause-customers as possible).

If the profit of the leader is greater than or equal to $n (M^1 + M^2) + 1 = n (5m + 5) + 1$, then there is at least one clause in the expression $E(x_{1}, \ldots, x_{n})$ that has not been satisfied -- otherwise, the associated clause-customer would have been captured by the follower.
We therefore look for an objective value $\Pi$ of the \ourproblem{} greater or equal to $\Pi = n (M^1 + M^2) + 1$.

\paragraph{Positive Direction.}
Assume a YES instance of the \ourproblem{} (\ie, $\exists (\hat{\boldsymbol{y}}, \hat{\boldsymbol{z}}) \in \mathcal{Y} \times \mathcal{Z}, \forall \boldsymbol{z}^{\prime} \in \mathcal{Z}: \pi^{L}(\hat{\boldsymbol{y}}, \hat{\boldsymbol{z}}) \geq n(M^1 + M^2) + 1 \land \pi^{F}(\hat{\boldsymbol{y}}, \hat{\boldsymbol{z}}) \geq \pi^{F}(\hat{\boldsymbol{y}}, \boldsymbol{z}^{\prime})$).
Since $\pi^{L}(\hat{\boldsymbol{y}}, \hat{\boldsymbol{z}}) \geq n(M^1 + M^2) + 1$ and $\pi^{F}(\hat{\boldsymbol{y}}, \hat{\boldsymbol{z}}) \geq \pi^{F}(\hat{\boldsymbol{y}}, \boldsymbol{z}^{\prime})$, there is no optimal reaction that could have been chosen by the follower to capture all clause-customers, so at least one clause is not satisfied in the expression $E(x_{1}, \ldots, x_{n})$. 
This implies a YES instance of the $\exists\forall$3SAT, so the positive direction holds.

\paragraph{Negative Direction.}
Assume a NO instance of the \ourproblem{} (\ie, $\forall (\hat{\boldsymbol{y}}, \hat{\boldsymbol{z}}) \in \mathcal{Y} \times \mathcal{Z}, \exists \boldsymbol{z}^{\prime} \in \mathcal{Z}: \pi^{L}(\hat{\boldsymbol{y}}, \hat{\boldsymbol{z}}) < n(M^1 + M^2) + 1 \lor \pi^{F}(\hat{\boldsymbol{y}}, \hat{\boldsymbol{z}}) < \pi^{F}(\hat{\boldsymbol{y}}, \boldsymbol{z}^{\prime})$).
If we fall into the second case (\ie, $\pi^{F}(\hat{\boldsymbol{y}}, \hat{\boldsymbol{z}}) < \pi^{F}(\hat{\boldsymbol{y}}, \boldsymbol{z}^{\prime})$), we can recursively look into the predicate for pair $(\hat{\boldsymbol{y}}, \boldsymbol{z}^{\prime})$ rather than pair $(\hat{\boldsymbol{y}}, \hat{\boldsymbol{z}})$.
This reasoning allow us to rewrite the statement as $\forall (\hat{\boldsymbol{y}}, \hat{\boldsymbol{z}}) \in \mathcal{Y} \times \mathcal{Z}: \pi^{L}(\hat{\boldsymbol{y}}, \hat{\boldsymbol{z}}) \leq n(M^1 + M^2)$.
Since $\pi^{L}(\hat{\boldsymbol{y}}, \hat{\boldsymbol{z}}) \leq n(M^1 + M^2)$, it is not possible for the leader to capture clause-customers, so all clauses are satisfied in expression $E(x_{1}, \ldots, x_{n})$.  This implies a NO instance of the $\exists\forall$3SAT, so the negative direction holds and the \ourproblem{} is $\Sigma^{p}_{2}$-hard. \Halmos

\subsection{Branch-and-Cut Proofs} \label{apx:branch-n-cut}

We prove Proposition~\ref{prp:termination}, Lemmas~\ref{lem:tailored}--\ref{lem:improved} and Theorem~\ref{thm:dominance} in Appendices~\ref{apx:termination}--\ref{apx:dominance}, respectively.

\subsubsection{Proof of Proposition~\ref{prp:termination}} \label{apx:termination}

There are $|\mathcal{Y} \times \mathcal{Z}|$ integer-feasible solutions for each instance of the \ourproblem{}, which is an exponentially large but finite number.
Once we visit an integer-feasible solution $(\boldsymbol{y}^{\prime}, \boldsymbol{z}^{\prime})$, we either \textit{(i)} accept it as bilevel feasible and stop descending on that node or \textit{(ii)} deem it bilevel infeasible and cut it from that node onwards.
We therefore never visit the same integer-feasible solution $(\boldsymbol{y}^{\prime}, \boldsymbol{z}^{\prime})$ more than once, and the branch-and-cut framework with Algorithm~\ref{alg:callback-optimistic} finishes in a finite number of steps.
Assume now, for the sake of contradiction, that this procedure does not output the optimal solution to the optimistic variant of the \ourproblem{}.
This means that the optimal bilevel solution $(\boldsymbol{y}^{\star}, \boldsymbol{z}^{\star})$ has been removed by the value-function cut generated for a bilevel infeasible solution $(\boldsymbol{y}^{\prime}, \boldsymbol{z}^{\prime})$, which is an absurd due to Lemmas~\ref{lem:tailored} and \ref{lem:improved}.

\subsubsection{Proof of Lemma~\ref{lem:tailored}} \label{apx:baseline}
First, we show that Inequality~\eqref{eq:tailored-cut}, built with a bilevel feasible solution $(\boldsymbol{y}^{\prime}, \boldsymbol{z}^{\star})$, cuts any bilevel-infeasible but integer-feasible solution of the form $(\boldsymbol{y}^{\prime}, \boldsymbol{z}^{\prime})$; in other words, this cut eliminates any non-optimal follower feasible reaction for $\boldsymbol{y}^\prime$.
Assume, for the sake of contradiction, that this is false.
This means that optimal values $v^{\ell{}t}_{ij} (\boldsymbol{y}^{\prime}, \boldsymbol{z}^{\prime})$ satisfy Inequality~\eqref{eq:tailored-cut}.
If $v^{\ell{}t}_{ij} (\boldsymbol{y}^{\prime}, \boldsymbol{z}^{\star}) > 0$, then $\sum_{\substack{k \in \mathcal{I}: \\ k \succ_{j} 0}} \sum_{\substack{s \in \mathcal{T}: \\ \ell{} < s < t}} {y^{\prime}}^{s}_{k} + \sum_{\substack{k \in \mathcal{I}: \\ k \succ_{j} i}} {y^{\prime}}^{t}_{k} = 0$, otherwise Constraints~\eqref{eq:third-detailed-optimistic-ct1}--\eqref{eq:third-detailed-optimistic-ct10} would not be respected. 
Due to the same argument, if $v^{\ell{}t}_{ij} (\boldsymbol{y}^{\prime}, \boldsymbol{z}^{\star}) = 1 \land \rho > 0$, then ${y^{\prime}}^{t}_{i} = 0$.
This leads to
$$
\sum_{t \in \mathcal{T}} \sum_{\substack{\ell{} \in \mathcal{T}^{S}: \\ \ell{} < t}} \sum_{j \in \mathcal{J}} \sum_{\substack{i \in \mathcal{I}: \\ i \succ_{j} 0}} r_{i} D^{\ell{}t}_{j}  v^{\ell{}t}_{ij} (\boldsymbol{y}^{\prime}, \boldsymbol{z}^{\prime}) \geq
\sum_{t \in \mathcal{T}} \sum_{\substack{\ell{} \in \mathcal{T}^{S}: \\ \ell{} < t}} \sum_{j \in \mathcal{J}} \sum_{\substack{i \in \mathcal{I}: \\ i \succ_{j} 0}} r_{i} D^{\ell{}t}_{j} v^{\ell{}t}_{ij} (\boldsymbol{y}^{\prime}, \boldsymbol{z}^{\star}),
$$
which is a contradiction, because solution $(\boldsymbol{y}^{\prime}, \boldsymbol{z}^{\prime})$ is not bilevel feasible.

We now show that Inequality~\eqref{eq:tailored-cut}, built with a bilevel feasible solution $(\boldsymbol{y}^{\prime}, \boldsymbol{z}^{\star})$, holds for all bilevel feasible solutions.
Assume, for the sake of contradiction, that this is false.
This implies that there is a bilevel feasible solution $(\boldsymbol{y}^{\prime\prime}, \boldsymbol{z}^{\prime\prime})$ such that
\begin{align*}
    \sum_{t \in \mathcal{T}} \sum_{\substack{\ell{} \in \mathcal{T}^{S}: \\ \ell{} < t}} \sum_{j \in \mathcal{J}} \sum_{\substack{i \in \mathcal{I}: \\ i \succ_{j} 0}} r_{i} D^{\ell{}t}_{j} v^{\ell{}t}_{ij} (\boldsymbol{y}^{\prime\prime}, \boldsymbol{z}^{\prime\prime}) < \\
        \sum_{t \in \mathcal{T}} \sum_{\substack{\ell{} \in \mathcal{T}^{S}: \\ \ell{} < t}} \sum_{j \in \mathcal{J}} \sum_{\substack{i \in \mathcal{I}: \\ i \succ_{j} 0}} r_{i} D^{\ell{}t}_{j} v^{\ell{}t}_{ij} (\boldsymbol{y}^{\prime}, \boldsymbol{z}^{\star}) \left( 1 - \sum_{\substack{k \in \mathcal{I}: \\ k \succ_{j} 0}} \sum_{\substack{s \in \mathcal{T}: \\ \ell{} < s < t}} {y^{\prime\prime}}^{s}_{k} - \sum_{\substack{k \in \mathcal{I}: \\ k \succ_{j} i}} {y^{\prime\prime}}^{t}_{k} - \mathbb{I} \left[ v^{\ell{}t}_{ij} (\boldsymbol{y}^{\prime}, \boldsymbol{z}^{\star}) = 1 \land \rho > 0 \right] {y^{\prime\prime}}^{t}_{i} \right).
\end{align*}
Since $(\boldsymbol{y}^{\prime\prime}, \boldsymbol{z}^{\prime\prime})$ is bilevel feasible, we can rewrite the left-hand side as follows
\begin{align*}
    \sum_{t \in \mathcal{T}} \sum_{\substack{\ell{} \in \mathcal{T}^{S}: \\ \ell{} < t}} \sum_{j \in \mathcal{J}} \sum_{\substack{i \in \mathcal{I}: \\ i \succ_{j} 0}} r_{i} D^{\ell{}t}_{j} v^{\ell{}t}_{ij} (\boldsymbol{y}^{\prime\prime}, \boldsymbol{z}^{\star}) \leq \sum_{t \in \mathcal{T}} \sum_{\substack{\ell{} \in \mathcal{T}^{S}: \\ \ell{} < t}} \sum_{j \in \mathcal{J}} \sum_{\substack{i \in \mathcal{I}: \\ i \succ_{j} 0}} r_{i} D^{\ell{}t}_{j} v^{\ell{}t}_{ij} (\boldsymbol{y}^{\prime\prime}, \boldsymbol{z}^{\prime\prime}) < \\
        \sum_{t \in \mathcal{T}} \sum_{\substack{\ell{} \in \mathcal{T}^{S}: \\ \ell{} < t}} \sum_{j \in \mathcal{J}} \sum_{\substack{i \in \mathcal{I}: \\ i \succ_{j} 0}} r_{i} D^{\ell{}t}_{j} v^{\ell{}t}_{ij} (\boldsymbol{y}^{\prime}, \boldsymbol{z}^{\star}) \left( 1 - \sum_{\substack{k \in \mathcal{I}: \\ k \succ_{j} 0}} \sum_{\substack{s \in \mathcal{T}: \\ \ell{} < s < t}} {y^{\prime\prime}}^{s}_{k} - \sum_{\substack{k \in \mathcal{I}: \\ k \succ_{j} i}} {y^{\prime\prime}}^{t}_{k} - \mathbb{I} \left[ v^{\ell{}t}_{ij} (\boldsymbol{y}^{\prime}, \boldsymbol{z}^{\star}) = 1 \land \rho > 0 \right] {y^{\prime\prime}}^{t}_{i} \right).
\end{align*}
There is therefore at least one combination of customer $j$ and period $t$ (with location $i$) such that
\begin{align*}
    \sum_{\substack{\ell{} \in \mathcal{T}^{S}: \\ \ell{} < t}} r_{i} D^{\ell{}t}_{j}v^{\ell{}t}_{ij} (\boldsymbol{y}^{\prime\prime}, \boldsymbol{z}^{\star}) < \\
        \sum_{\substack{\ell{} \in \mathcal{T}^{S}: \\ \ell{} < t}}  r_{i} D^{\ell{}t}_{j}v^{\ell{}t}_{ij} (\boldsymbol{y}^{\prime}, \boldsymbol{z}^{\star}) \left( 1 - \sum_{\substack{k \in \mathcal{I}: \\ k \succ_{j} 0}} \sum_{\substack{s \in \mathcal{T}: \\ \ell{} < s < t}} {y^{\prime\prime}}^{s}_{k} - \sum_{\substack{k \in \mathcal{I}: \\ k \succ_{j} i}} {y^{\prime\prime}}^{t}_{k} - \mathbb{I} \left[ v^{\ell{}t}_{ij} (\boldsymbol{y}^{\prime}, \boldsymbol{z}^{\star}) = 1 \land \rho > 0 \right] {y^{\prime\prime}}^{t}_{i} \right).
\end{align*}
We abuse the notation to simply write location $i$ as the location that captures customer $j$ at period $t$ within location schedule $\boldsymbol{z}^{\star}$.
We further analyze whether this inequality can ever hold.

First, if $\sum_{\substack{\ell{} \in \mathcal{T}^{S}: \\ \ell{} < t}} r_{i} D^{\ell{}t}_{j}v^{\ell{}t}_{ij} (\boldsymbol{y}^{\prime\prime}, \boldsymbol{z}^{\star}) = 0$, then $\exists {y^{\prime\prime}}^{t}_{k} = 1: k \succ_{j} i$, so the right-hand side becomes zero (or strictly negative), and the inequality does not hold.
If $\sum_{\substack{\ell{} \in \mathcal{T}^{S}: \\ \ell{} < t}} r_{i} D^{\ell{}t}_{j}v^{\ell{}t}_{ij} (\boldsymbol{y}^{\prime\prime}, \boldsymbol{z}^{\star}) > 0$, we can rewrite the inequality with periods $\ell{}_1$ and $\ell{}_2$ as follows (for more details, see Section~\ref{sec:decomposition}):
\begin{align*}
    r_{i} D^{\ell{}_1t}_{j} v^{\ell{}_1t}_{ij} (\boldsymbol{y}^{\prime\prime}, \boldsymbol{z}^{\star}) < \\
        r_{i} D^{\ell{}_2t}_{j}v^{\ell{}_2t}_{ij} (\boldsymbol{y}^{\prime}, \boldsymbol{z}^{\star}) \left( 1 - \sum_{\substack{k \in \mathcal{I}: \\ k \succ_{j} 0}} \sum_{\substack{s \in \mathcal{T}: \\ \ell{}_2 < s < t}} {y^{\prime\prime}}^{s}_{k} - \sum_{\substack{k \in \mathcal{I}: \\ k \succ_{j} i}} {y^{\prime\prime}}^{t}_{k} - \mathbb{I} \left[ v^{\ell{}_2t}_{ij} (\boldsymbol{y}^{\prime}, \boldsymbol{z}^{\star}) = 1 \land \rho > 0 \right] {y^{\prime\prime}}^{t}_{i} \right).
\end{align*}

Consider first the case where $v^{\ell{}_1t}_{ij} (\boldsymbol{y}^{\prime\prime}, \boldsymbol{z}^{\star}) = 1 - \rho < 1$. 
If $v^{\ell{}_1t}_{ij} (\boldsymbol{y}^{\prime\prime}, \boldsymbol{z}^{\star}) = 1 - \rho$, then $\rho > 0$,  ${y^{\prime\prime}}^{t}_{i} = 1$, and $\sum_{\substack{k \in \mathcal{I}: \\ k \succ_{j} 0}} \sum_{\substack{s \in \mathcal{T}: \\ \ell{}_1 < s < t}} {y^{\prime\prime}}^{s}_{k} + \sum_{\substack{k \in \mathcal{I}: \\ k \succ_{j} i}} {y^{\prime\prime}}^{t}_{k} = 0$.
The inequality then forces $v^{\ell{}_2t}_{ij} (\boldsymbol{y}^{\prime}, \boldsymbol{z}^{\star}) > 0$.
If $v^{\ell{}_2t}_{ij} (\boldsymbol{y}^{\prime}, \boldsymbol{z}^{\star}) = 1$, then the right-hand side becomes zero (or strictly negative) due to ${y^{\prime\prime}}^{t}_{i} = 1$ and $\rho > 0$, and the inequality does not hold.
If $v^{\ell{}_2t}_{ij} (\boldsymbol{y}^{\prime}, \boldsymbol{z}^{\star}) = 1 - \rho$, then the inequality simplifies to
\begin{align}
    D^{\ell{}_1t}_{j} < 
        D^{\ell{}_2t}_{j} \left( 1 - \sum_{\substack{k \in \mathcal{I}: \\ k \succ_{j} 0}} \sum_{\substack{s \in \mathcal{T}: \\ \ell{}_2 < s \leq \ell{}_1}} {y^{\prime\prime}}^{s}_{k} \right). \label{eq:simplification}
\end{align}
If $\ell{}_2 \geq \ell{}_1$, then $\sum_{\substack{k \in \mathcal{I}: \\ k \succ_{j} 0}} \sum_{\substack{s \in \mathcal{T}: \\ \ell{}_2 < s \leq \ell{}_1}} {y^{\prime\prime}}^{s}_{k} = 0$, and the inequality simplifies to $D^{\ell{}_1t}_{j} < D^{\ell{}_2t}_{j}$, which is an absurd.
If $\ell{}_2 < \ell{}_1$, then $\exists {y^{\prime\prime}}^{s}_{k} = 1 : \ell{}_2 < s \leq \ell{}_1, k \succ_{j} 0$, so the right-hand side becomes zero (or strictly negative), and the inequality does not hold.

Consider now the case where $v^{\ell{}_1t}_{ij} (\boldsymbol{y}^{\prime\prime}, \boldsymbol{z}^{\star}) = 1$. 
If $v^{\ell{}_1t}_{ij} (\boldsymbol{y}^{\prime\prime}, \boldsymbol{z}^{\star}) = 1$, then $\rho = 0$ or $y^{t}_{i} = 0$,  and $\sum_{\substack{k \in \mathcal{I}: \\ k \succ_{j} 0}} \sum_{\substack{s \in \mathcal{T}: \\ \ell{}_1 < s < t}} {y^{\prime\prime}}^{s}_{k} + \sum_{\substack{k \in \mathcal{I}: \\ k \succ_{j} i}} {y^{\prime\prime}}^{t}_{k} = 0$.
The inequality then forces $v^{\ell{}_2t}_{ij} (\boldsymbol{y}^{\prime}, \boldsymbol{z}^{\star}) > 0$.
If $v^{\ell{}_2t}_{ij} (\boldsymbol{y}^{\prime}, \boldsymbol{z}^{\star}) = 1$, then the inequality simplifies to Inequality \eqref{eq:simplification} again, and we follow the same reasoning to arrive at a contradiction.
If $v^{\ell{}_2t}_{ij} (\boldsymbol{y}^{\prime}, \boldsymbol{z}^{\star}) = 1 - \rho$, then we can upper bound the right-hand side by $D^{\ell{}_2t}_{j} \left( 1 - \sum_{\substack{k \in \mathcal{I}: \\ k \succ_{j} 0}} \sum_{\substack{s \in \mathcal{T}: \\ \ell{}_2 < s \leq \ell{}_1}} {y^{\prime\prime}}^{s}_{k} \right)$ and arrive at Inequality \eqref{eq:simplification} again, which result in a contradiction.

Since this inequality can never hold, we have shown that the lemma holds by contradiction. \Halmos

\subsubsection{Proof of Lemma~\ref{lem:improved}} \label{apx:improved}

First, we show that Inequality~\eqref{eq:improved}, built with a bilevel feasible solution $(\boldsymbol{y}^{\prime}, \boldsymbol{z}^{\star})$, cuts any bilevel-infeasible but integer-feasible solution of the form $(\boldsymbol{y}^{\prime}, \boldsymbol{z}^{\prime})$; in other words, this cut eliminates any non-optimal follower feasible reaction for $\boldsymbol{y}^\prime$.
Assume, for the sake of contradiction, that this is false. 
This means that optimal values $v^{\ell{}t}_{ij} (\boldsymbol{y}^{\prime}, \boldsymbol{z}^{\prime})$ satisfy Inequality~\eqref{eq:improved}.
Note that, for an empty location schedule $\boldsymbol{0}$, which we employ to  build the improved cut, $o^{\ell{}st}_{ij} (\boldsymbol{0}) = 0, \forall s \in \mathcal{T} : \ell{} <s \leq t$.
This leads to
\begin{align*}
\sum_{t \in \mathcal{T}} \sum_{\substack{\ell{} \in \mathcal{T}^{S}: \\ \ell{} < t}} \sum_{j \in \mathcal{J}} \sum_{\substack{i \in \mathcal{I}: \\ i \succ_{j} 0}} r_{i} D^{\ell{}t}_{j}  v^{\ell{}t}_{ij} (\boldsymbol{y}^{\prime}, \boldsymbol{z}^{\prime}) \geq 
    \sum_{t \in \mathcal{T}} \sum_{\substack{\ell{} \in \mathcal{T}^{S}: \\ \ell{} < t}} \sum_{j \in \mathcal{J}} \sum_{\substack{i \in \mathcal{I}: \\ i \succ_{j} 0}} r_{i} D^{\ell{}t}_{j} v^{\ell{}t}_{ij} (\boldsymbol{0}, \boldsymbol{z}^{\star}) \geq \\
        \sum_{t \in \mathcal{T}} \sum_{\substack{\ell{} \in \mathcal{T}^{S}: \\ \ell{} < t}} \sum_{j \in \mathcal{J}} \sum_{\substack{i \in \mathcal{I}: \\ i \succ_{j} 0}} r_{i} D^{\ell{}t}_{j} v^{\ell{}t}_{ij} (\boldsymbol{y}^{\prime}, \boldsymbol{z}^{\star}),
\end{align*}
which is an absurd, because solution $(\boldsymbol{y}^{\prime}, \boldsymbol{z}^{\prime})$ is not bilevel feasible. 
We remark that the right-most transition holds because the profit of the follower can either remain the same value or decrease once the leader opens some location over the planning horizon rather than none.

We now show that Inequality~\eqref{eq:improved}, built with a bilevel feasible solution $(\boldsymbol{y}^{\prime}, \boldsymbol{z}^{\star})$, holds for all bilevel feasible solutions.
Assume, for the sake of contradiction, that this is false.
This implies that there is a bilevel feasible solution $(\boldsymbol{y}^{\prime\prime}, \boldsymbol{z}^{\prime\prime})$ such that
\begin{align*}
    \sum_{t \in \mathcal{T}} \sum_{\substack{\ell{} \in \mathcal{T}^{S}: \\ \ell{} < t}} \sum_{j \in \mathcal{J}} \sum_{\substack{i \in \mathcal{I}: \\ i \succ_{j} 0}} r_{i} D^{\ell{}t}_{j} v^{\ell{}t}_{ij} (\boldsymbol{y}^{\prime\prime}, \boldsymbol{z}^{\prime\prime}) < \\
        \sum_{t \in \mathcal{T}} \sum_{\substack{\ell{} \in \mathcal{T}^{S} :\\ \ell{} < t}} \sum_{j \in \mathcal{J}} \sum_{\substack{i \in \mathcal{I}: \\ i \succ_{j} 0}}  r_{i} v^{\ell{}t}_{ij} (\boldsymbol{0}, \boldsymbol{z}^{\star}) \left( D^{\ell{}t}_{j} - \sum_{\substack{s \in \mathcal{T}: \\ \ell{} < s \leq t}} d^{s}_{j} o^{\ell{}st}_{ij} (\boldsymbol{y}^{\prime\prime}) - \rho {y^{\prime\prime}}^{t}_{i} \sum_{\substack{s \in \mathcal{T}: \\ \ell{} < s \leq t}} d^{s}_{j} (1 - o^{\ell{}st}_{ij} (\boldsymbol{y}^{\prime\prime})) \right).
\end{align*}
Since $(\boldsymbol{y}^{\prime\prime}, \boldsymbol{z}^{\prime\prime})$ is bilevel feasible, we can rewrite the left-hand side as follows
\begin{align*}
    \sum_{t \in \mathcal{T}} \sum_{\substack{\ell{} \in \mathcal{T}^{S}: \\ \ell{} < t}} \sum_{j \in \mathcal{J}} \sum_{\substack{i \in \mathcal{I}: \\ i \succ_{j} 0}} r_{i} D^{\ell{}t}_{j} v^{\ell{}t}_{ij} (\boldsymbol{y}^{\prime\prime}, \boldsymbol{z}^{\star}) \leq \sum_{t \in \mathcal{T}} \sum_{\substack{\ell{} \in \mathcal{T}^{S}: \\ \ell{} < t}} \sum_{j \in \mathcal{J}} \sum_{\substack{i \in \mathcal{I}: \\ i \succ_{j} 0}} r_{i} D^{\ell{}t}_{j} v^{\ell{}t}_{ij} (\boldsymbol{y}^{\prime\prime}, \boldsymbol{z}^{\prime\prime}) < \\
        \sum_{t \in \mathcal{T}} \sum_{\substack{\ell{} \in \mathcal{T}^{S} :\\ \ell{} < t}} \sum_{j \in \mathcal{J}} \sum_{\substack{i \in \mathcal{I}: \\ i \succ_{j} 0}}  r_{i} v^{\ell{}t}_{ij} (\boldsymbol{0}, \boldsymbol{z}^{\star}) \left( D^{\ell{}t}_{j} - \sum_{\substack{s \in \mathcal{T}: \\ \ell{} < s \leq t}} d^{s}_{j} o^{\ell{}st}_{ij} (\boldsymbol{y}^{\prime\prime}) - \rho {y^{\prime\prime}}^{t}_{i} \sum_{\substack{s \in \mathcal{T}: \\ \ell{} < s \leq t}} d^{s}_{j} (1 - o^{\ell{}st}_{ij} (\boldsymbol{y}^{\prime\prime})) \right).
\end{align*}
There is therefore at least one combination of customer $j$ and period $t$ (with location $i$) such that 
\begin{align*}
    \sum_{\substack{\ell{} \in \mathcal{T}^{S}: \\ \ell{} < t}} r_{i} D^{\ell{}t}_{j} v^{\ell{}t}_{ij} (\boldsymbol{y}^{\prime\prime}, \boldsymbol{z}^{\star}) < 
        \sum_{\substack{\ell{} \in \mathcal{T}^{S}: \\ \ell{} < t}} r_{i} v^{\ell{}t}_{ij} (\boldsymbol{0}, \boldsymbol{z}^{\star}) \left( D^{\ell{}t}_{j} - \sum_{\substack{s \in \mathcal{T}: \\ \ell{} < s \leq t}} d^{s}_{j} o^{\ell{}st}_{ij} (\boldsymbol{y}^{\prime\prime}) - \rho {y^{\prime\prime}}^{t}_{i} \sum_{\substack{s \in \mathcal{T}: \\ \ell{} < s \leq t}} d^{s}_{j} (1 - o^{\ell{}st}_{ij} (\boldsymbol{y}^{\prime\prime})) \right).
\end{align*}
We abuse the notation to simply write location $i$ as the location that captures customer $j$ at period $t$ within location schedule $\boldsymbol{z}^{\star}$.
We further analyze whether this inequality can ever hold.

First, if $\sum_{\substack{\ell{} \in \mathcal{T}^{S}: \\ \ell{} < t}} r_{i} D^{\ell{}t}_{j}v^{\ell{}t}_{ij} (\boldsymbol{y}^{\prime\prime}, \boldsymbol{z}^{\star}) = 0$, then $\exists {y^{\prime\prime}}^{t}_{k} = 1: k \succ_{j} i$ and $o^{\ell{}st}_{ij} (\boldsymbol{y}^{\prime\prime}) = 1, \forall \ell{} \in \mathcal{T}^{S}: \ell{} < t, \forall s \in \mathcal{T}: \ell{} < s \leq t$, so the right-hand side becomes zero, and the inequality does not hold.
If $\sum_{\substack{\ell{} \in \mathcal{T}^{S}: \\ \ell{} < t}} r_{i} D^{\ell{}t}_{j}v^{\ell{}t}_{ij} (\boldsymbol{y}^{\prime\prime}, \boldsymbol{z}^{\star}) > 0$, we can rewrite the inequality with periods $\ell{}_1$ and $\ell{}_2$ as follows, where it must hold that $\ell{}_2 \leq \ell{}_1$ (for more details, see Section~\ref{sec:decomposition}):
\begin{align*}
    r_{i} D^{\ell{}_1t}_{j} v^{\ell{}_1t}_{ij} (\boldsymbol{y}^{\prime\prime}, \boldsymbol{z}^{\star}) < 
        r_{i} v^{\ell{}_2t}_{ij} (\boldsymbol{0}, \boldsymbol{z}^{\star}) \left( D^{\ell{}_2t}_{j} - \sum_{\substack{s \in \mathcal{T}: \\ \ell{}_2 < s \leq t}} d^{s}_{j} o^{\ell{}_2st}_{ij} (\boldsymbol{y}^{\prime\prime}) - \rho {y^{\prime\prime}}^{t}_{i} \sum_{\substack{s \in \mathcal{T}: \\ \ell{}_2 < s \leq t}} d^{s}_{j} (1 - o^{\ell{}_2st}_{ij} (\boldsymbol{y}^{\prime\prime})) \right).
\end{align*}

Consider first the case where $v^{\ell{}_1t}_{ij} (\boldsymbol{y}^{\prime\prime}, \boldsymbol{z}^{\star}) = 1 - \rho < 1$. 
If $v^{\ell{}_1t}_{ij} (\boldsymbol{y}^{\prime\prime}, \boldsymbol{z}^{\star}) = 1 - \rho$, then $\rho > 0$,  ${y^{\prime\prime}}^{t}_{i} = 1$, $o^{\ell{}_2st}_{ij} (\boldsymbol{y}^{\prime\prime}) = 1, \forall s \in \mathcal{T}: \ell{}_2 < s \leq \ell{}_1$, and $o^{\ell{}_2st}_{ij} (\boldsymbol{y}^{\prime\prime}) = 0,  \forall s \in \mathcal{T}: \ell{}_1 < s \leq t$.
The inequality then forces $v^{\ell{}_2t}_{ij} (\boldsymbol{0}, \boldsymbol{z}^{\star}) = 1$ and results in a contradiction
\begin{align}
    D^{\ell{}_1t}_{j} (1 - \rho) < 
        \left( D^{\ell{}_2t}_{j}  - \sum_{\substack{s \in \mathcal{T}: \\ \ell{}_2 < s \leq \ell{}_1}} d^{s}_{j} - \rho \sum_{\substack{s \in \mathcal{T}: \\ \ell{}_1 < s \leq t}} d^{s}_{j}  \right) = D^{\ell{}_1t}_{j} (1 - \rho).
\end{align}

Consider now the case where $v^{\ell{}_1t}_{ij} (\boldsymbol{y}^{\prime\prime}, \boldsymbol{z}^{\star}) = 1$. 
If $v^{\ell{}_1t}_{ij} (\boldsymbol{y}^{\prime\prime}, \boldsymbol{z}^{\star}) = 1$, then $\rho = 0$ or ${y^{\prime\prime}}^{t}_{i} = 0$, $o^{\ell{}_2st}_{ij} (\boldsymbol{y}^{\prime\prime}) = 1, \forall s \in \mathcal{T}: \ell{}_2 < s \leq \ell{}_1$, and $o^{\ell{}_2st}_{ij} (\boldsymbol{y}^{\prime\prime}) = 0,  \forall s \in \mathcal{T}: \ell{}_1 < s \leq t$.
The inequality then forces $v^{\ell{}_2t}_{ij} (\boldsymbol{0}, \boldsymbol{z}^{\star}) = 1$ and results in a contradiction
\begin{align}
    D^{\ell{}_1t}_{j} < 
        \left( D^{\ell{}_2t}_{j}  - \sum_{\substack{s \in \mathcal{T}: \\ \ell{}_2 < s \leq \ell{}_1}} d^{s}_{j} \right) = D^{\ell{}_1t}_{j}.
\end{align}

Since this inequality can never hold, we have shown that the lemma holds by contradiction. \Halmos

\subsubsection{Proof of Theorem~\ref{thm:dominance}} \label{apx:dominance}

We show that Inequality~\eqref{eq:improved} dominates Inequality~\eqref{eq:tailored-cut} for the follower problem $\mathcal{F} (\boldsymbol{y}^{\prime})$ when built with the bilevel feasible solution $(\boldsymbol{y}^{\prime}, \boldsymbol{z}^{\star})$.
Assume, for the sake of contradiction, that this is false.
This implies that there is a bilevel feasible solution $(\boldsymbol{y}^{\prime\prime}, \boldsymbol{z}^{\prime\prime})$ such that 
\begin{align*}
    \sum_{t \in \mathcal{T}} \sum_{\substack{\ell{} \in \mathcal{T}^{S} :\\ \ell{} < t}} \sum_{j \in \mathcal{J}} \sum_{\substack{i \in \mathcal{I}: \\ i \succ_{j} 0}}  r_{i} v^{\ell{}t}_{ij} (\boldsymbol{0}, \boldsymbol{z}^{\star}) \left( D^{\ell{}t}_{j} - \sum_{\substack{s \in \mathcal{T}: \\ \ell{} < s \leq t}} d^{s}_{j} o^{\ell{}st}_{ij} (\boldsymbol{y})- \rho y^{t}_{i}\sum_{\substack{s \in \mathcal{T}: \\ \ell{} < s \leq t}} d^{s}_{j} (1 - o^{\ell{}st}_{ij} (\boldsymbol{y})) \right) < \\
        \sum_{t \in \mathcal{T}} \sum_{\substack{\ell{} \in \mathcal{T}^{S}: \\ \ell{} < t}} \sum_{j \in \mathcal{J}} \sum_{\substack{i \in \mathcal{I}: \\ i \succ_{j} 0}}  r_{i} D^{\ell{}t}_{j} v^{\ell{}t}_{ij} (\boldsymbol{y}^{\prime}, \boldsymbol{z}^{\star}) \left( 1 - \sum_{\substack{k \in \mathcal{I}: \\ k \succ_{j} 0}} \sum_{\substack{s \in \mathcal{T}: \\ \ell{} < s < t}} y^{s}_{k} - \sum_{\substack{k \in \mathcal{I}: \\ k \succ_{j} i}} y^{t}_{k} - \mathbb{I} \left[ v^{\ell{}t}_{ij} (\boldsymbol{y}^{\prime}, \boldsymbol{z}^{\star}) = 1 \land \rho > 0 \right] y^{t}_{i} \right).
\end{align*}
There is therefore at least one combination of customer $j$ and period $t$ (with location $i$) such that 
\begin{align*}
    \sum_{\substack{\ell{} \in \mathcal{T}^{S} :\\ \ell{} < t}} r_{i} v^{\ell{}t}_{ij} (\boldsymbol{0}, \boldsymbol{z}^{\star}) \left( D^{\ell{}t}_{j} - \sum_{\substack{s \in \mathcal{T}: \\ \ell{} < s \leq t}} d^{s}_{j} o^{\ell{}st}_{ij} (\boldsymbol{y}) - \rho y^{t}_{i}\sum_{\substack{s \in \mathcal{T}: \\ \ell{} < s \leq t}} d^{s}_{j} (1 - o^{\ell{}st}_{ij}(\boldsymbol{y}))  \right) < \\
        \sum_{\substack{\ell{} \in \mathcal{T}^{S}: \\ \ell{} < t}} r_{i} D^{\ell{}t}_{j} v^{\ell{}t}_{ij} (\boldsymbol{y}^{\prime}, \boldsymbol{z}^{\star}) \left( 1 - \sum_{\substack{k \in \mathcal{I}: \\ k \succ_{j} 0}} \sum_{\substack{s \in \mathcal{T}: \\ \ell{} < s < t}} y^{s}_{k} - \sum_{\substack{k \in \mathcal{I}: \\ k \succ_{j} i}} y^{t}_{k} - \mathbb{I} \left[ v^{\ell{}t}_{ij} (\boldsymbol{y}^{\prime}, \boldsymbol{z}^{\star}) = 1 \land \rho > 0 \right] y^{t}_{i} \right).
\end{align*}
We abuse the notation to simply write location $i$ as the location that captures customer $j$ at period $t$ within location schedule $\boldsymbol{z}^{\star}$.
We further analyze whether this inequality can ever hold.

The left-hand side is strictly greater than zero, so we can rewrite the inequality with periods $\ell{}_1$ and $\ell{}_2$ as follows, where $v^{\ell{}_1t}_{ij} (\boldsymbol{0}, \boldsymbol{z}^{\star}) = 1$ and $\ell{}_1 \leq \ell{}_2$ (for more details, see Section~\ref{sec:decomposition}):
\begin{align*}
    r_{i} v^{\ell{}_1t}_{ij} (\boldsymbol{0}, \boldsymbol{z}^{\star}) \left( D^{\ell{}_1t}_{j} - \sum_{\substack{s \in \mathcal{T}: \\ \ell{} < s \leq t}} d^{s}_{j} o^{\ell{}_1st}_{ij} (\boldsymbol{y}) - \rho y^{t}_{i}\sum_{\substack{s \in \mathcal{T}: \\ \ell{}_1 < s \leq t}} d^{s}_{j} (1 - o^{\ell{}_1st}_{ij}(\boldsymbol{y}))  \right) < \\
        r_{i} D^{\ell{}_2t}_{j} v^{\ell{}_2t}_{ij} (\boldsymbol{y}^{\prime}, \boldsymbol{z}^{\star}) \left( 1 - \sum_{\substack{k \in \mathcal{I}: \\ k \succ_{j} 0}} \sum_{\substack{s \in \mathcal{T}: \\ \ell{}_2 < s < t}} y^{s}_{k} - \sum_{\substack{k \in \mathcal{I}: \\ k \succ_{j} i}} y^{t}_{k} - \mathbb{I} \left[ v^{\ell{_2}t}_{ij} (\boldsymbol{y}^{\prime}, \boldsymbol{z}^{\star}) = 1 \land \rho > 0 \right] y^{t}_{i} \right).
\end{align*}

Assume that $\not\exists y^{s}_{k} = 1, s \in \mathcal{T} : \ell{}_2 < s \leq t, k \in \mathcal{I} : k \succ_{j} 0$ and $\not\exists y^{t}_{k} = 1, k \in \mathcal{I} : k \succ_{j} i$. 
If $y^{t}_{i} = 0$, the inequality simplifies to $D^{\ell{}_3t}_{j} < D^{\ell{}_2t}_{j}  v^{\ell{}_2t}_{ij} (\boldsymbol{y}^{\prime}, \boldsymbol{z}^{\star})$, where $\ell{}_1 \leq \ell{}_3 \leq \ell{}_2$.
Since $D^{\ell{}_3t}_{j} \geq D^{\ell{}_2t}_{j}$, there are no feasible values for $ v^{\ell{}_2t}_{ij} (\boldsymbol{y}^{\prime}, \boldsymbol{z}^{\star})$ so that the inequality holds.
If $y^{t}_{i} = 1$, then the inequality simplifies to $D^{\ell{}_3t}_{j} (1 -\rho) < D^{\ell{}_2t}_{j}  v^{\ell{}_2t}_{ij} (\boldsymbol{y}^{\prime}, \boldsymbol{z}^{\star}) \left(1  - \mathbb{I} \left[ v^{\ell{_2}t}_{ij} (\boldsymbol{y}^{\prime}, \boldsymbol{z}^{\star}) = 1 \land \rho > 0 \right] \right)$, where $\ell{}_1 \leq \ell{}_3 \leq \ell{}_2$.
If $v^{\ell{_2}t}_{ij} (\boldsymbol{y}^{\prime}, \boldsymbol{z}^{\star}) = 1 - \rho$ or $\rho = 0$, then the inequality becomes $D^{\ell{}_3t}_{j} (1 -\rho) < D^{\ell{}_2t}_{j} (1 - \rho)$ and cannot be satisfied since $D^{\ell{}_3t}_{j} \geq D^{\ell{}_2t}_{j}$.
If $v^{\ell{_2}t}_{ij} (\boldsymbol{y}^{\prime}, \boldsymbol{z}^{\star}) = 1$ and $\rho > 0$, then the inequality becomes $D^{\ell{}_3t}_{j} (1 -\rho) < 0$, which is absurd.

Assume that $\exists y^{s}_{k} = 1, s \in \mathcal{T}: \ell{}_2 < s \leq t$ or $\exists y^{t}_{k} = 1, k \in \mathcal{I} : k \succ_{j} i$.
In this context, let us compute $\ell{}_3 = \max \{ \max_{\substack{s \in \mathcal{T}: \\ \ell{}_2 < s \leq t \\ k \in \mathcal{I}: \\ k \succ_{j} 0}} \{ s y^{s}_{k}\}, \max_{\substack{k \in \mathcal{I}: \\ k \succ_{j} i}} \{ t y^{t}_{k}\} \}$ (\ie, the latest intervention of the leader in this scenario).
The right-hand side simplifies to $D^{\ell{}_3t}_{j}$, whereas the right-hand side simplifies to zero (or strictly negative), which is absurd.

Since this inequality can never hold, we have shown that the theorem holds by contradiction. \Halmos

\section{Computational Benchmark} \label{apx:benchmark}

We provide a detailed explanation of how we create our computational benchmark. Recall that some parameters are generated randomly with seed values $S \in \mathset{1,2,3,4,5}$.

We considered each federal electoral district in Quebec as a customer.
There are $78$ districts in the Quebec province, among which $40$ districts are in the Montreal region.
We obtained the population size $P_{j}$ of each district $j$ from Statistics Canada, catalogue number \texttt{98-401-X2021029}, and computed a centroid $(x,y)_{j}$ based on geographical data from OpenStreetMap \citep{OpenStreetMap}.
%
More specifically, we calculated the centroid of a district as the average coordinate of the buildings within the district boundary.
To avoid large coefficients in the mathematical formulations, we scaled down the original population sizes $P^{old}$ by computing $P^{new}_{j} = \ceil{\frac{P^{old}_{j}}{10^4}}$.
We randomly sampled half of targeted customers to become candidate locations according to the seed value.
We then computed the driving distance $DD_{ij}$ between customer $j$ and location $i$ based on centroids $(x,y)_{j}$ and $(x,y)_{i}$ with an open-source OpenStreetMap router\footnote{https://routing.openstreetmap.de/}, and built customer rankings such that $i \succ_{j} k$ if and only if $DD_{ij} < DD_{kj}$ and $i \succ_{j} 0$ if and only if $DD_{ij} < MAX\_DD$, where $MAX\_DD$ is the maximum driving distance of the instance.
The remaining parameters are set as explained in Section~\ref{sec:benchmark}.

\section{Additional Experiments} \label{apx:results}

We present here additional numerical results to support our claims.
Appendix~\ref{apx:preliminary} presents the comparison of our solution methods with the state-of-the-art general-purpose bilevel solver, Appendix~\ref{apx:supplementary} features supplementary performance results for both variants of the \ourproblem{}, and Appendix~\ref{apx:managerial} contains supplementary managerial insights for both variants of the \ourproblem{}.

\subsection{General-Purpose Bilevel Solver} \label{apx:preliminary}

We compare here our solution methods with the state-of-the-art general-purpose bilevel solver MIX++ \citep{fischettiNewGeneralPurposeAlgorithm2017}.
Recall that we only consider the optimistic variant of the \ourproblem{} in these experiments, as the MIX++ Solver cannot handle the pessimistic variant.
Since we were unable to deploy the MIX++ Solver on the Nibi server, we ran it together with our solution methods on a local machine equipped with two Intel(R) Xeon(R) Gold 6226 CPUs (12 cores each).
Each solution method had a time limit of $1$ hour and was limited to a single thread to avoid bias related to computational resources.
Our solution methods employ Python (version \texttt{3.10}) and Gurobi (version \texttt{12.0}), whereas the MIX++ Solver uses CPLEX (version \texttt{22.1}).
We consider a subset of $720$ instances (out of $1800$) where $\rho \in \mathset{0,1}$, as the remaining instances do not satisfy a key assumption required by the MIX++ Solver \citep[for more details, see][Assumption 3]{fischettiNewGeneralPurposeAlgorithm2017}.

\begin{table}[!ht]
    \centering
    \caption{Percentage of the benchmark solved to optimality by each exact method, grouped by dimensional instance attributes (\ie, number of locations $|\mathcal{I}|$, customers $|\mathcal{J}|$, and periods $|\mathcal{T}|$). [Optimistic variant of the \ourproblem{}]}
    \small
    \begin{tabular}{ccccccc}
    \toprule
        \multirow{2}{*}{$|\mathcal{T}|$} & \multicolumn{3}{c}{Montreal ($|\mathcal{I}| = 20, |\mathcal{J}| = 40$)} & \multicolumn{3}{c}{Quebec ($|\mathcal{I}| = 39, |\mathcal{J}| = 78$)}  \\ \cmidrule(lr){2-7}
        & MIX++ Solver & Tailored B\&C& Tightened B\&C& MIX++ Solver & Tailored B\&C & Tightened B\&C\\ \midrule
        3&5.00\%&100.00\%&97.50\%&4.17\%&100.00\%&97.50\%\\ \midrule
5&0.00\%&34.17\%&50.00\%&0.00\%&46.67\%&54.17\%\\ \midrule
7&0.00\%&2.50\%&10.00\%&0.00\%&6.67\%&10.00\%\\ \bottomrule
    \end{tabular}
    \label{tab:summary-preliminary}
\end{table}
Table~\ref{tab:summary-preliminary} presents the percentage of the (filtered) benchmark solved to optimality by each exact methods, grouped by dimensional instance attributes.
The MIX++ Solver only proves optimality for a small percentage of the instances with $|\mathcal{T}| = 3$ periods, hitting the time limit without a proven optimal solution for the remaining instances.
This outcome suggests that our solution methods have a considerably better performance than the MIX++ Solver, showcasing the advantage of exploring the problem structure.
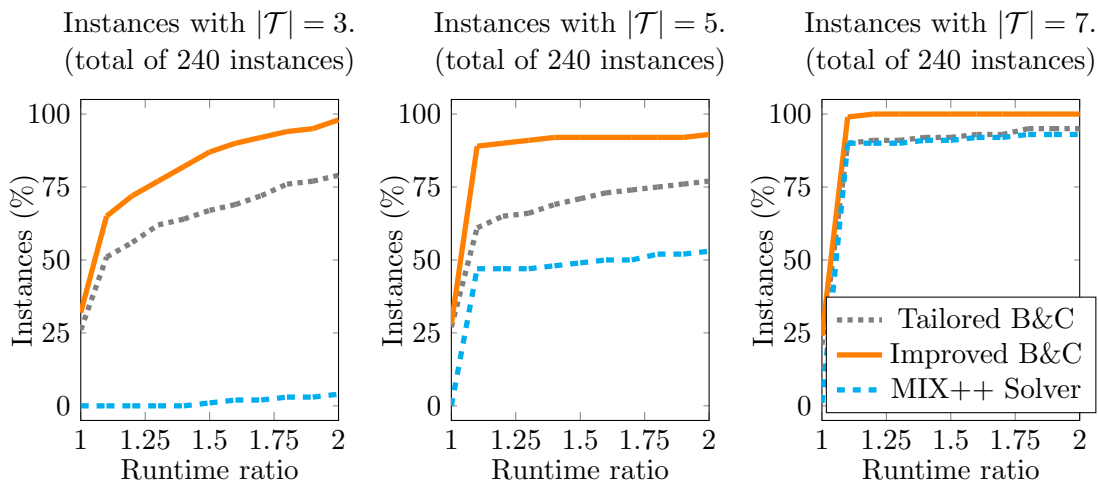
\begin{figure}[!ht]
     \centering
     \caption{Runtime ratios of both exact methods for exact methods, where the $y$ axis presents the percentage of instances with a ratio smaller than or equal to the reference value on the $x$ axis. [Optimistic variant of the \ourproblem{}]}
    \begin{tikzpicture}
\begin{axis}[name = 3periods, xmin = 1, xmax = 2, xscale = 0.5, ymin = -5, ymax = 105, yscale = 0.75, xtick = {1,1.25,1.5,1.75,2}, ytick = {0, 25, 50, 75, 100}, xlabel = {Runtime ratio}, ylabel = {Instances (\%)}, legend pos = south west]
\addplot[line width = 2pt, color = gray, style = dotted] coordinates {(1,26) (1.1,51)};
\addplot[line width = 2pt, color = gray, style = dotted] coordinates {(1.1,51) (1.2,56)};
\addplot[line width = 2pt, color = gray, style = dotted] coordinates {(1.2,56) (1.3,62)};
\addplot[line width = 2pt, color = gray, style = dotted] coordinates {(1.3,62) (1.4,64)};
\addplot[line width = 2pt, color = gray, style = dotted] coordinates {(1.4,64) (1.5,67)};
\addplot[line width = 2pt, color = gray, style = dotted] coordinates {(1.5,67) (1.6,69)};
\addplot[line width = 2pt, color = gray, style = dotted] coordinates {(1.6,69) (1.7,72)};
\addplot[line width = 2pt, color = gray, style = dotted] coordinates {(1.7,72) (1.8,76)};
\addplot[line width = 2pt, color = gray, style = dotted] coordinates {(1.8,76) (1.9,77)};
\addplot[line width = 2pt, color = gray, style = dotted] coordinates {(1.9,77) (2.0,79)};
\addplot[line width = 2pt, color = orange, style = solid] coordinates {(1,32) (1.1,65)};
\addplot[line width = 2pt, color = orange, style = solid] coordinates {(1.1,65) (1.2,72)};
\addplot[line width = 2pt, color = orange, style = solid] coordinates {(1.2,72) (1.3,77)};
\addplot[line width = 2pt, color = orange, style = solid] coordinates {(1.3,77) (1.4,82)};
\addplot[line width = 2pt, color = orange, style = solid] coordinates {(1.4,82) (1.5,87)};
\addplot[line width = 2pt, color = orange, style = solid] coordinates {(1.5,87) (1.6,90)};
\addplot[line width = 2pt, color = orange, style = solid] coordinates {(1.6,90) (1.7,92)};
\addplot[line width = 2pt, color = orange, style = solid] coordinates {(1.7,92) (1.8,94)};
\addplot[line width = 2pt, color = orange, style = solid] coordinates {(1.8,94) (1.9,95)};
\addplot[line width = 2pt, color = orange, style = solid] coordinates {(1.9,95) (2.0,98)};
\addplot[line width = 2pt, color = cyan, style = dashed] coordinates {(1,0) (1.1,0)};
\addplot[line width = 2pt, color = cyan, style = dashed] coordinates {(1.1,0) (1.2,0)};
\addplot[line width = 2pt, color = cyan, style = dashed] coordinates {(1.2,0) (1.3,0)};
\addplot[line width = 2pt, color = cyan, style = dashed] coordinates {(1.3,0) (1.4,0)};
\addplot[line width = 2pt, color = cyan, style = dashed] coordinates {(1.4,0) (1.5,1)};
\addplot[line width = 2pt, color = cyan, style = dashed] coordinates {(1.5,1) (1.6,2)};
\addplot[line width = 2pt, color = cyan, style = dashed] coordinates {(1.6,2) (1.7,2)};
\addplot[line width = 2pt, color = cyan, style = dashed] coordinates {(1.7,2) (1.8,3)};
\addplot[line width = 2pt, color = cyan, style = dashed] coordinates {(1.8,3) (1.9,3)};
\addplot[line width = 2pt, color = cyan, style = dashed] coordinates {(1.9,3) (2.0,4)};

\end{axis}
\begin{axis}[name = 5periods, at={(3periods.east)}, anchor=west, xshift = 1.5cm, xmin = 1, xmax = 2, xscale = 0.5, ymin = -5, ymax = 105, yscale = 0.75, xtick = {1,1.25,1.5,1.75,2}, ytick = {0, 25, 50, 75, 100}, xlabel = {Runtime ratio}, ylabel = {Instances (\%)}, legend pos = south west]
\addplot[line width = 2pt, color = gray, style = dotted] coordinates {(1,27) (1.1,61)};
\addplot[line width = 2pt, color = gray, style = dotted] coordinates {(1.1,61) (1.2,65)};
\addplot[line width = 2pt, color = gray, style = dotted] coordinates {(1.2,65) (1.3,66)};
\addplot[line width = 2pt, color = gray, style = dotted] coordinates {(1.3,66) (1.4,69)};
\addplot[line width = 2pt, color = gray, style = dotted] coordinates {(1.4,69) (1.5,71)};
\addplot[line width = 2pt, color = gray, style = dotted] coordinates {(1.5,71) (1.6,73)};
\addplot[line width = 2pt, color = gray, style = dotted] coordinates {(1.6,73) (1.7,74)};
\addplot[line width = 2pt, color = gray, style = dotted] coordinates {(1.7,74) (1.8,75)};
\addplot[line width = 2pt, color = gray, style = dotted] coordinates {(1.8,75) (1.9,76)};
\addplot[line width = 2pt, color = gray, style = dotted] coordinates {(1.9,76) (2.0,77)};
\addplot[line width = 2pt, color = orange, style = solid] coordinates {(1,28) (1.1,89)};
\addplot[line width = 2pt, color = orange, style = solid] coordinates {(1.1,89) (1.2,90)};
\addplot[line width = 2pt, color = orange, style = solid] coordinates {(1.2,90) (1.3,91)};
\addplot[line width = 2pt, color = orange, style = solid] coordinates {(1.3,91) (1.4,92)};
\addplot[line width = 2pt, color = orange, style = solid] coordinates {(1.4,92) (1.5,92)};
\addplot[line width = 2pt, color = orange, style = solid] coordinates {(1.5,92) (1.6,92)};
\addplot[line width = 2pt, color = orange, style = solid] coordinates {(1.6,92) (1.7,92)};
\addplot[line width = 2pt, color = orange, style = solid] coordinates {(1.7,92) (1.8,92)};
\addplot[line width = 2pt, color = orange, style = solid] coordinates {(1.8,92) (1.9,92)};
\addplot[line width = 2pt, color = orange, style = solid] coordinates {(1.9,92) (2.0,93)};
\addplot[line width = 2pt, color = cyan, style = dashed] coordinates {(1,0) (1.1,47)};
\addplot[line width = 2pt, color = cyan, style = dashed] coordinates {(1.1,47) (1.2,47)};
\addplot[line width = 2pt, color = cyan, style = dashed] coordinates {(1.2,47) (1.3,47)};
\addplot[line width = 2pt, color = cyan, style = dashed] coordinates {(1.3,47) (1.4,48)};
\addplot[line width = 2pt, color = cyan, style = dashed] coordinates {(1.4,48) (1.5,49)};
\addplot[line width = 2pt, color = cyan, style = dashed] coordinates {(1.5,49) (1.6,50)};
\addplot[line width = 2pt, color = cyan, style = dashed] coordinates {(1.6,50) (1.7,50)};
\addplot[line width = 2pt, color = cyan, style = dashed] coordinates {(1.7,50) (1.8,52)};
\addplot[line width = 2pt, color = cyan, style = dashed] coordinates {(1.8,52) (1.9,52)};
\addplot[line width = 2pt, color = cyan, style = dashed] coordinates {(1.9,52) (2.0,53)};

\end{axis}
\begin{axis}[name = 7periods, at={(5periods.east)}, anchor=west, xshift = 1.5cm, xmin = 1, xmax = 2, xscale = 0.5, ymin = -5, ymax = 105, yscale = 0.75, xtick = {1,1.25,1.5,1.75,2}, ytick = {0, 25, 50, 75, 100}, xlabel = {Runtime ratio}, ylabel = {Instances (\%)}, legend pos = south west]
\addlegendimage{line width = 2, color = gray, style = dotted}
\addlegendentry{Tailored B\&C}
\addlegendimage{line width = 2, color = orange, style = solid}
\addlegendentry{Improved B\&C}
\addlegendimage{line width = 2, color = cyan, style = dashed}
\addlegendentry{MIX++ Solver}
\addplot[line width = 2pt, color = gray, style = dotted] coordinates {(1,21) (1.1,90)};
\addplot[line width = 2pt, color = gray, style = dotted] coordinates {(1.1,90) (1.2,91)};
\addplot[line width = 2pt, color = gray, style = dotted] coordinates {(1.2,91) (1.3,91)};
\addplot[line width = 2pt, color = gray, style = dotted] coordinates {(1.3,91) (1.4,92)};
\addplot[line width = 2pt, color = gray, style = dotted] coordinates {(1.4,92) (1.5,92)};
\addplot[line width = 2pt, color = gray, style = dotted] coordinates {(1.5,92) (1.6,93)};
\addplot[line width = 2pt, color = gray, style = dotted] coordinates {(1.6,93) (1.7,93)};
\addplot[line width = 2pt, color = gray, style = dotted] coordinates {(1.7,93) (1.8,95)};
\addplot[line width = 2pt, color = gray, style = dotted] coordinates {(1.8,95) (1.9,95)};
\addplot[line width = 2pt, color = gray, style = dotted] coordinates {(1.9,95) (2.0,95)};
\addplot[line width = 2pt, color = orange, style = solid] coordinates {(1,24) (1.1,99)};
\addplot[line width = 2pt, color = orange, style = solid] coordinates {(1.1,99) (1.2,100)};
\addplot[line width = 2pt, color = orange, style = solid] coordinates {(1.2,100) (1.3,100)};
\addplot[line width = 2pt, color = orange, style = solid] coordinates {(1.3,100) (1.4,100)};
\addplot[line width = 2pt, color = orange, style = solid] coordinates {(1.4,100) (1.5,100)};
\addplot[line width = 2pt, color = orange, style = solid] coordinates {(1.5,100) (1.6,100)};
\addplot[line width = 2pt, color = orange, style = solid] coordinates {(1.6,100) (1.7,100)};
\addplot[line width = 2pt, color = orange, style = solid] coordinates {(1.7,100) (1.8,100)};
\addplot[line width = 2pt, color = orange, style = solid] coordinates {(1.8,100) (1.9,100)};
\addplot[line width = 2pt, color = orange, style = solid] coordinates {(1.9,100) (2.0,100)};
\addplot[line width = 2pt, color = cyan, style = dashed] coordinates {(1,1) (1.1,90)};
\addplot[line width = 2pt, color = cyan, style = dashed] coordinates {(1.1,90) (1.2,90)};
\addplot[line width = 2pt, color = cyan, style = dashed] coordinates {(1.2,90) (1.3,90)};
\addplot[line width = 2pt, color = cyan, style = dashed] coordinates {(1.3,90) (1.4,91)};
\addplot[line width = 2pt, color = cyan, style = dashed] coordinates {(1.4,91) (1.5,91)};
\addplot[line width = 2pt, color = cyan, style = dashed] coordinates {(1.5,91) (1.6,92)};
\addplot[line width = 2pt, color = cyan, style = dashed] coordinates {(1.6,92) (1.7,92)};
\addplot[line width = 2pt, color = cyan, style = dashed] coordinates {(1.7,92) (1.8,93)};
\addplot[line width = 2pt, color = cyan, style = dashed] coordinates {(1.8,93) (1.9,93)};
\addplot[line width = 2pt, color = cyan, style = dashed] coordinates {(1.9,93) (2.0,93)};

\end{axis}
\node[align=center, yshift=1cm] at (3periods.north) {Instances with $|\mathcal{T}| = 3$.};
\node[align=center, yshift=0.5cm] at (3periods.north) {(total of $240$ instances)};
\node[align=center, yshift=1cm] at (5periods.north) {Instances with $|\mathcal{T}| = 5$.};
\node[align=center, yshift=0.5cm] at (5periods.north) {(total of $240$ instances)};
\node[align=center, yshift=1cm] at (7periods.north) {Instances with $|\mathcal{T}| = 7$.};
\node[align=center, yshift=0.5cm] at (7periods.north) {(total of $240$ instances)};
\end{tikzpicture}
    \label{fig:runtimes-preliminary}
\end{figure}
\begin{figure}[!ht]
     \centering
     \caption{Objective ratios of both exact methods for exact methods, where the $y$ axis presents the percentage of instances with a ratio smaller than or equal to the reference value on the $x$ axis. [Optimistic variant of the \ourproblem{}]}    
    \begin{tikzpicture}
\begin{axis}[name = 3periods, xmin = 1, xmax = 2, xscale = 0.5, ymin = 45, ymax = 105, yscale = 0.75, xtick = {1,1.25,1.5,1.75,2}, ytick = {50, 75, 100}, xlabel = {Objective ratio}, ylabel = {Instances (\%)}, legend pos = south west]
\addplot[line width = 2pt, color = gray, style = dotted] coordinates {(1,100) (1.1,100)};
\addplot[line width = 2pt, color = gray, style = dotted] coordinates {(1.1,100) (1.2,100)};
\addplot[line width = 2pt, color = gray, style = dotted] coordinates {(1.2,100) (1.3,100)};
\addplot[line width = 2pt, color = gray, style = dotted] coordinates {(1.3,100) (1.4,100)};
\addplot[line width = 2pt, color = gray, style = dotted] coordinates {(1.4,100) (1.5,100)};
\addplot[line width = 2pt, color = gray, style = dotted] coordinates {(1.5,100) (1.6,100)};
\addplot[line width = 2pt, color = gray, style = dotted] coordinates {(1.6,100) (1.7,100)};
\addplot[line width = 2pt, color = gray, style = dotted] coordinates {(1.7,100) (1.8,100)};
\addplot[line width = 2pt, color = gray, style = dotted] coordinates {(1.8,100) (1.9,100)};
\addplot[line width = 2pt, color = gray, style = dotted] coordinates {(1.9,100) (2.0,100)};
\addplot[line width = 2pt, color = orange, style = solid] coordinates {(1,99) (1.1,100)};
\addplot[line width = 2pt, color = orange, style = solid] coordinates {(1.1,100) (1.2,100)};
\addplot[line width = 2pt, color = orange, style = solid] coordinates {(1.2,100) (1.3,100)};
\addplot[line width = 2pt, color = orange, style = solid] coordinates {(1.3,100) (1.4,100)};
\addplot[line width = 2pt, color = orange, style = solid] coordinates {(1.4,100) (1.5,100)};
\addplot[line width = 2pt, color = orange, style = solid] coordinates {(1.5,100) (1.6,100)};
\addplot[line width = 2pt, color = orange, style = solid] coordinates {(1.6,100) (1.7,100)};
\addplot[line width = 2pt, color = orange, style = solid] coordinates {(1.7,100) (1.8,100)};
\addplot[line width = 2pt, color = orange, style = solid] coordinates {(1.8,100) (1.9,100)};
\addplot[line width = 2pt, color = orange, style = solid] coordinates {(1.9,100) (2.0,100)};
\addplot[line width = 2pt, color = cyan, style = dashed] coordinates {(1,34) (1.1,66)};
\addplot[line width = 2pt, color = cyan, style = dashed] coordinates {(1.1,66) (1.2,82)};
\addplot[line width = 2pt, color = cyan, style = dashed] coordinates {(1.2,82) (1.3,89)};
\addplot[line width = 2pt, color = cyan, style = dashed] coordinates {(1.3,89) (1.4,92)};
\addplot[line width = 2pt, color = cyan, style = dashed] coordinates {(1.4,92) (1.5,94)};
\addplot[line width = 2pt, color = cyan, style = dashed] coordinates {(1.5,94) (1.6,96)};
\addplot[line width = 2pt, color = cyan, style = dashed] coordinates {(1.6,96) (1.7,98)};
\addplot[line width = 2pt, color = cyan, style = dashed] coordinates {(1.7,98) (1.8,98)};
\addplot[line width = 2pt, color = cyan, style = dashed] coordinates {(1.8,98) (1.9,100)};
\addplot[line width = 2pt, color = cyan, style = dashed] coordinates {(1.9,100) (2.0,100)};

\end{axis}
\begin{axis}[name = 5periods, at={(3periods.east)}, anchor=west, xshift = 1.5cm, xmin = 1, xmax = 2, xscale = 0.5, ymin = 45, ymax = 105, yscale = 0.75, xtick = {1,1.25,1.5,1.75,2}, ytick = {50, 75, 100}, xlabel = {Objective ratio}, ylabel = {Instances (\%)}, legend pos = south west]
\addplot[line width = 2pt, color = gray, style = dotted] coordinates {(1,83) (1.1,97)};
\addplot[line width = 2pt, color = gray, style = dotted] coordinates {(1.1,97) (1.2,99)};
\addplot[line width = 2pt, color = gray, style = dotted] coordinates {(1.2,99) (1.3,100)};
\addplot[line width = 2pt, color = gray, style = dotted] coordinates {(1.3,100) (1.4,100)};
\addplot[line width = 2pt, color = gray, style = dotted] coordinates {(1.4,100) (1.5,100)};
\addplot[line width = 2pt, color = gray, style = dotted] coordinates {(1.5,100) (1.6,100)};
\addplot[line width = 2pt, color = gray, style = dotted] coordinates {(1.6,100) (1.7,100)};
\addplot[line width = 2pt, color = gray, style = dotted] coordinates {(1.7,100) (1.8,100)};
\addplot[line width = 2pt, color = gray, style = dotted] coordinates {(1.8,100) (1.9,100)};
\addplot[line width = 2pt, color = gray, style = dotted] coordinates {(1.9,100) (2.0,100)};
\addplot[line width = 2pt, color = orange, style = solid] coordinates {(1,71) (1.1,92)};
\addplot[line width = 2pt, color = orange, style = solid] coordinates {(1.1,92) (1.2,96)};
\addplot[line width = 2pt, color = orange, style = solid] coordinates {(1.2,96) (1.3,98)};
\addplot[line width = 2pt, color = orange, style = solid] coordinates {(1.3,98) (1.4,99)};
\addplot[line width = 2pt, color = orange, style = solid] coordinates {(1.4,99) (1.5,99)};
\addplot[line width = 2pt, color = orange, style = solid] coordinates {(1.5,99) (1.6,100)};
\addplot[line width = 2pt, color = orange, style = solid] coordinates {(1.6,100) (1.7,100)};
\addplot[line width = 2pt, color = orange, style = solid] coordinates {(1.7,100) (1.8,100)};
\addplot[line width = 2pt, color = orange, style = solid] coordinates {(1.8,100) (1.9,100)};
\addplot[line width = 2pt, color = orange, style = solid] coordinates {(1.9,100) (2.0,100)};
\addplot[line width = 2pt, color = cyan, style = dashed] coordinates {(1,3) (1.1,19)};
\addplot[line width = 2pt, color = cyan, style = dashed] coordinates {(1.1,19) (1.2,32)};
\addplot[line width = 2pt, color = cyan, style = dashed] coordinates {(1.2,32) (1.3,47)};
\addplot[line width = 2pt, color = cyan, style = dashed] coordinates {(1.3,47) (1.4,61)};
\addplot[line width = 2pt, color = cyan, style = dashed] coordinates {(1.4,61) (1.5,74)};
\addplot[line width = 2pt, color = cyan, style = dashed] coordinates {(1.5,74) (1.6,83)};
\addplot[line width = 2pt, color = cyan, style = dashed] coordinates {(1.6,83) (1.7,87)};
\addplot[line width = 2pt, color = cyan, style = dashed] coordinates {(1.7,87) (1.8,90)};
\addplot[line width = 2pt, color = cyan, style = dashed] coordinates {(1.8,90) (1.9,94)};
\addplot[line width = 2pt, color = cyan, style = dashed] coordinates {(1.9,94) (2.0,94)};

\end{axis}
\begin{axis}[name = 7periods, at={(5periods.east)}, anchor=west, xshift = 1.5cm, xmin = 1, xmax = 2, xscale = 0.5, ymin = 45, ymax = 105, yscale = 0.75, xtick = {1,1.25,1.5,1.75,2}, ytick = {50, 75, 100}, xlabel = {Objective ratio}, ylabel = {Instances (\%)}, legend pos = south west]
\addlegendimage{line width = 2, color = gray, style = dotted}
\addlegendentry{Tailored B\&C}
\addlegendimage{line width = 2, color = orange, style = solid}
\addlegendentry{Improved B\&C}
\addlegendimage{line width = 2, color = cyan, style = dashed}
\addlegendentry{MIX++ Solver}
\addplot[line width = 2pt, color = gray, style = dotted] coordinates {(1,61) (1.1,95)};
\addplot[line width = 2pt, color = gray, style = dotted] coordinates {(1.1,95) (1.2,98)};
\addplot[line width = 2pt, color = gray, style = dotted] coordinates {(1.2,98) (1.3,99)};
\addplot[line width = 2pt, color = gray, style = dotted] coordinates {(1.3,99) (1.4,100)};
\addplot[line width = 2pt, color = gray, style = dotted] coordinates {(1.4,100) (1.5,100)};
\addplot[line width = 2pt, color = gray, style = dotted] coordinates {(1.5,100) (1.6,100)};
\addplot[line width = 2pt, color = gray, style = dotted] coordinates {(1.6,100) (1.7,100)};
\addplot[line width = 2pt, color = gray, style = dotted] coordinates {(1.7,100) (1.8,100)};
\addplot[line width = 2pt, color = gray, style = dotted] coordinates {(1.8,100) (1.9,100)};
\addplot[line width = 2pt, color = gray, style = dotted] coordinates {(1.9,100) (2.0,100)};
\addplot[line width = 2pt, color = orange, style = solid] coordinates {(1,54) (1.1,76)};
\addplot[line width = 2pt, color = orange, style = solid] coordinates {(1.1,76) (1.2,87)};
\addplot[line width = 2pt, color = orange, style = solid] coordinates {(1.2,87) (1.3,92)};
\addplot[line width = 2pt, color = orange, style = solid] coordinates {(1.3,92) (1.4,94)};
\addplot[line width = 2pt, color = orange, style = solid] coordinates {(1.4,94) (1.5,94)};
\addplot[line width = 2pt, color = orange, style = solid] coordinates {(1.5,94) (1.6,95)};
\addplot[line width = 2pt, color = orange, style = solid] coordinates {(1.6,95) (1.7,97)};
\addplot[line width = 2pt, color = orange, style = solid] coordinates {(1.7,97) (1.8,98)};
\addplot[line width = 2pt, color = orange, style = solid] coordinates {(1.8,98) (1.9,99)};
\addplot[line width = 2pt, color = orange, style = solid] coordinates {(1.9,99) (2.0,99)};
\addplot[line width = 2pt, color = cyan, style = dashed] coordinates {(1,1) (1.1,6)};
\addplot[line width = 2pt, color = cyan, style = dashed] coordinates {(1.1,6) (1.2,16)};
\addplot[line width = 2pt, color = cyan, style = dashed] coordinates {(1.2,16) (1.3,25)};
\addplot[line width = 2pt, color = cyan, style = dashed] coordinates {(1.3,25) (1.4,35)};
\addplot[line width = 2pt, color = cyan, style = dashed] coordinates {(1.4,35) (1.5,44)};
\addplot[line width = 2pt, color = cyan, style = dashed] coordinates {(1.5,44) (1.6,54)};
\addplot[line width = 2pt, color = cyan, style = dashed] coordinates {(1.6,54) (1.7,62)};
\addplot[line width = 2pt, color = cyan, style = dashed] coordinates {(1.7,62) (1.8,72)};
\addplot[line width = 2pt, color = cyan, style = dashed] coordinates {(1.8,72) (1.9,77)};
\addplot[line width = 2pt, color = cyan, style = dashed] coordinates {(1.9,77) (2.0,82)};

\end{axis}
\node[align=center, yshift=1cm] at (3periods.north) {Instances with $|\mathcal{T}| = 3$.};
\node[align=center, yshift=0.5cm] at (3periods.north) {(total of $240$ instances)};
\node[align=center, yshift=1cm] at (5periods.north) {Instances with $|\mathcal{T}| = 5$.};
\node[align=center, yshift=0.5cm] at (5periods.north) {(total of $240$ instances)};
\node[align=center, yshift=1cm] at (7periods.north) {Instances with $|\mathcal{T}| = 7$.};
\node[align=center, yshift=0.5cm] at (7periods.north) {(total of $240$ instances)};
\end{tikzpicture}
    \label{fig:objectives-preliminary}
\end{figure}
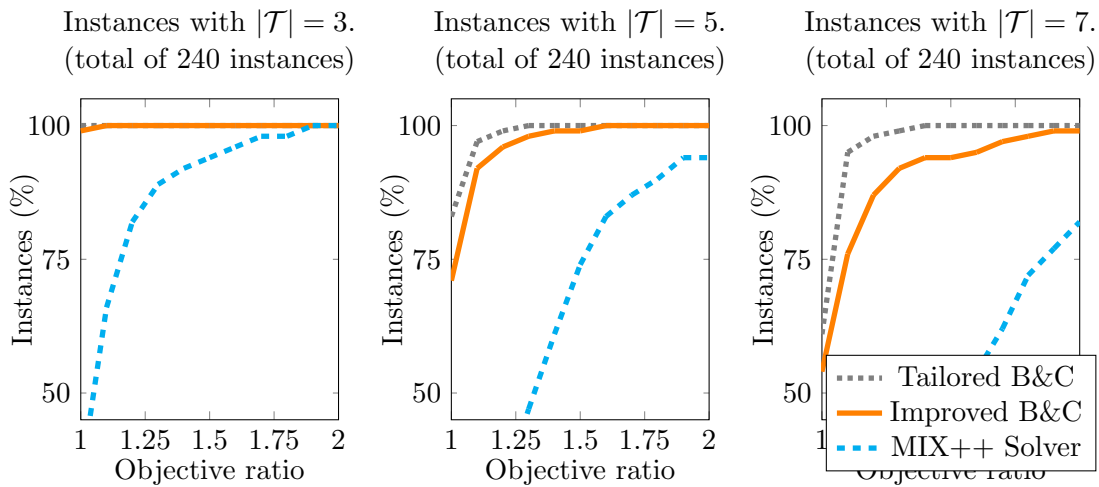
We now compare the solutions obtained by the three exact methods through runtime and objective ratios as done in Section~\ref{sec:performance}. Figures~\ref{fig:runtimes-preliminary} and ~\ref{fig:objectives-preliminary} present these ratios for the three exact methods.
Indeed, the MIX++ solver is considerably slower than our solution methods, and finds solutions that are significantly worse than those found by our branch-and-cut algorithms.

Figure~\ref{fig:objectives-preliminary} also reveals that the Tailored B\&C finds solutions of higher quality than the Tightened B\&C within a time limit of $1$ hour, particularly as the number of periods increase.
This outcome is rather unexpected, since Theorem~\ref{thm:dominance} states that Inequality~\eqref{eq:improved} dominates Inequality~\eqref{eq:tailored-cut}, one would therefore expect the Tightened B\&C to find better solutions than the Tailored B\&C within the same time limit.
The main reason behind this outcome seems to be the overhead of additional variables and constraints required by Inequality~\eqref{eq:improved}.
As the number of periods increases, each node of the Tightened B\&C solves a larger linear program, thus exploring less nodes within the same time limit and, consequently, visiting less bilevel feasible solutions than the Tailored B\&C.
Nevertheless, once we increase the time limit (\eg, to $24$ hours as in Section~\ref{sec:performance}), the Tightened B\&C starts to find solutions of as high-quality as the Tailored B\&C, if not better, within the same time limit.

\subsection{Supplementary Performance Results} \label{apx:supplementary}

We here present supplementary performance results for both variants of the \ourproblem{}.

\subsubsection{Optimistic Variant of the \ourproblem{}}  \label{apx:supplementary-optimistic}
Table~\ref{tab:soltimes-optimistic} reports average solution times and their standard deviations for the 1092 instances solved to optimality by both exact methods.
%
%
%
The Tightened B\&C is approximately three times faster, on average, than the Tailored B\&C ($36.90 \pm 130.25$ versus $97.55 \pm 238.11$ minutes) for $1092$ instances solved to optimality by both exact methods.
%
%
%
Instances with large planning horizons, long maximum travel times, identical rewards, and constant spawning demands seem harder to solve by our solution methods than their counterparts.
This outcome is expected, since these characteristics directly affect the size and density of the respective mathematical formulations.
Nevertheless, the Tightened B\&C presents smaller solution times and proves optimality for more instances than the Tailored B\&C across different instance attributes.
\begin{table}[!ht]
    \centering
    \caption{Average solution times of Tailored B\&C and Tightened B\&C, as well as their standard deviations, for instances solved to optimality by both exact methods. [Optimistic variant of the \ourproblem{}]}
    \begin{tabular}{crrr}
        \toprule
        \multirow{2}*{Instance attributes} &  \# instances & Tailored B\&C & Tightened B\&C \\ \cmidrule{3-4}
        & considered & time (min) & time (min)  \\ \midrule
        Benchmark&1092 (1800)&$97.55\pm238.11$&$36.90\pm130.25$\\ \midrule
$|\mathcal{T}|$ = 3 periods&600 (600)&$3.85\pm9.16$&$2.97\pm7.93$\\$|\mathcal{T}|$ = 5 periods&412 (600)&$178.31\pm289.75$&$67.24\pm167.04$\\$|\mathcal{T}|$ = 7 periods&80 (600)&$384.36\pm401.92$&$135.12\pm255.65$\\ \midrule
Montreal region&537 (900)&$99.31\pm222.21$&$26.11\pm104.76$\\Quebec province&555 (900)&$95.84\pm252.73$&$47.34\pm150.24$\\ \midrule
Travel $\leq$ 15 minutes&463 (600)&$37.78\pm135.68$&$13.79\pm89.62$\\Travel $\leq$ 30 minutes&365 (600)&$185.63\pm319.10$&$72.92\pm179.19$\\Travel $\leq$ 45 minutes&264 (600)&$80.59\pm211.13$&$27.65\pm94.89$\\ \midrule
Identical rewards&493 (900)&$84.83\pm214.25$&$7.89\pm29.89$\\Inverse rewards&599 (900)&$108.01\pm255.79$&$60.78\pm170.15$\\ \midrule
Constant demands&489 (900)&$104.48\pm247.75$&$55.12\pm171.77$\\Sparse demands&603 (900)&$91.93\pm230.04$&$22.13\pm79.64$\\ \midrule
Factor $\rho = 0.0$&203 (360)&$61.89\pm198.64$&$56.58\pm179.91$\\Factor $\rho = 0.25$&198 (360)&$75.89\pm219.25$&$30.35\pm117.05$\\Factor $\rho = 0.5$&229 (360)&$124.71\pm261.19$&$34.43\pm121.61$\\Factor $\rho = 0.75$&232 (360)&$114.59\pm260.26$&$37.73\pm134.40$\\Factor $\rho = 1.0$&230 (360)&$103.43\pm234.74$&$26.80\pm84.00$\\ \bottomrule

    \end{tabular}
    \label{tab:soltimes-optimistic}
\end{table}

Table~\ref{tab:optgaps-optimistic} presents the number of instances solved to optimality, average optimality gaps and their standard deviation for the $708$ instances that were not solved to optimality by both exact methods. The Tightened B\&C proves optimality gaps that are almost half as large, on average, than those proven by the Tailored B\&C ($ 16.97\% \pm 20.91\%$ versus $28.94\% \pm 19.10\%$).
We highlight that instances 
with splitting factor $\rho = 0$ present larger optimality gaps than those with splitting factor $\rho > 0$.
Since \textit{(i)} small $\rho$ values increase the set of feasible reactions for the follower (in practical terms, co-location favours the follower) and (\textit{ii)} the high-point relaxation focuses on the leader's profit, our branch-and-cut algorithm likely requires the enumeration of several reactions of the follower that interfere as little as possible with the location schedule of the leader to close the gap and thus prove optimality. 
This outcome indicates that, at least when solved by a branch-and-cut framework based on value-function cuts, flexible competitive (location) problems, where the follower has more feasible reactions due to the splitting factor, are considerably harder to tackle than their interdiction counterparts, where these reactions would be simply infeasible for the follower.
\begin{table}[!ht]
    \centering
    \caption{Number of instances solved to optimality, average optimality gaps and their standard deviations, of Tailored B\&C and Tightened B\&C for instances not solved to optimality by both exact methods. [Optimistic variant of the \ourproblem{}]}
    \begin{tabular}{crrrrr}
        \toprule
        \multirow{2}*{Instance attributes} &  \# instances & \multicolumn{2}{c}{Tailored B\&C} & \multicolumn{2}{c}{Tightened B\&C} \\ \cmidrule{3-6}
        & considered & \# opt. & opt. gap (\%) & \# opt. & opt. gap (\%)  \\
        \midrule
        Benchmark&708 (1800)&$3$&$28.94\pm19.10$&$230$&$16.97\pm20.91$\\ \midrule
$|\mathcal{T}|$ = 3 periods&0 (600)&$0$&$$-$$&$0$&$$-$$\\$|\mathcal{T}|$ = 5 periods&188 (600)&$0$&$26.21\pm18.14$&$95$&$15.57\pm21.22$\\$|\mathcal{T}|$ = 7 periods&520 (600)&$3$&$29.92\pm19.36$&$135$&$17.47\pm20.79$\\ \midrule
Montreal region&363 (900)&$1$&$34.02\pm19.00$&$132$&$20.73\pm23.12$\\Quebec province&345 (900)&$2$&$23.59\pm17.73$&$98$&$13.01\pm17.48$\\ \midrule
Travel $\leq$ 15 minutes&137 (600)&$3$&$19.50\pm13.62$&$82$&$1.93\pm3.58$\\Travel $\leq$ 30 minutes&235 (600)&$0$&$31.32\pm17.45$&$56$&$18.59\pm18.88$\\Travel $\leq$ 45 minutes&336 (600)&$0$&$31.12\pm20.91$&$92$&$21.96\pm23.44$\\ \midrule
Identical rewards&407 (900)&$0$&$33.55\pm18.04$&$180$&$17.53\pm23.36$\\Inverse rewards&301 (900)&$3$&$22.70\pm18.76$&$50$&$16.21\pm17.05$\\ \midrule
Constant demands&411 (900)&$3$&$27.89\pm19.75$&$115$&$16.76\pm20.63$\\Sparse demands&297 (900)&$0$&$30.38\pm18.11$&$115$&$17.26\pm21.32$\\ \midrule
Factor $\rho = 0.0$&157 (360)&$3$&$45.80\pm23.34$&$18$&$39.41\pm27.13$\\Factor $\rho = 0.25$&162 (360)&$0$&$32.03\pm15.77$&$46$&$16.82\pm15.35$\\Factor $\rho = 0.5$&131 (360)&$0$&$21.72\pm12.53$&$59$&$7.23\pm10.38$\\Factor $\rho = 0.75$&128 (360)&$0$&$20.60\pm12.61$&$53$&$8.19\pm10.90$\\Factor $\rho = 1.0$&130 (360)&$0$&$20.18\pm12.62$&$54$&$8.50\pm11.34$\\ \bottomrule

    \end{tabular}
    \label{tab:optgaps-optimistic}
\end{table}
%

\subsubsection{Pessimistic Variant of the \ourproblem{}} \label{apx:supplementary-pessimistic}

Table~\ref{tab:summary-pessimistic} and Figures~\ref{fig:runtimes-pessimistic}--~\ref{fig:objectives-pessimistic} have the same meaning as Table~\ref{tab:summary-optimistic} and Figures~\ref{fig:runtimes-optimistic}--~\ref{fig:objectives-optimistic}, but for the pessimistic variant of the \ourproblem{}.
Similarly, Tables~\ref{tab:soltimes-pessimistic} and \ref{tab:optgaps-pessimistic} have the same meaning as Tables~\ref{tab:soltimes-optimistic} and \ref{tab:optgaps-optimistic}, but for the pessimistic variant of the \ourproblem{}.
We draw similar conclusions about the computational performance of our solution methods when faced with the pessimistic variant of the \ourproblem{}.

\begin{table}[!ht]
    \centering
    \caption{Percentage of the benchmark solved to optimality by each exact method, grouped by dimensional instance attributes (\ie, number of locations $|\mathcal{I}|$, customers $|\mathcal{J}|$, and periods $|\mathcal{T}|$). [Pessimistic variant of the \ourproblem{}]}
    \begin{tabular}{ccccc}
    \toprule
        \multirow{2}{*}{$|\mathcal{T}|$}& \multicolumn{2}{c}{Montreal ($|\mathcal{I}| = 20, |\mathcal{J}| = 40$)} & \multicolumn{2}{c}{Quebec ($|\mathcal{I}| = 39, |\mathcal{J}| = 78$)}  \\ \cmidrule(lr){2-5}
        & Tailored B\&C& Tightened B\&C& Tailored B\&C & Tightened B\&C\\ \midrule
        3&100.00\%&100.00\%&97.33\%&96.67\%\\ \midrule
5&58.67\%&72.33\%&48.00\%&62.67\%\\ \midrule
7&7.00\%&38.00\%&6.33\%&26.33\%\\ \bottomrule

    \end{tabular}
    \label{tab:summary-pessimistic}
\end{table}
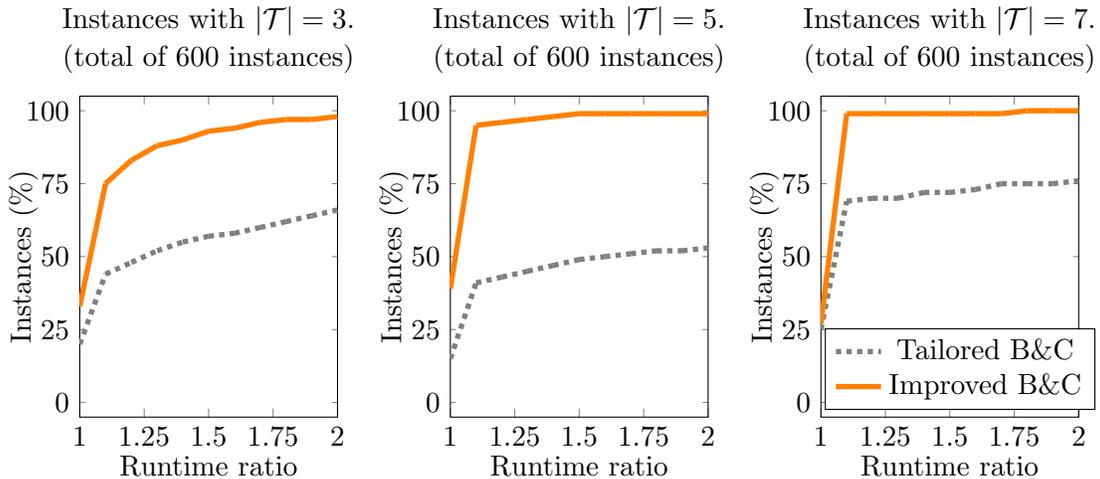
\begin{figure}[!ht]
     \centering
     \caption{Runtime ratios of both exact methods for exact methods, where the $y$ axis presents the percentage of instances with a ratio smaller than or equal to the reference value on the $x$ axis. [Pessimistic variant of the \ourproblem{}]}
    \begin{tikzpicture}
\begin{axis}[name = 3periods, xmin = 1, xmax = 2, xscale = 0.5, ymin = -5, ymax = 105, yscale = 0.75, xtick = {1,1.25,1.5,1.75,2}, ytick = {0, 25, 50, 75, 100}, xlabel = {Runtime ratio}, ylabel = {Instances (\%)}, legend pos = south west]
\addplot[line width = 2pt, color = gray, style = dotted] coordinates {(1,20) (1.1,44)};
\addplot[line width = 2pt, color = gray, style = dotted] coordinates {(1.1,44) (1.2,48)};
\addplot[line width = 2pt, color = gray, style = dotted] coordinates {(1.2,48) (1.3,52)};
\addplot[line width = 2pt, color = gray, style = dotted] coordinates {(1.3,52) (1.4,55)};
\addplot[line width = 2pt, color = gray, style = dotted] coordinates {(1.4,55) (1.5,57)};
\addplot[line width = 2pt, color = gray, style = dotted] coordinates {(1.5,57) (1.6,58)};
\addplot[line width = 2pt, color = gray, style = dotted] coordinates {(1.6,58) (1.7,60)};
\addplot[line width = 2pt, color = gray, style = dotted] coordinates {(1.7,60) (1.8,62)};
\addplot[line width = 2pt, color = gray, style = dotted] coordinates {(1.8,62) (1.9,64)};
\addplot[line width = 2pt, color = gray, style = dotted] coordinates {(1.9,64) (2.0,66)};
\addplot[line width = 2pt, color = orange, style = solid] coordinates {(1,33) (1.1,75)};
\addplot[line width = 2pt, color = orange, style = solid] coordinates {(1.1,75) (1.2,83)};
\addplot[line width = 2pt, color = orange, style = solid] coordinates {(1.2,83) (1.3,88)};
\addplot[line width = 2pt, color = orange, style = solid] coordinates {(1.3,88) (1.4,90)};
\addplot[line width = 2pt, color = orange, style = solid] coordinates {(1.4,90) (1.5,93)};
\addplot[line width = 2pt, color = orange, style = solid] coordinates {(1.5,93) (1.6,94)};
\addplot[line width = 2pt, color = orange, style = solid] coordinates {(1.6,94) (1.7,96)};
\addplot[line width = 2pt, color = orange, style = solid] coordinates {(1.7,96) (1.8,97)};
\addplot[line width = 2pt, color = orange, style = solid] coordinates {(1.8,97) (1.9,97)};
\addplot[line width = 2pt, color = orange, style = solid] coordinates {(1.9,97) (2.0,98)};

\end{axis}
\begin{axis}[name = 5periods, at={(3periods.east)}, anchor=west, xshift = 1.5cm, xmin = 1, xmax = 2, xscale = 0.5, ymin = -5, ymax = 105, yscale = 0.75, xtick = {1,1.25,1.5,1.75,2}, ytick = {0, 25, 50, 75, 100}, xlabel = {Runtime ratio}, ylabel = {Instances (\%)}, legend pos = south west]
\addplot[line width = 2pt, color = gray, style = dotted] coordinates {(1,15) (1.1,41)};
\addplot[line width = 2pt, color = gray, style = dotted] coordinates {(1.1,41) (1.2,43)};
\addplot[line width = 2pt, color = gray, style = dotted] coordinates {(1.2,43) (1.3,45)};
\addplot[line width = 2pt, color = gray, style = dotted] coordinates {(1.3,45) (1.4,47)};
\addplot[line width = 2pt, color = gray, style = dotted] coordinates {(1.4,47) (1.5,49)};
\addplot[line width = 2pt, color = gray, style = dotted] coordinates {(1.5,49) (1.6,50)};
\addplot[line width = 2pt, color = gray, style = dotted] coordinates {(1.6,50) (1.7,51)};
\addplot[line width = 2pt, color = gray, style = dotted] coordinates {(1.7,51) (1.8,52)};
\addplot[line width = 2pt, color = gray, style = dotted] coordinates {(1.8,52) (1.9,52)};
\addplot[line width = 2pt, color = gray, style = dotted] coordinates {(1.9,52) (2.0,53)};
\addplot[line width = 2pt, color = orange, style = solid] coordinates {(1,39) (1.1,95)};
\addplot[line width = 2pt, color = orange, style = solid] coordinates {(1.1,95) (1.2,96)};
\addplot[line width = 2pt, color = orange, style = solid] coordinates {(1.2,96) (1.3,97)};
\addplot[line width = 2pt, color = orange, style = solid] coordinates {(1.3,97) (1.4,98)};
\addplot[line width = 2pt, color = orange, style = solid] coordinates {(1.4,98) (1.5,99)};
\addplot[line width = 2pt, color = orange, style = solid] coordinates {(1.5,99) (1.6,99)};
\addplot[line width = 2pt, color = orange, style = solid] coordinates {(1.6,99) (1.7,99)};
\addplot[line width = 2pt, color = orange, style = solid] coordinates {(1.7,99) (1.8,99)};
\addplot[line width = 2pt, color = orange, style = solid] coordinates {(1.8,99) (1.9,99)};
\addplot[line width = 2pt, color = orange, style = solid] coordinates {(1.9,99) (2.0,99)};

\end{axis}
\begin{axis}[name = 7periods, at={(5periods.east)}, anchor=west, xshift = 1.5cm, xmin = 1, xmax = 2, xscale = 0.5, ymin = -5, ymax = 105, yscale = 0.75, xtick = {1,1.25,1.5,1.75,2}, ytick = {0, 25, 50, 75, 100}, xlabel = {Runtime ratio}, ylabel = {Instances (\%)}, legend pos = south west]
\addlegendimage{line width = 2, color = gray, style = dotted}
\addlegendentry{Tailored B\&C}
\addlegendimage{line width = 2, color = orange, style = solid}
\addlegendentry{Improved B\&C}
\addplot[line width = 2pt, color = gray, style = dotted] coordinates {(1,25) (1.1,69)};
\addplot[line width = 2pt, color = gray, style = dotted] coordinates {(1.1,69) (1.2,70)};
\addplot[line width = 2pt, color = gray, style = dotted] coordinates {(1.2,70) (1.3,70)};
\addplot[line width = 2pt, color = gray, style = dotted] coordinates {(1.3,70) (1.4,72)};
\addplot[line width = 2pt, color = gray, style = dotted] coordinates {(1.4,72) (1.5,72)};
\addplot[line width = 2pt, color = gray, style = dotted] coordinates {(1.5,72) (1.6,73)};
\addplot[line width = 2pt, color = gray, style = dotted] coordinates {(1.6,73) (1.7,75)};
\addplot[line width = 2pt, color = gray, style = dotted] coordinates {(1.7,75) (1.8,75)};
\addplot[line width = 2pt, color = gray, style = dotted] coordinates {(1.8,75) (1.9,75)};
\addplot[line width = 2pt, color = gray, style = dotted] coordinates {(1.9,75) (2.0,76)};
\addplot[line width = 2pt, color = orange, style = solid] coordinates {(1,27) (1.1,99)};
\addplot[line width = 2pt, color = orange, style = solid] coordinates {(1.1,99) (1.2,99)};
\addplot[line width = 2pt, color = orange, style = solid] coordinates {(1.2,99) (1.3,99)};
\addplot[line width = 2pt, color = orange, style = solid] coordinates {(1.3,99) (1.4,99)};
\addplot[line width = 2pt, color = orange, style = solid] coordinates {(1.4,99) (1.5,99)};
\addplot[line width = 2pt, color = orange, style = solid] coordinates {(1.5,99) (1.6,99)};
\addplot[line width = 2pt, color = orange, style = solid] coordinates {(1.6,99) (1.7,99)};
\addplot[line width = 2pt, color = orange, style = solid] coordinates {(1.7,99) (1.8,100)};
\addplot[line width = 2pt, color = orange, style = solid] coordinates {(1.8,100) (1.9,100)};
\addplot[line width = 2pt, color = orange, style = solid] coordinates {(1.9,100) (2.0,100)};

\end{axis}
\node[align=center, yshift=1cm] at (3periods.north) {Instances with $|\mathcal{T}| = 3$.};
\node[align=center, yshift=0.5cm] at (3periods.north) {(total of $600$ instances)};
\node[align=center, yshift=1cm] at (5periods.north) {Instances with $|\mathcal{T}| = 5$.};
\node[align=center, yshift=0.5cm] at (5periods.north) {(total of $600$ instances)};
\node[align=center, yshift=1cm] at (7periods.north) {Instances with $|\mathcal{T}| = 7$.};
\node[align=center, yshift=0.5cm] at (7periods.north) {(total of $600$ instances)};
\end{tikzpicture}
    \label{fig:runtimes-pessimistic}
\end{figure}
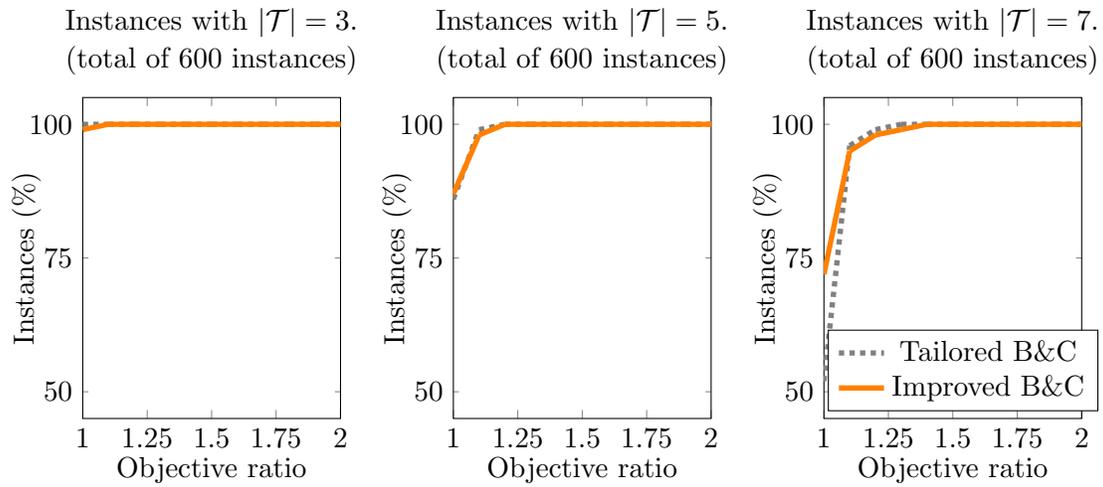
\begin{figure}[!ht]
     \centering
     \caption{Objective ratios of both exact methods for exact methods, where the $y$ axis presents the percentage of instances with a ratio smaller than or equal to the reference value on the $x$ axis. [Pessimistic variant of the \ourproblem{}]}    
    \begin{tikzpicture}
\begin{axis}[name = 3periods, xmin = 1, xmax = 2, xscale = 0.5, ymin = 45, ymax = 105, yscale = 0.75, xtick = {1,1.25,1.5,1.75,2}, ytick = {50, 75, 100}, xlabel = {Objective ratio}, ylabel = {Instances (\%)}, legend pos = south west]
\addplot[line width = 2pt, color = gray, style = dotted] coordinates {(1,100) (1.1,100)};
\addplot[line width = 2pt, color = gray, style = dotted] coordinates {(1.1,100) (1.2,100)};
\addplot[line width = 2pt, color = gray, style = dotted] coordinates {(1.2,100) (1.3,100)};
\addplot[line width = 2pt, color = gray, style = dotted] coordinates {(1.3,100) (1.4,100)};
\addplot[line width = 2pt, color = gray, style = dotted] coordinates {(1.4,100) (1.5,100)};
\addplot[line width = 2pt, color = gray, style = dotted] coordinates {(1.5,100) (1.6,100)};
\addplot[line width = 2pt, color = gray, style = dotted] coordinates {(1.6,100) (1.7,100)};
\addplot[line width = 2pt, color = gray, style = dotted] coordinates {(1.7,100) (1.8,100)};
\addplot[line width = 2pt, color = gray, style = dotted] coordinates {(1.8,100) (1.9,100)};
\addplot[line width = 2pt, color = gray, style = dotted] coordinates {(1.9,100) (2.0,100)};
\addplot[line width = 2pt, color = orange, style = solid] coordinates {(1,99) (1.1,100)};
\addplot[line width = 2pt, color = orange, style = solid] coordinates {(1.1,100) (1.2,100)};
\addplot[line width = 2pt, color = orange, style = solid] coordinates {(1.2,100) (1.3,100)};
\addplot[line width = 2pt, color = orange, style = solid] coordinates {(1.3,100) (1.4,100)};
\addplot[line width = 2pt, color = orange, style = solid] coordinates {(1.4,100) (1.5,100)};
\addplot[line width = 2pt, color = orange, style = solid] coordinates {(1.5,100) (1.6,100)};
\addplot[line width = 2pt, color = orange, style = solid] coordinates {(1.6,100) (1.7,100)};
\addplot[line width = 2pt, color = orange, style = solid] coordinates {(1.7,100) (1.8,100)};
\addplot[line width = 2pt, color = orange, style = solid] coordinates {(1.8,100) (1.9,100)};
\addplot[line width = 2pt, color = orange, style = solid] coordinates {(1.9,100) (2.0,100)};

\end{axis}
\begin{axis}[name = 5periods, at={(3periods.east)}, anchor=west, xshift = 1.5cm, xmin = 1, xmax = 2, xscale = 0.5, ymin = 45, ymax = 105, yscale = 0.75, xtick = {1,1.25,1.5,1.75,2}, ytick = {50, 75, 100}, xlabel = {Objective ratio}, ylabel = {Instances (\%)}, legend pos = south west]
\addplot[line width = 2pt, color = gray, style = dotted] coordinates {(1,86) (1.1,99)};
\addplot[line width = 2pt, color = gray, style = dotted] coordinates {(1.1,99) (1.2,100)};
\addplot[line width = 2pt, color = gray, style = dotted] coordinates {(1.2,100) (1.3,100)};
\addplot[line width = 2pt, color = gray, style = dotted] coordinates {(1.3,100) (1.4,100)};
\addplot[line width = 2pt, color = gray, style = dotted] coordinates {(1.4,100) (1.5,100)};
\addplot[line width = 2pt, color = gray, style = dotted] coordinates {(1.5,100) (1.6,100)};
\addplot[line width = 2pt, color = gray, style = dotted] coordinates {(1.6,100) (1.7,100)};
\addplot[line width = 2pt, color = gray, style = dotted] coordinates {(1.7,100) (1.8,100)};
\addplot[line width = 2pt, color = gray, style = dotted] coordinates {(1.8,100) (1.9,100)};
\addplot[line width = 2pt, color = gray, style = dotted] coordinates {(1.9,100) (2.0,100)};
\addplot[line width = 2pt, color = orange, style = solid] coordinates {(1,87) (1.1,98)};
\addplot[line width = 2pt, color = orange, style = solid] coordinates {(1.1,98) (1.2,100)};
\addplot[line width = 2pt, color = orange, style = solid] coordinates {(1.2,100) (1.3,100)};
\addplot[line width = 2pt, color = orange, style = solid] coordinates {(1.3,100) (1.4,100)};
\addplot[line width = 2pt, color = orange, style = solid] coordinates {(1.4,100) (1.5,100)};
\addplot[line width = 2pt, color = orange, style = solid] coordinates {(1.5,100) (1.6,100)};
\addplot[line width = 2pt, color = orange, style = solid] coordinates {(1.6,100) (1.7,100)};
\addplot[line width = 2pt, color = orange, style = solid] coordinates {(1.7,100) (1.8,100)};
\addplot[line width = 2pt, color = orange, style = solid] coordinates {(1.8,100) (1.9,100)};
\addplot[line width = 2pt, color = orange, style = solid] coordinates {(1.9,100) (2.0,100)};

\end{axis}
\begin{axis}[name = 7periods, at={(5periods.east)}, anchor=west, xshift = 1.5cm, xmin = 1, xmax = 2, xscale = 0.5, ymin = 45, ymax = 105, yscale = 0.75, xtick = {1,1.25,1.5,1.75,2}, ytick = {50, 75, 100}, xlabel = {Objective ratio}, ylabel = {Instances (\%)}, legend pos = south west]
\addlegendimage{line width = 2, color = gray, style = dotted}
\addlegendentry{Tailored B\&C}
\addlegendimage{line width = 2, color = orange, style = solid}
\addlegendentry{Improved B\&C}
\addplot[line width = 2pt, color = gray, style = dotted] coordinates {(1,52) (1.1,96)};
\addplot[line width = 2pt, color = gray, style = dotted] coordinates {(1.1,96) (1.2,99)};
\addplot[line width = 2pt, color = gray, style = dotted] coordinates {(1.2,99) (1.3,100)};
\addplot[line width = 2pt, color = gray, style = dotted] coordinates {(1.3,100) (1.4,100)};
\addplot[line width = 2pt, color = gray, style = dotted] coordinates {(1.4,100) (1.5,100)};
\addplot[line width = 2pt, color = gray, style = dotted] coordinates {(1.5,100) (1.6,100)};
\addplot[line width = 2pt, color = gray, style = dotted] coordinates {(1.6,100) (1.7,100)};
\addplot[line width = 2pt, color = gray, style = dotted] coordinates {(1.7,100) (1.8,100)};
\addplot[line width = 2pt, color = gray, style = dotted] coordinates {(1.8,100) (1.9,100)};
\addplot[line width = 2pt, color = gray, style = dotted] coordinates {(1.9,100) (2.0,100)};
\addplot[line width = 2pt, color = orange, style = solid] coordinates {(1,72) (1.1,95)};
\addplot[line width = 2pt, color = orange, style = solid] coordinates {(1.1,95) (1.2,98)};
\addplot[line width = 2pt, color = orange, style = solid] coordinates {(1.2,98) (1.3,99)};
\addplot[line width = 2pt, color = orange, style = solid] coordinates {(1.3,99) (1.4,100)};
\addplot[line width = 2pt, color = orange, style = solid] coordinates {(1.4,100) (1.5,100)};
\addplot[line width = 2pt, color = orange, style = solid] coordinates {(1.5,100) (1.6,100)};
\addplot[line width = 2pt, color = orange, style = solid] coordinates {(1.6,100) (1.7,100)};
\addplot[line width = 2pt, color = orange, style = solid] coordinates {(1.7,100) (1.8,100)};
\addplot[line width = 2pt, color = orange, style = solid] coordinates {(1.8,100) (1.9,100)};
\addplot[line width = 2pt, color = orange, style = solid] coordinates {(1.9,100) (2.0,100)};

\end{axis}
\node[align=center, yshift=1cm] at (3periods.north) {Instances with $|\mathcal{T}| = 3$.};
\node[align=center, yshift=0.5cm] at (3periods.north) {(total of $600$ instances)};
\node[align=center, yshift=1cm] at (5periods.north) {Instances with $|\mathcal{T}| = 5$.};
\node[align=center, yshift=0.5cm] at (5periods.north) {(total of $600$ instances)};
\node[align=center, yshift=1cm] at (7periods.north) {Instances with $|\mathcal{T}| = 7$.};
\node[align=center, yshift=0.5cm] at (7periods.north) {(total of $600$ instances)};
\end{tikzpicture}
    \label{fig:objectives-pessimistic}
\end{figure}

\begin{table}[!ht]
    \centering
    \caption{Average solution times of Tailored B\&C and Tightened B\&C, as well as their standard deviations, for instances solved to optimality by both exact methods. [Pessimistic variant of the \ourproblem{}]}
    \begin{tabular}{crrr}
        \toprule
        \multirow{2}*{Instance attributes} &  \# instances & Tailored B\&C & Tightened B\&C \\ \cmidrule{3-4}
        & considered & time (min) & time (min)  \\ \midrule
        Benchmark&948 (1800)&$111.93\pm265.98$&$47.01\pm153.51$\\ \midrule
$|\mathcal{T}|$ = 3 periods&590 (600)&$35.07\pm154.32$&$34.77\pm154.84$\\$|\mathcal{T}|$ = 5 periods&318 (600)&$226.21\pm346.89$&$59.47\pm140.73$\\$|\mathcal{T}|$ = 7 periods&40 (600)&$337.12\pm363.60$&$128.63\pm196.91$\\ \midrule
Montreal region&497 (900)&$120.74\pm278.13$&$30.12\pm101.73$\\Quebec province&451 (900)&$102.22\pm251.85$&$65.64\pm193.70$\\ \midrule
Travel $\leq$ 15 minutes&377 (600)&$45.96\pm163.02$&$21.36\pm91.23$\\Travel $\leq$ 30 minutes&313 (600)&$169.75\pm308.38$&$62.18\pm165.34$\\Travel $\leq$ 45 minutes&258 (600)&$138.18\pm308.54$&$66.11\pm199.65$\\ \midrule
Identical rewards&482 (900)&$86.09\pm230.45$&$8.44\pm20.51$\\Inverse rewards&466 (900)&$138.66\pm296.19$&$86.91\pm210.76$\\ \midrule
Constant demands&392 (900)&$145.80\pm318.36$&$70.46\pm209.69$\\Sparse demands&556 (900)&$88.05\pm218.91$&$30.49\pm92.57$\\ \midrule
Factor $\rho = 0.0$&175 (360)&$90.63\pm243.74$&$79.73\pm221.32$\\Factor $\rho = 0.25$&169 (360)&$68.48\pm199.04$&$29.09\pm108.50$\\Factor $\rho = 0.5$&201 (360)&$152.79\pm314.29$&$41.95\pm143.20$\\Factor $\rho = 0.75$&202 (360)&$120.99\pm274.11$&$43.10\pm137.18$\\Factor $\rho = 1.0$&201 (360)&$117.04\pm268.21$&$42.61\pm135.61$\\ \bottomrule

    \end{tabular}
    \label{tab:soltimes-pessimistic}
\end{table}
\begin{table}[!ht]
    \centering
    \caption{Number of instances solved to optimality, average optimality gaps and their standard deviations, of Tailored B\&C and Tightened B\&C for instances not solved to optimality by both exact methods. [Pessimistic variant of the \ourproblem{}]}
    \begin{tabular}{crrrrr}
        \toprule
        \multirow{2}*{Instance attributes} &  \# instances & \multicolumn{2}{c}{Tailored B\&C} & \multicolumn{2}{c}{Tightened B\&C} \\ \cmidrule{3-6}
        & considered & \# opt. & opt. gap (\%) & \# opt. & opt. gap (\%)  \\
        \midrule
        Benchmark&852 (1800)&$4$&$26.83\pm18.73$&$240$&$16.65\pm19.48$\\ \midrule
$|\mathcal{T}|$ = 3 periods&10 (600)&$2$&$1.79\pm2.61$&$0$&$2.25\pm2.46$\\$|\mathcal{T}|$ = 5 periods&282 (600)&$2$&$20.26\pm17.41$&$87$&$13.69\pm18.03$\\$|\mathcal{T}|$ = 7 periods&560 (600)&$0$&$30.58\pm18.29$&$153$&$18.39\pm20.07$\\ \midrule
Montreal region&403 (900)&$0$&$32.53\pm19.07$&$134$&$20.12\pm22.49$\\Quebec province&449 (900)&$4$&$21.71\pm16.86$&$106$&$13.53\pm15.69$\\ \midrule
Travel $\leq$ 15 minutes&223 (600)&$2$&$18.36\pm12.49$&$95$&$6.57\pm8.31$\\Travel $\leq$ 30 minutes&287 (600)&$0$&$28.07\pm18.20$&$67$&$17.13\pm17.92$\\Travel $\leq$ 45 minutes&342 (600)&$2$&$31.30\pm20.68$&$78$&$22.82\pm22.97$\\ \midrule
Identical rewards&418 (900)&$0$&$33.98\pm18.17$&$183$&$17.38\pm23.22$\\Inverse rewards&434 (900)&$4$&$19.94\pm16.56$&$57$&$15.95\pm15.01$\\ \midrule
Constant demands&508 (900)&$2$&$25.72\pm18.78$&$108$&$17.15\pm18.73$\\Sparse demands&344 (900)&$2$&$28.45\pm18.56$&$132$&$15.91\pm20.53$\\ \midrule
Factor $\rho = 0.0$&185 (360)&$1$&$42.22\pm23.98$&$15$&$37.13\pm25.78$\\Factor $\rho = 0.25$&191 (360)&$2$&$29.40\pm16.58$&$46$&$16.12\pm14.65$\\Factor $\rho = 0.5$&159 (360)&$1$&$20.47\pm12.32$&$60$&$8.22\pm10.01$\\Factor $\rho = 0.75$&158 (360)&$0$&$19.58\pm12.16$&$59$&$9.00\pm10.51$\\Factor $\rho = 1.0$&159 (360)&$0$&$19.38\pm12.47$&$60$&$9.48\pm11.33$\\ \bottomrule

    \end{tabular}
    \label{tab:optgaps-pessimistic}
\end{table}

\subsection{Supplementary Managerial Insights} \label{apx:managerial}

We here  present supplementary managerial results for both variants of the \ourproblem{}.

\subsubsection{Optimistic Variant of the \ourproblem{}} \label{apx:managerial-optimistic}

Figure~\ref{fig:services-optimistic} presents the average number of captures and the average percentage of captured demand in the optimal solution under competition (solid lines) and under cooperation (dotted lines), grouped by different instance attributes.
On average, we highlight that competition leads to customers being visited more often throughout the planning horizon, whereas cooperation induces a higher percentage of demand satisfied by the players at the end of the planning horizon, and that across different instance attributes. 
\begin{figure}[!ht]
     \centering
     \caption{Service quality under competition (solid lines) and cooperation (dotted lines) grouped by different instance attributes. [Optimistic  variant of the \ourproblem{}]}
     \include{graphs_optimistic/services}
     \label{fig:services-optimistic}
\end{figure}

\subsubsection{Pessimistic Variant of the \ourproblem{}} \label{apx:managerial-pessimistic}

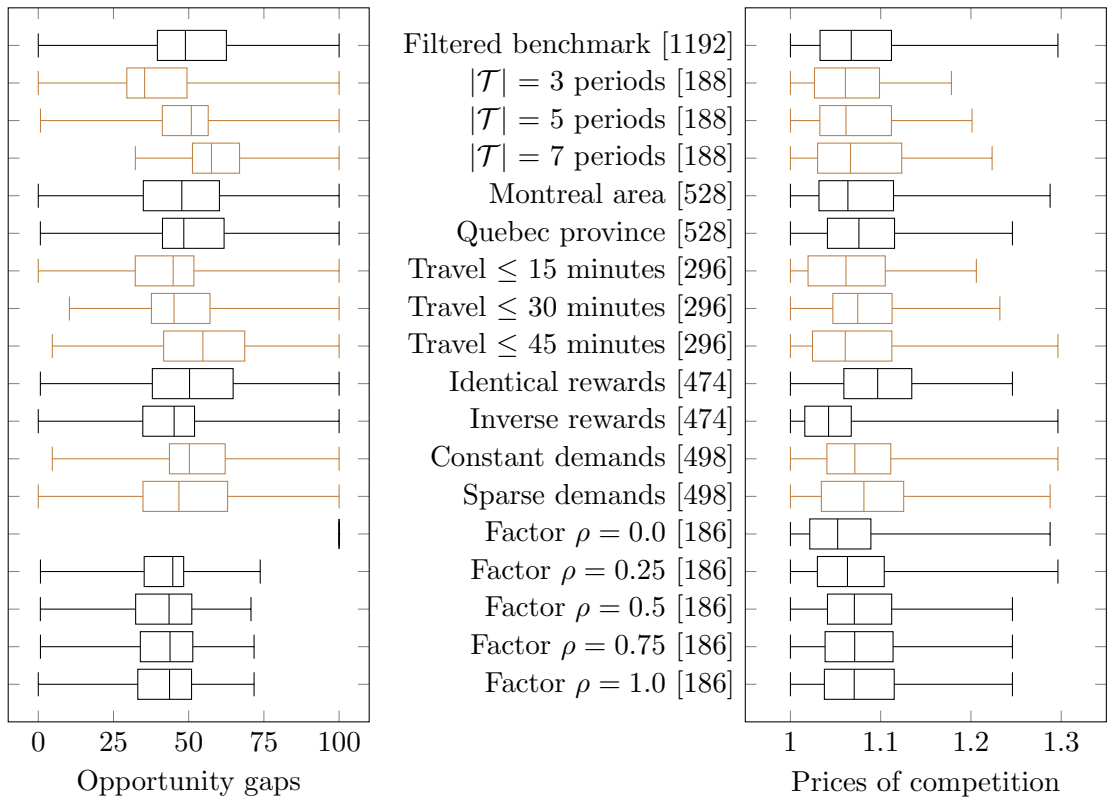
\begin{figure}
    \centering
    \caption{Opportunity gaps under the monopolistic assumption (left) and prices of competition (right) grouped by different instance attributes, with the number of instances between square brackets. [Pessimistic variant of the \ourproblem{}]}
    \begin{tikzpicture}
\begin{axis}[name = oppgaps, y=0.5cm, ymin=0, ymax=19, x = 0.04cm, xmin = -10, xmax = 110, xlabel = {Opportunity gaps}, xtick = {0,25,50,75,100},  ytick = {1,2,3,4,5,6,7,8,9,10,11,12,13,14,15,16,17,18}, yticklabels = {}]
\addplot+[
	boxplot prepared={
	lower whisker=0.0,
	lower quartile=33.074025,
	median=43.635099999999994,
	upper quartile=50.948750000000004,
	upper whisker=71.7325
}, color = black, style = solid
] coordinates {};
\addplot+[
	boxplot prepared={
	lower whisker=0.7236,
	lower quartile=33.924525,
	median=43.79625,
	upper quartile=51.3663,
	upper whisker=71.7325
}, color = black, style = solid
] coordinates {};
\addplot+[
	boxplot prepared={
	lower whisker=0.6944,
	lower quartile=32.36455,
	median=43.4744,
	upper quartile=51.034775,
	upper whisker=70.7087
}, color = black, style = solid
] coordinates {};
\addplot+[
	boxplot prepared={
	lower whisker=0.7236,
	lower quartile=35.2031,
	median=44.70415,
	upper quartile=48.3294,
	upper whisker=73.7762
}, color = black, style = solid
] coordinates {};
\addplot+[
	boxplot prepared={
	lower whisker=100.0,
	lower quartile=100.0,
	median=100.0,
	upper quartile=100.0,
	upper whisker=100.0
}, color = black, style = solid
] coordinates {};
\addplot+[
	boxplot prepared={
	lower whisker=0.0,
	lower quartile=34.810675,
	median=46.7532,
	upper quartile=62.93965,
	upper whisker=100.0
}, color = brown, style = solid
] coordinates {};
\addplot+[
	boxplot prepared={
	lower whisker=4.668,
	lower quartile=43.5986,
	median=50.23365,
	upper quartile=62.126025,
	upper whisker=100.0
}, color = brown, style = solid
] coordinates {};
\addplot+[
	boxplot prepared={
	lower whisker=0.0,
	lower quartile=34.7444,
	median=45.1519,
	upper quartile=51.938725,
	upper whisker=100.0
}, color = black, style = solid
] coordinates {};
\addplot+[
	boxplot prepared={
	lower whisker=0.6944,
	lower quartile=37.91715,
	median=50.27915,
	upper quartile=64.770325,
	upper whisker=100.0
}, color = black, style = solid
] coordinates {};
\addplot+[
	boxplot prepared={
	lower whisker=4.668,
	lower quartile=41.669525,
	median=54.7168,
	upper quartile=68.63352499999999,
	upper whisker=100.0
}, color = brown, style = solid
] coordinates {};
\addplot+[
	boxplot prepared={
	lower whisker=10.3519,
	lower quartile=37.609475,
	median=45.119749999999996,
	upper quartile=57.0866,
	upper whisker=100.0
}, color = brown, style = solid
] coordinates {};
\addplot+[
	boxplot prepared={
	lower whisker=0.0,
	lower quartile=32.239225,
	median=44.831149999999994,
	upper quartile=51.7316,
	upper whisker=100.0
}, color = brown, style = solid
] coordinates {};
\addplot+[
	boxplot prepared={
	lower whisker=0.6944,
	lower quartile=41.287525,
	median=48.3699,
	upper quartile=61.763799999999996,
	upper whisker=100.0
}, color = black, style = solid
] coordinates {};
\addplot+[
	boxplot prepared={
	lower whisker=0.0,
	lower quartile=34.9169,
	median=47.7071,
	upper quartile=60.21505,
	upper whisker=100.0
}, color = black, style = solid
] coordinates {};
\addplot+[
	boxplot prepared={
	lower whisker=32.2687,
	lower quartile=51.2577,
	median=57.5758,
	upper quartile=66.906475,
	upper whisker=100.0
}, color = brown, style = solid
] coordinates {};
\addplot+[
	boxplot prepared={
	lower whisker=0.7236,
	lower quartile=41.225274999999996,
	median=50.8771,
	upper quartile=56.485325,
	upper whisker=100.0
}, color = brown, style = solid
] coordinates {};
\addplot+[
	boxplot prepared={
	lower whisker=0.0,
	lower quartile=29.434525,
	median=35.35185,
	upper quartile=49.437375,
	upper whisker=100.0
}, color = brown, style = solid
] coordinates {};
\addplot+[
	boxplot prepared={
	lower whisker=0.0,
	lower quartile=39.6,
	median=48.90745,
	upper quartile=62.5328,
	upper whisker=100.0
}, color = black, style = solid
] coordinates {};
\end{axis}

\begin{axis}[at={(oppgaps.east)}, anchor=west, xshift = 5cm,y=0.5cm, ymin=0, ymax=19, x=12cm, xmin = 0.95, xmax = 1.35, xlabel = {Prices of competition}, xtick = {1,1.1,1.2,1.3},  ytick = {1,2,3,4,5,6,7,8,9,10,11,12,13,14,15,16,17,18}, yticklabels = {Factor $\rho = 1.0$ [186],Factor $\rho = 0.75$ [186],Factor $\rho = 0.5$ [186],Factor $\rho = 0.25$ [186],Factor $\rho = 0.0$ [186],Sparse demands [498],Constant demands [498],Inverse rewards [474],Identical rewards [474],Travel $\leq$ 45 minutes [296],Travel $\leq$ 30 minutes [296],Travel $\leq$ 15 minutes [296],Quebec province [528],Montreal area [528],$|\mathcal{T}|$ = 7 periods [188],$|\mathcal{T}|$ = 5 periods [188],$|\mathcal{T}|$ = 3 periods [188],Filtered benchmark [1192]}]
\addplot+[
	boxplot prepared={
	lower whisker=1.0,
	lower quartile=1.0376174767886401,
	median=1.0706708076265516,
	upper quartile=1.1150047141281532,
	upper whisker=1.2459677419354838
}, color = black, style = solid
] coordinates {};
\addplot+[
	boxplot prepared={
	lower whisker=1.0,
	lower quartile=1.0384144321727502,
	median=1.0711956625448438,
	upper quartile=1.11371632996633,
	upper whisker=1.2459677419354838
}, color = black, style = solid
] coordinates {};
\addplot+[
	boxplot prepared={
	lower whisker=1.0,
	lower quartile=1.041050347237297,
	median=1.0708437428986683,
	upper quartile=1.112173698151596,
	upper whisker=1.2459677419354838
}, color = black, style = solid
] coordinates {};
\addplot+[
	boxplot prepared={
	lower whisker=1.0,
	lower quartile=1.0298709244703128,
	median=1.0631882303400357,
	upper quartile=1.103944158868845,
	upper whisker=1.2963619021665076
}, color = black, style = solid
] coordinates {};
\addplot+[
	boxplot prepared={
	lower whisker=1.0,
	lower quartile=1.0214066485668063,
	median=1.052323507800392,
	upper quartile=1.0891467251522235,
	upper whisker=1.2877008652657602
}, color = black, style = solid
] coordinates {};
\addplot+[
	boxplot prepared={
	lower whisker=1.0,
	lower quartile=1.0342645932536378,
	median=1.0814168807120244,
	upper quartile=1.12547283823796,
	upper whisker=1.2877008652657602
}, color = brown, style = solid
] coordinates {};
\addplot+[
	boxplot prepared={
	lower whisker=1.0,
	lower quartile=1.040390802832118,
	median=1.0713685978169605,
	upper quartile=1.111174841357089,
	upper whisker=1.2963619021665076
}, color = brown, style = solid
] coordinates {};
\addplot+[
	boxplot prepared={
	lower whisker=1.0,
	lower quartile=1.0159271108638197,
	median=1.0422558922558922,
	upper quartile=1.0674006908462867,
	upper whisker=1.2963619021665076
}, color = black, style = solid
] coordinates {};
\addplot+[
	boxplot prepared={
	lower whisker=1.0,
	lower quartile=1.0592062001850513,
	median=1.096646942800789,
	upper quartile=1.1344463971880492,
	upper whisker=1.2459677419354838
}, color = black, style = solid
] coordinates {};
\addplot+[
	boxplot prepared={
	lower whisker=1.0,
	lower quartile=1.0244160760309278,
	median=1.0607912662125596,
	upper quartile=1.1123683900688295,
	upper whisker=1.2963619021665076
}, color = brown, style = solid
] coordinates {};
\addplot+[
	boxplot prepared={
	lower whisker=1.0,
	lower quartile=1.0470274443742542,
	median=1.074588403722262,
	upper quartile=1.1126237284030358,
	upper whisker=1.2320441988950277
}, color = brown, style = solid
] coordinates {};
\addplot+[
	boxplot prepared={
	lower whisker=1.0,
	lower quartile=1.019537699504678,
	median=1.0615449356934326,
	upper quartile=1.1050420168067228,
	upper whisker=1.2059846903270703
}, color = brown, style = solid
] coordinates {};
\addplot+[
	boxplot prepared={
	lower whisker=1.0,
	lower quartile=1.0410207161476197,
	median=1.0758671614984552,
	upper quartile=1.1153442203009927,
	upper whisker=1.2459677419354838
}, color = black, style = solid
] coordinates {};
\addplot+[
	boxplot prepared={
	lower whisker=1.0,
	lower quartile=1.031634446397188,
	median=1.063664596273292,
	upper quartile=1.1141975308641976,
	upper whisker=1.2877008652657602
}, color = black, style = solid
] coordinates {};
\addplot+[
	boxplot prepared={
	lower whisker=1.0,
	lower quartile=1.0301160808507184,
	median=1.0665476656392627,
	upper quartile=1.1233183856502242,
	upper whisker=1.2234392113910186
}, color = brown, style = solid
] coordinates {};
\addplot+[
	boxplot prepared={
	lower whisker=1.0,
	lower quartile=1.032640879648458,
	median=1.0615560072084884,
	upper quartile=1.1119023397761953,
	upper whisker=1.201219512195122
}, color = brown, style = solid
] coordinates {};
\addplot+[
	boxplot prepared={
	lower whisker=1.0,
	lower quartile=1.0266682330827068,
	median=1.0609137055837563,
	upper quartile=1.098672839546373,
	upper whisker=1.178343949044586
}, color = brown, style = solid
] coordinates {};
\addplot+[
	boxplot prepared={
	lower whisker=1.0,
	lower quartile=1.0328205195669211,
	median=1.0674006908462867,
	upper quartile=1.1119491794384984,
	upper whisker=1.2963619021665076
}, color = black, style = solid
] coordinates {};
\end{axis}
\end{tikzpicture}
    \label{fig:oppgaps-prices-pessimistic}
\end{figure}

\begin{figure}[!ht]
     \centering
     \caption{Resilience ratios averaged over our computational benchmark. [Pessimistic variant of the \ourproblem{}]}
     \label{fig:heatmaps-p}
     \begin{subfigure}{0.44\textwidth}
        \caption{Under an optimistic ground truth.}
        \begin{tikzpicture}
\begin{axis}[
    colormap/blackwhite,
    xlabel = {Ground truth $(optimistic, \rho_2)$},
    ylabel = {Assumption $(pessimistic, \rho_1)$},
    point meta min=0,
    point meta max=100,
    xmin=-0.125, xmax=1.125, xscale=0.8,
    xtick style={color=black},
    xtick={0,0.25,0.5,0.75,1},
    ymin=-0.125, ymax=1.125, yscale=0.8,
    ytick style={color=black},
    ytick={0,0.25,0.5,0.75,1},
    every node near coord/.append style={font=\footnotesize, color=black, anchor=center,
        /pgf/number format/.cd, fixed, fixed zerofill, precision=2
    }
    ]
    \addplot[%
        matrix plot*,
        point meta=explicit,
        mesh/cols=5,
        visualization depends on=\thisrow{Deviation}\as\std,,
        nodes near coords={\pgfmathprintnumber{\pgfplotspointmeta}},
    ] table[meta=Mean] {
    x y Mean Deviation
0.0 0.0 96.82 4.08
0.25 0.0 86.14 9.42
0.5 0.0 79.10 15.68
0.75 0.0 78.30 16.57
1.0 0.0 78.28 16.59
0.0 0.25 49.15 33.34
0.25 0.25 97.54 3.11
0.5 0.25 90.05 9.23
0.75 0.25 88.89 10.45
1.0 0.25 88.86 10.47
0.0 0.5 44.78 34.14
0.25 0.5 82.65 15.99
0.5 0.5 97.31 3.41
0.75 0.5 96.26 4.02
1.0 0.5 96.23 4.05
0.0 0.75 46.82 33.96
0.25 0.75 82.99 15.73
0.5 0.75 95.75 5.77
0.75 0.75 97.30 3.31
1.0 0.75 97.30 3.31
0.0 1.0 46.18 34.25
0.25 1.0 82.96 16.01
0.5 1.0 95.80 5.79
0.75 1.0 97.38 3.20
1.0 1.0 97.38 3.20
    };
    \end{axis}
    \end{tikzpicture}
        \label{fig:heatmap-po}
     \end{subfigure}
     \hfill
     \begin{subfigure}{0.55\textwidth}
        \caption{Under a pessimistic ground truth.}
        \begin{tikzpicture}
\begin{axis}[
    colormap/blackwhite,
    colorbar,
    xlabel = {Ground truth $(pessimistic, \rho_2)$},
    ylabel = {Assumption $(pessimistic, \rho_1)$},
    point meta min=0,
    point meta max=100,
    xmin=-0.125, xmax=1.125, xscale=0.8,
    xtick style={color=black},
    xtick={0,0.25,0.5,0.75,1},
    ymin=-0.125, ymax=1.125, yscale=0.8,
    ytick style={color=black},
    ytick={0,0.25,0.5,0.75,1},
    every node near coord/.append style={font=\footnotesize, color=black, anchor=center,
        /pgf/number format/.cd, fixed, fixed zerofill, precision=2
    }
    ]
    \addplot[%
        matrix plot*,
        point meta=explicit,
        mesh/cols=5,
        visualization depends on=\thisrow{Deviation}\as\std,,
        nodes near coords={\pgfmathprintnumber{\pgfplotspointmeta}},
    ] table[meta=Mean] {
    x y Mean Deviation
0.0 0.0 100.00 0.00
0.25 0.0 88.19 10.07
0.5 0.0 81.18 16.93
0.75 0.0 80.41 17.72
1.0 0.0 80.39 17.70
0.0 0.25 43.51 32.20
0.25 0.25 100.00 0.00
0.5 0.25 92.37 10.57
0.75 0.25 91.23 11.83
1.0 0.25 91.21 11.81
0.0 0.5 38.71 32.19
0.25 0.5 84.51 16.90
0.5 0.5 100.00 0.00
0.75 0.5 99.04 2.94
1.0 0.5 99.02 2.92
0.0 0.75 40.86 32.44
0.25 0.75 84.81 16.67
0.5 0.75 98.24 5.64
0.75 0.75 100.00 0.00
1.0 0.75 99.99 0.43
0.0 1.0 40.81 33.25
0.25 1.0 84.69 16.94
0.5 1.0 98.17 5.71
0.75 1.0 100.02 0.44
1.0 1.0 100.00 0.00
    };
    \end{axis}
    \end{tikzpicture}
        \label{fig:heatmap-pp}
    \end{subfigure}
\end{figure}
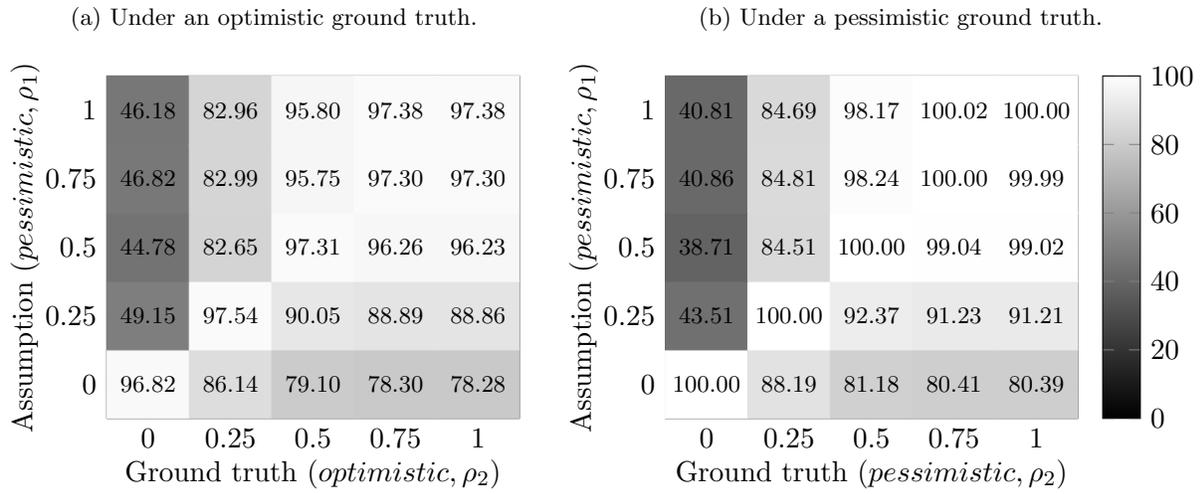

We consider here $1192$ instances solved to optimality by at least one of the exact methods.
Figures~\ref{fig:oppgaps-prices-pessimistic}, \ref{fig:heatmaps-p}, and \ref{fig:services-pessimistic} have the same meaning as Figures~\ref{fig:oppgaps-prices-optimistic}, \ref{fig:heatmaps-o}, and  \ref{fig:services-optimistic}, but for the pessimistic variant of the \ourproblem{}.
We observe trends similar to those reported for the optimistic variant of the \ourproblem{}.



\begin{figure}[!ht]
     \centering
     \caption{Service quality under competition (solid lines) and cooperation (dotted lines) grouped by different instance attributes. [Pessimistic variant of the \ourproblem{}]}
     \include{graphs_pessimistic/services}
     \label{fig:services-pessimistic}
\end{figure}

\end{document}